\documentclass[10pt,journal,comsoc]{IEEEtran}
\usepackage[T1]{fontenc}
\usepackage[linesnumbered,lined,ruled,commentsnumbered]{algorithm2e}
\usepackage{cite}
\usepackage{multirow}
\usepackage{microtype}
\usepackage{graphicx}
\usepackage{booktabs}
\usepackage{color}
\usepackage{hyperref}
\usepackage{amssymb,amsfonts, mathtools}
\usepackage{amsthm}
\usepackage{amsfonts}
\usepackage{threeparttable}
\usepackage{footnote}
\makesavenoteenv{tabular}
\usepackage{amsmath}
\usepackage{amsthm}
\usepackage{bm}
\usepackage[small]{subfigure}
\usepackage{caption}
\usepackage{tabularx}
\usepackage{multirow}
\allowdisplaybreaks[4]
 \usepackage[misc]{ifsym}

\let\oldnl\nl
\newcommand{\nonl}{\renewcommand{\nl}{\let\nl\oldnl}}

\newtheorem{assumption}{Assumption}
\newtheorem{lemma}{Lemma}
\newtheorem{remark}{Remark}
\newtheorem{theorem}{Theorem}
\newtheorem{corollary}{Corollary}[theorem]

\allowdisplaybreaks

\ifCLASSINFOpdf
\else
\fi

\ifCLASSINFOpdf
\else
\fi

\begin{document}


\title{On the Convergence of Quantized Parallel Restarted SGD for Central Server Free Distributed Training}

\author{Feijie~Wu, 
        Shiqi~He, 
        Yutong~Yang,
        Haozhao~Wang,
        Zhihao~Qu,
        \Letter Song~Guo,~\IEEEmembership{Fellow,~IEEE}
        and~Weihua~Zhuang,~\IEEEmembership{Fellow,~IEEE}
        \thanks{Song Guo is the corresponding author.}

\thanks{Feijie Wu and Song Guo are with the Department of Computing, The Hong Kong Polytechnic University, Hong Kong, China (e-mail: harli.wu@connect.polyu.hk, song.guo@polyu.edu.hk).} 
\thanks{Shiqi He is with the Department of Computer Science, The University of British Columbia, Vancouver, Canada (e-mail: shiqihe@cs.ubc.ca).}
\thanks{Yutong Yang is with National University of Singapore, Singapore (e-mail: e0575792@u.nus.edu).}
\thanks{Haozhao Wang is with School of Computer Science and Technology, Huazhong University of Science and Technology, Wuhan, China and the Department of Computing, The Hong Kong Polytechnic University, Hong Kong, China (e-mail: hz\_wang@hust.edu.cn).}
\thanks{Zhihao Qu is with the School of Computer and Information, Hohai University, Nanjing, China and Department of Computing, The Hong Kong Polytechnic University, Hong Kong, China (e-mail: quzhihao@hhu.edu.cn).}
\thanks{Weihua Zhuang is with the Department of Electrical and Computer Engineering, University of Waterloo, Waterloo, Canada (e-mail: wzhuang@uwaterloo.ca).}
}
\markboth{IEEE JOURNAL ON SELECTED AREAS IN COMMUNICATIONS, manuscript}%
{Wang \MakeLowercase{\textit{et al.}}: LOSP: Overlapping Computation and Quantized Communication in Parameter Server for Fast Distributed Training}

\maketitle

\begin{abstract}

Communication is a crucial phase in the context of distributed training. Because parameter server (PS) frequently experiences network congestion, recent studies have found that training paradigms without a centralized server outperform the traditional server-based paradigms in terms of communication efficiency. However, with the increasing growth of model sizes, these server-free paradigms are also confronted with substantial communication overhead that seriously deteriorates the performance of distributed training. In this paper, we focus on communication efficiency of two serverless paradigms, i.e., Ring All-Reduce (RAR) and gossip, by proposing the Quantized Parallel Restarted Stochastic Gradient Descent (QPRSGD), an algorithm that allows multiple local SGD updates before a global synchronization, in synergy with the quantization to significantly reduce the communication overhead. We establish the bound of accumulative errors according to the synchronization mode and the network topology, which is essential to ensure the convergence property. Under both aggregation paradigms, the algorithm achieves the linear speedup property with respect to the number of local updates as well as the number of workers. Remarkably, the proposed algorithm achieves a convergence rate $O(1/\sqrt{NK^2M})$ under the gossip paradigm and outperforms all existing compression methods, where $N$ is the times of global synchronizations, and $K$ is the number of local updates, while $M$ is the number of nodes. An empirical study on various machine learning models demonstrates that the communication overhead is reduced by 90\%, and the convergence speed is boosted by up to 18.6 times, in a low bandwidth network, in comparison with Parallel SGD. 

\end{abstract}


\IEEEpeerreviewmaketitle

\section{Introduction}


\IEEEPARstart{W}{ith} the growing data volume and the increasing concerns of data privacy, distributed machine learning is thriving with unprecedented prosperity. Its scalability makes it possible to expand the computational capability with the help of numerous nodes. Parallel Stochastic Gradient Descent (PSGD) is a general method for distributed machine learning, which allows clients to compute the gradients with their own data in parallel \cite{psgd_1_zinkevich,psgd_2_dekel,psgd_3_li}.  

\begin{table*}[]
\centering
\begin{tabularx}{0.98\textwidth}{cccccc}
\toprule
Reference  & PS or AR & Gossip & \begin{tabular}[c]{@{}c@{}}Bounded Gradients\end{tabular} & \begin{tabular}[c]{@{}c@{}}Periodical Averaging\end{tabular} & \begin{tabular}[c]{@{}c@{}}Communication Quantization\end{tabular} \\ 
\hline
\hline
DOUBLESQUEEZE \cite{tang2019doublesqueeze} &  $O\left(1/\sqrt{NM}\right)$   & NA    &    Yes  &  No &     All communication data      \\
Local SGD \cite{li2019communicationefficient}  & NA &  $O\left(1/\sqrt{NKM}\right)$  &    No  &   Yes  &   No quantization              \\
DCD-PSGD \cite{dcdpsgd}  &  NA  & $O\left(1/\sqrt{NM}\right)$ & No & No  & All communication data \\
QSparse-Local-SGD \cite{basu2019qsparse} & $O\left(1/\sqrt{NKM}\right)$ & NA & Yes & Yes & Not data from server to worker \\
FedPAQ \cite{FedPAQ} & $O\left(1/\sqrt{NK}\right)$ & NA &  No & Yes  & Not data from server to worker \\
\textbf{This paper} & $O\left(1/\sqrt{NKM}\right)$  &  $O\left(1/\sqrt{NK^2M}\right)$ & No & Yes  &  All communication data  \\ 
\bottomrule
\end{tabularx}
\caption{Convergence rate for non-convex objectives under the best case after $N$ global synchronizations}
\label{table:1}
\end{table*}



As a server-based architecture, Parameter Server (PS) paradigm is an aggregation model and has been widely accepted because algorithms built on it are easy to implement and maintain \cite{hzwang_2,xie2019asynchronous}. PS-empowered classical PSGD\footnote{In the classical PSGD, workers compute local stochastic gradients and follow Vanilla SGD steps to update the model parameters in accordance with others' gradients. Workers can receive all gradients from other nodes under PS/AR paradigm, while learning only the information from neighbours under gossip paradigm.} achieves the convergence rate of $O\left(1/\sqrt{NM}\right)$ for non-convex objectives, where $N$ and $M$ refers to the times of global synchronizations and the number of workers, respectively \cite{ps_1,dekel2012optimal,ps_2}. 

While the PS paradigm seems to be efficient, it is not friendly in terms of communication as occurrences of network congestion degrade the training performance \cite{dmlsurvey}. To overcome the problem, a serverless architecture is an alternative choice. All-Reduce (AR) paradigm is one of the serverless models which releases the burden of the central node while achieving obtaining aggregation result same as the one under PS paradigm \cite{yu2019distributed}. Ring AR (RAR) paradigm \cite{baidu,horovod}, one of its successful examples, utilizes ring network topology, where all clients simultaneously process data whose total volume is $M$ times smaller than the one handled by the dedicated parameter server \cite{ar_analysis}. Besides, the gossip paradigm offers a solution for an arbitrary network topology. Gossip-based classical PSGD, also known as decentralized PSGD (D-PSGD) \cite{dpsgd}, achieves the same convergence rate as the one under PS, although it does not utilize all gradients throughout the network. 

How to reduce communication overhead is a crucial consideration for serverless paradigms. Classical PSGD requires model synchronization at every iteration. As a result, the algorithm not only consumes a large proportion of bandwidth throughout the training process, but also requires a great amount of time in the communication phase. Generally, there are two ways to increase communication efficiency: (1) reducing the frequency of synchronization and (2) compressing the traffic data in each transmission. The first approach, referring to parallel restarted (PR) SGD, is equivalent to periodical averaging SGD \cite{kavg_1,kavg_2,local_sgd}. Instead of exchanging model updates at each iteration, workers synchronize averaged results only once after aggregating individual solutions for multiple iterations. Under the PS paradigm, the best algorithm can converge at a rate of $O(1/\sqrt{NKM})$, where $K$ denotes the number of local updates \cite{prsgd}. Apparently, it achieves a linear speedup with respect to both the number of local updates and the number of workers. A similar result is obtained under the gossip paradigm \cite{li2019communicationefficient}. Gradient quantization \cite{qsgd,qsgd_tern_wen,qsgd_mean_suresh,yu2018gradiveq} is one of the common strategies for the second approach. It sacrifices the gradient precision to reduce bandwidth consumption as well as communication overhead. DoubleSqueeze \cite{tang2019doublesqueeze} and DCD-PSGD \cite{dcdpsgd} achieve gradient quantization at every transmission under the PS and gossip paradigms, respectively. Their convergence results are shown in Table \ref{table:1}. Under the AR paradigm, ECQSGD in \cite{wu2018error} is shown to converge when training a deep neural network. 

To further reduce consumption of communication resources, in this paper, we elaborately design Quantized Parallel Restarted SGD (QPRSGD) for two serverless aggregation paradigms -- RAR and gossip -- named AR-QPRSGD and G-QPRSGD, respectively. These two algorithms support precision-loss gradients synchronization after periodical averaging, in which the information exchanged between any two nodes is compressed and distorted. Theoretical analysis presents that our proposed algorithms retain the best convergence rate, while empirical studies indicate that they perform well under a network-intensive environment. Our contributions are listed as follows:

\begin{itemize} 
    \item Under the RAR paradigm, we prove that AR-QPRSGD for non-convex objectives achieves a convergence rate of $O(1/\sqrt{NKM})$, which indicates a linear speedup with respect to the number of workers and the number of local updates. To the best of our knowledge, this is the first work that investigates the PR-SGD with considering quantization in practical implementation of the RAR paradigm; 
    \item We propose the G-QPRSGD algorithm to evaluate PR-SGD with lossy-compression in the gossip paradigm. Our theoretical analysis shows that G-QPRSGD achieves an improved convergence rate over the state-of-the-art compression methods in the gossip training. Moreover, G-QPRSGD preserves the linear speedup with respect to the number of workers and the number of local updates, which ensures the effectiveness of local updates and the scalability of G-QPRSGD in the gossip paradigm; 
    \item We conduct an empirical study to illustrate the effect of our proposed algorithms. In terms of convergence rate, QPRSGD achieves up to 3.8 times and 4.7 times convergence efficiency in comparison with QSGD and PR-SGD, respectively. In addition, it reduces the communication cost by more than 90\% as compared with the PR-SGD. 
\end{itemize}

The rest of the paper is organized as follows: In Section \ref{section:related_work}, related work is introduced to provide an overview of three paradigms. Section \ref{section:preliminary} introduces the two basic SGD models and key notations. We analyze the convergence rate and the communication cost of Quantized-PR-SGD under the PS, RAR and Gossip paradigms in Section \ref{section:QPrSGD}. An empirical study is presented in Section \ref{section:experiment} to validate our theoretical analysis. Section \ref{section:conclusion} concludes this study. 

\section{Related Work} \label{section:related_work}



Aggregation paradigm plays an important role in distributed machine learning as it affects the performance of computation and communication. Both PS and AR paradigms require the specific underlying network topology, while gossip paradigm can work on arbitrary network topology. This section reviews how classical PSGD runs under three aggregation models, i.e., PS, AR and gossip paradigms. 

\noindent \textbf{Parameter Server (PS).} \quad
The PS is one of the most common centralized paradigms for large-scale distributed training. It typically consists of one or more server nodes and multiple worker nodes, each of which carries a subset of training data. The worker nodes firstly compute the stochastic gradients in parallel based on the local dataset. Then, the server node aggregates and averages the gradients sent from the workers. The worker nodes subsequently update their parameters using the averaged gradients. Repeat these three steps until the model converges. Since server nodes handle all communications, the performance of PS is largely determined by the communication resources of server nodes. 

Arguably, there is a work \cite{FedPAQ} achieving the similar features to ours under PS paradigm. It indicates the impact is trivial that the periodical averaging process will amplify the precision loss generated by quantization. However, each worker synchronizes the compressed model only with neighbours, leading to recursive compression and knowledge missing under RAR and gossip paradigms, respectively. Therefore, the compression errors possibly accumulate and spread over the whole network, which seriously affects the convergence ability. In addition, although \cite{FedPAQ} eventually converges at a stable stage, its theoretical convergence rate $O(1/\sqrt{NK})$ is not the best result compared to the existing works \cite{basu2019qsparse,prsgd}. Our proposed algorithms significantly improve the result, where they further realize linear speedup with respect to the number of workers. 






\begin{figure*}[htbp]
\centering

\begin{minipage}[t]{1.0\textwidth}
\centering
\includegraphics[width=1.0\textwidth]{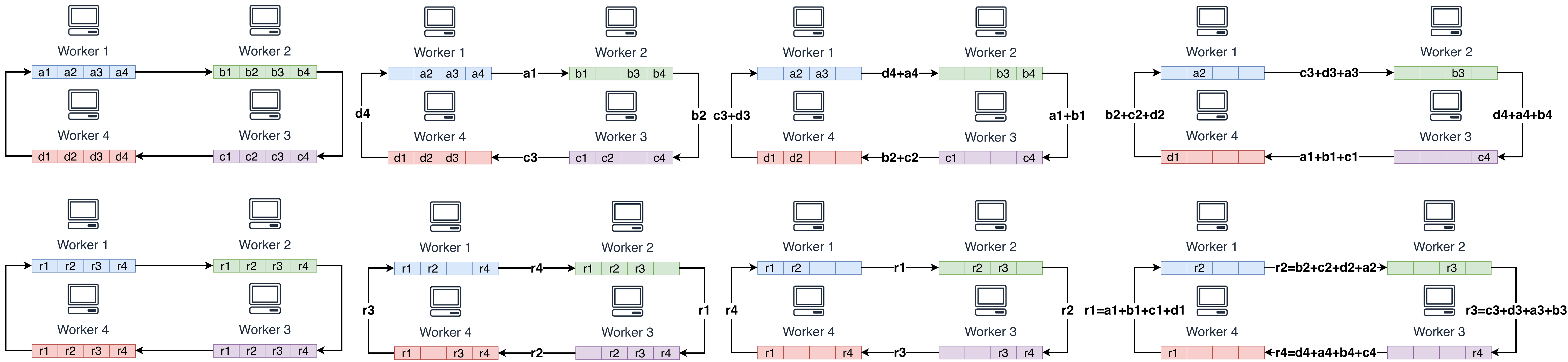}
\end{minipage}
\centering
\caption{Ring All-Reduce Paradigm with 4 workers. The process of \textit{reduce-scatter} is from left to right at the top, while the workflow of \textit{all-gather} is from right to left at the bottom.}
\label{fig:RAR}
\end{figure*}


\noindent \textbf{All-Reduce (AR).} \quad
The workers are able to preserve a consistent model using the AR paradigm without introducing central nodes \cite{ar_1_patarsuk}. Different from PS paradigm, its aggregation process constitutes with two phases -- \textit{reduce-scatter} and \textit{all-gather}. \textit{Reduce-scatter} is the process that some nodes maintain the sum of partial model parameters (or gradients), while \textit{all-gather} lets all nodes acquire updated parameters (or gradients). Their implementation details depend on the underlying network topology. Commonly, there are three types of underlying network topology, i.e., complete graph \cite{yu2019distributed}, ring \cite{horovod} and tree \cite{wan2020rat}. Without parameters (or gradients) compression, the convergence rate should be the same as the one under PS paradigm without extra communication overhead. Besides, clients process the same amount of data at a time whatever the topology uses. 


%



\noindent \textbf{Gossip.} \quad
As a fully distributed model, gossip paradigm possesses higher scalability than a centralized model and eliminates the risk of single points of failure. Peers under gossip paradigm compute a local stochastic gradient based on the holding data. Then, every worker exchanges the parameter with its neighbours and updates the local parameters with the existing knowledge. Apparently, the parameters of a model vary among clients. Compared to PS paradigm, gossip paradigm retains the same convergence rate with less communication overhead \cite{dpsgd}. 

\section{Preliminary} \label{section:preliminary}


\noindent \textbf{Problem Setting.} \quad
Generally, the objective of distributed machine learning is to minimize the cumulative expected loss over all workers, which can be formulated as 
\begin{equation} \label{equation:problem}
    \min_{\bm{x} \in \mathbb{R}^d} \quad F(\bm{x}) = \frac{1}{M} \sum_{m=1}^{M} \underbrace{\mathbb{E}_{\xi_m \sim \mathcal{D}_m} \left[f_m(\bm{x}, \xi_m)\right]}_{:= F_m(\bm{x})},
\end{equation}
where $M$ is the number of workers (or nodes), $\mathcal{D}_m$ is the local data distribution for worker $m$, $f_m(\bm{x}, \xi_m)$ is the empirical loss by the given parameter $\bm{x}$ and by the stochastic sample $\xi_m$ from $\mathcal{D}_m$, and $F_m(\cdot)$ is an objective function. Since the objective function is abstract over a given data distribution, it is a common practise that the bias does not exist between the expected loss and the empirical one. 






\noindent \textbf{Quantization Method.} \quad
Quantization compresses gradients that are exchanged through the network, while generally preserving the model convergence performance of optimization. In this paper, we adopt generally used QSGD \cite{qsgd} as the quantized function. For the quantization of any vector $\bm{v} \in \mathbb{R}^d$, QSGD relies on a finite set $\mathcal{S}$ which is defined as 
\begin{equation}
    \mathcal{S}=\{\bm{v}\|_2 \cdot \ell/s | \ell\ \text{is an integer and}\ \ell \in [-s, s]\},
\end{equation}
where $s$ is a pre-defined hyper-parameter. It can be seen that the size of the finite set $\mathcal{S}$ is $2s+1$, and thus each element in the set can be precisely presented with only several bits via an encoded function when $s$ is small. Given the set $\mathcal{S}$, the quantization function $Q_s(\cdot)$ of QSGD is to map each dimension $\bm{v}_i$ from a $32$-bits float to an element in the set:

\begin{equation} \label{equation:quantization_function}
\begin{split}
    Q_s(\bm{v}_i) = \|\bm{v}\|_2 \cdot \text{sgn}(v_i) \cdot \zeta(v_i, s) / s
\end{split}
\end{equation}


\noindent where $\text{sgn}(\cdot) \in \{-1, 1\}$ represents the sign bit of a real number. In (\ref{equation:quantization_function}) $\zeta(v_i, s)$ is defined as follows: Let $\ell$ be an integer such that $|v_i|/\|\bm{v}\|_2 \in [\ell/s, (\ell+1)/s]$ and 

\begin{equation} \label{equation:definition_ell}
    \zeta(v_i, s) = 
     \begin{cases}
        \ell + 1, & \text{with probability } p = \frac{|v_i|}{\|\bm{v}\|_2}s - \ell \\
        \ell, & \text{otherwise.}
     \end{cases}
\end{equation}

\noindent \textbf{Definition and Notation. }\quad
We use the following notations throughout this paper: 

\begin{itemize} \setlength\itemsep{0em}
    \item $\partial F(X) := \left[\nabla F_1\left(\bm{x}^{(1)}\right)  \quad ... \quad \nabla F_M\left(\bm{x}^{(M)}\right)\right]$ denotes the gradient tuple of the local loss function of each worker
    \item $F_*$ denotes optimal solution to Problem (\ref{equation:problem}) 
    \item $\|\cdot\|_n$ denotes $\ell_n$ norm of a vector in $\mathbb{R}^d$ 
    \item $\|\cdot\|_F$ denotes the Frobenius norm in a matrix 
    \item $\bm{1}_n$ denotes an $n$-dimension column vector filled with $1$
    \item $\lambda_i(\cdot)$ denotes the $i$-th largest eigenvalue of a matrix 
\end{itemize}
\section{Quantized Parallel Restarted SGD} \label{section:QPrSGD}

In this section, we propose QPRSGD that minimizes the communication overhead by reducing the frequency of synchronization and compressing information exchanged among nodes. We design two algorithms, named AR-QPRSGD and G-QPRSGD, for the RAR paradigm and gossip paradigm, respectively. 

\subsection{Ring All-Reduced (RAR) Paradigm} 


\subsubsection{Model Design}

Ring network topology is one of the common architecture for AR paradigm. Fig.~\ref{fig:RAR} illustrates the procedure of RAR paradigm with four workers. To elaborate the figure, we assume the system is to aggregate the gradients. As is described in the caption, the entire process can be divided into two parts -- \textit{reduce-scatter} and \textit{all-gather}. Initially, these four workers split the local stochastic gradient into four segments. Then, a worker repeats the following steps for three times in \textit{reduce-scatter}: receiving a segment from the last worker, adding it to specified segment, and passing the summation result to the next worker. Eventually, every worker at a specific position possesses a segment that is the sum of all workers. Next in \textit{all-gather}, a worker obtains the gathered segment from the last worker and sends it to the next worker until the worker gets all segments. 

Next, we further explore a more general case, where it makes up with $M$ workers and supports gradient compression. Firstly, a worker divides a stochastic gradient into $M$ segments. During \textit{reduce-scatter} period, every transmission is quantized such that every segment is recursively summed and compressed. In \textit{all-gather} period, every worker transmits a segment without compression as it has followed quantization format. Finally, all workers obtain the same compressed summation gradient. 

In order to let the aforementioned case support multiple local updates, we design the following steps between two successive synchronizations for a node: 

\begin{itemize} 
    \item \textbf{(Pull):} pull the parameter $\tilde{\bm{x}}$ from the last update as initial model parameter $\bm{x}_{0}$;
    \item \textbf{(Compute):} iterate the steps from $k=0$ to $k=K-1$: 
    \begin{itemize}
        \item Generate a realization of random samples $\xi_k$,
        \item Compute the gradient $\nabla f(\bm{x}_k; \xi_k)$ with samples;
        \item Update the model parameter with step size $\gamma$ and gradient: $\bm{x}_{k+1} = \bm{x}_{k} - \gamma \nabla f(\bm{x}_k; \xi_k)$; 
    \end{itemize}
    \item \textbf{(Push):} push the quantized update $g$ by $g = Q_s(\bm{x}_0 - \bm{x}_K)$;
    \item \textbf{(Aggregate):} aggregate the averaging stochastic gradients $g$ from all other nodes and summarize them into $\Delta$;
    \item \textbf{(Update):} update the parameter $\tilde{\bm{x}}$ by $\tilde{\bm{x}} = \tilde{\bm{x}} - \Delta$. 
\end{itemize}
The full AR-QPRSGD algorithm is presented in Algorithm \ref{algo:QPRSGD_AR}. 


\begin{algorithm}[t]
\SetKwData{Left}{left}\SetKwData{This}{this}\SetKwData{Up}{up}
\SetKwFunction{Union}{Union}\SetKwFunction{FindCompress}{FindCompress}
\SetKwInOut{Input}{Input}\SetKwInOut{Output}{Output}
\Input{Initial Point $\tilde{\bm{x}}_1$, stepsize series $\{\gamma_n\}$, the interval value $K$, and the number of total iterations $N$}
\BlankLine
\For{$  n\leftarrow 1$ \KwTo $N$}{
Initiate the first parameter of an epoch by $\bm{x}_{n;0}^{(m)} \leftarrow \tilde{\bm{x}}_n$\;
\For{$k\leftarrow 0$ \KwTo $K-1$}{
Randomly sample $\xi_k^{(m)}$ from local data $\mathcal{D}_m$\;
Compute local stochastic gradient $\nabla f_m(\bm{x}_{n;k}^{(m)}; \xi_{k}^{(m)})$ and update the local parameters via $\bm{x}_{n;k+1}^{(m)} \leftarrow \bm{x}_{n;k}^{(m)} - \gamma_n \nabla f_m(\bm{x}_{n;k}^{(m)}; \xi_{k}^{(m)})$\;
}

Calculate the update by $g_n^{(m)} \leftarrow \bm{x}_{n;0}^{(m)}-\bm{x}_{n;K}^{(m)}$\;
Initiate the starting assembler $\Delta_{m} \leftarrow \bm{0}$\;
\For{$i \leftarrow 0$ \KwTo $M-1$}{
$j \leftarrow (m + i)$ mod $M$\;
 Select $j$-th part of $g_n^{(m)}$ (a.k.a. $g_{n, j}^{(m)}$)\;
Receive $\Delta_{j}$ from last node and update it with $g_{n, j}^{(m)}$: $\Delta_{j} \leftarrow Q_s(\Delta_{j} + g_{n, j}^{(m)})$ \;
Send $\Delta_{j}$ to the next node\;
}
Broadcast $\Delta_m$ to other nodes\;
Assemble $\Delta_m$ for $m \in \{1, ..., M\}$ into $\Delta$ and average it with $\Delta \leftarrow \Delta/M$\;
Update the parameters through $\tilde{\bm{x}}_{n+1} \leftarrow \tilde{\bm{x}}_n - \Delta$\;
}
\caption{AR-QPRSGD (Worker $m$)}\label{algo:QPRSGD_AR}
\end{algorithm}


\subsubsection{Theoretical Analysis} 

Prior to evaluating the convergence property of Algorithm \ref{algo:QPRSGD_AR}, Lemma \ref{lemma:1} presents the update between the recursive compression and the average of the segments at the specific positive under a general case. 
\begin{lemma} \label{lemma:1}
For any $i \in \{1, 2, ..., M\}$, $m \in \{1, 2, ..., M\}$ and vector $\bm{w}_i^{(m)} \in \mathbb{R}^{d/M}$ which means $i$-th segment of $m$-th worker, under RAR paradigm, we denote the recursive compression as 
\begin{equation*}
    X = \frac{1}{M} Q_s\left(...Q_s\left(Q_s\left(\bm{w}_i^{(1)}\right)+\bm{w}_i^{(2)}\right)+...+\bm{w}_i^{(M)}\right). 
\end{equation*}
Then, we have 
\begin{equation} \label{eq:4}
\begin{split}
    \mathbb{E} \left\| X - \frac{1}{M} \sum_{m=1}^M \bm{w}_i^{(m)} \right\|_2^2 \leq \frac{2C_1}{M} \sum_{m=1}^M \left\|\bm{w}_i^{(m)}\right\|_2^2, 
\end{split}
\end{equation}
where 
\begin{equation*}
    C_1 := \frac{1}{M}\exp\left({\frac{d}{4s^2}}\right) + \frac{4s^2}{d}\left(\frac{d}{4s^2M}+1\right)^{M-1} - \left(\frac{4s^2}{d} + 1\right). 
\end{equation*}
\end{lemma}

Lemma \ref{lemma:1} indicates that permutation and division operations do not influence the convergence result, because the coefficient on the right hand side does not depend on the segments. The following remark is to draw a conclusion from the perspective of all segments. 

\begin{remark} \label{remark:1}
For any $m \in \{1, 2, ..., M\}$ and vector $\bm{w}^{(m)} \in \mathbb{R}^{d}$ which is independent with others, under RAR paradigm, we represent the noise as: 
\begin{equation*}
\begin{split}
    \Delta X &= \frac{1}{M} Q_s\left(...Q_s\left(Q_s\left(\bm{w}^{(1)}\right)+\bm{w}^{(2)}\right)+...+\bm{w}^{(M)}\right)\\
    & \quad - \frac{1}{M} \sum_{m=1}^M \bm{w}^{(m)}
\end{split}
\end{equation*}
Then, the second moment of the noise should be 
\begin{equation} \label{eq:5}
\begin{split}
    \mathbb{E} \left\| \Delta X \right\|_2^2 \leq \frac{2C_1}{M} \sum_{m=1}^M \left\|\bm{w}^{(m)}\right\|_2^2
\end{split}
\end{equation}
where $C_1$ is the same as the definition in Lemma \ref{lemma:1}. 
\end{remark}

In Equation \ref{eq:5}, the complexity of the bounded noise is $O(1/M)$. The following remark achieves $O\left(1/M^2\right)$ by bounding the number of workers, which can improve the performance of the convergence result. 

\begin{remark}
Under Remark \ref{remark:1}, given that 
\begin{equation} \label{eq:worker_bound}
    M \leq d \left(1 + \frac{d}{4s^2}\right) \left(\exp\left(\frac{d}{4s^2}\right) - \frac{d}{4s^2} - 1 \right)^{-1},
\end{equation}
we have 
\begin{equation}
\begin{split}
    & \quad \mathbb{E} \left\| \frac{1}{M} Q_s\left(...Q_s\left(Q_s\left(\bm{w}^{(1)}\right)+\bm{w}^{(2)}\right)+...+\bm{w}^{(M)}\right) \right\|_2^2 \\
    & \leq \frac{2}{M^2} \exp\left({\frac{d}{4s^2}}\right) \sum_{m=1}^M \left\|\bm{w}^{(m)}\right\|_2^2 + \left\| \frac{1}{M} \sum_{m=1}^M \bm{w}^{(m)} \right\|_2^2.
\end{split}
\end{equation}
\end{remark}

The quantization method we use is unbiased \cite{qsgd}. Knowing that the compression order does not lead to a substantial impact, we can consider a general recursion formula for Algorithm \ref{algo:QPRSGD_AR} as: 

\begin{gather*}
    \bm{x}_{n;0}^{(m)}  = \tilde{\bm{x}}_n, \quad
    \bm{x}_{n;t}^{(m)}  = \tilde{\bm{x}}_n - \gamma_n \sum_{j=0}^{t-1} \nabla f_m\left(\bm{x}_{n;j}^{(m)}, \xi_j^{(m)}\right), \\
    \tilde{\bm{x}}_{n+1} = \tilde{\bm{x}}_n - \frac{1}{M} \underbrace{Q_s\left(...Q_s\left(Q_s\left(g_n^{(1)}\right)+g_n^{(2)}\right)+...\right)}_{\text{$M$ gradients from $M$ workers}}, \\
    g_n^{(i)} = \gamma_n \sum_{k=0}^{K-1} \nabla f_i\left(\bm{x}_{n;k}^{(i)}, \xi_k^{(i)}\right), \quad \Bar{\bm{x}}_{n;t} = \frac{1}{M} \sum_{m=1}^M \bm{x}_{n;t}^{(m)}
\end{gather*}

We make the following assumptions, which are ubiquitously used for SGD-based distributed optimization \cite{optimization_method,dpsgd}: 

\begin{assumption} \label{assumption:1}
Problem (\ref{equation:problem}) satisfies the following constraints: 

\begin{enumerate}
    \item \textbf{Smoothness:} All function $F_m(\cdot)$'s are continuous differentiable and their gradient functions are $L$-Lipschitz continuous with $L>0$; 
    \item \textbf{Bounded variance:} For any worker $m$ and vector $\bm{x} \in \mathbb{R}^d$, there exist scalars $\sigma \geq 0$ and $\kappa \geq 0$ such that
    \begin{align*}
        \mathbb{E}_{\xi \sim \mathcal{D}_m} \|\nabla f_m(\bm{x}, \xi) - \nabla F_m (\bm{x})\|^2_2 \leq \sigma^2, \\
        \frac{1}{M} \sum_{m=1}^M \|\nabla F(\bm{x}) - \nabla F_m(\bm{x})\|_2^2 \leq \kappa^2. 
    \end{align*}
\end{enumerate}

\end{assumption}

According to the recursion formula, we figure out the convergence rate for non-convex objectives with a fixed stepsize under Assumption 1. 

\begin{theorem} \label{theorem:1}
Consider Problem (\ref{equation:problem}) for non-convex optimization. Suppose Algorithm \ref{algo:QPRSGD_AR} is running with a constant stepsize $\Bar{\gamma}$ satisfying the three inequalities to avoid the probability that any denominators are being 0: 
\begin{gather*}
 \exists \epsilon \in (0, 1),  \frac{(8C_1L^2\Bar{\gamma}^2 + (K-1)L\Bar{\gamma} + 2C_1)KL\Bar{\gamma}}{1-2(K+1)(K-2)L^2\Bar{\gamma}^2}\leq \frac{1 - \epsilon}{2},  \\
 \exists \hspace{0.1cm} \delta \in (0, 1), \quad 2(K+1)(K-2)L^2\Bar{\gamma}^2 \leq 1 - \delta, \\
 1 - L\Bar{\gamma} - \frac{\left(4KL\Bar{\gamma}C_1 + 1\right)(K-1)KL^2\Bar{\gamma}^2}{1-2(K+1)(K-2)L^2\Bar{\gamma}^2} > 0. 
\end{gather*}
Then, under Assumption \ref{assumption:1}, for all $N \geq 1$, we have 
\begin{equation} \label{eq:9}
\begin{split}
    & \quad \frac{1}{NK} \sum_{n=1}^N \sum_{k=0}^{K-1} \mathbb{E} \left\| \nabla F(\Bar{\bm{x}}_{n, k}) \right\|_2^2 \\
    & \leq \frac{\left(4KL\Bar{\gamma}C_1 + 1\right)}{\delta \epsilon} \left(\frac{\sigma^2(M+1)}{2M}+ 2K\kappa^2\right)(K-1)L^2\Bar{\gamma}^2 \\
    & \quad + \frac{1}{\epsilon} \left(2C_1\sigma^2 + 4KC_1\kappa^2 + \frac{\sigma^2}{M}\right) L\Bar{\gamma}  + \frac{2\left[F(\tilde{\bm{x}}_1) - F_*\right]}{NK\Bar{\gamma}\epsilon}, 
\end{split}
\end{equation}
where $C_1$ is the same as the definition in Lemma \ref{lemma:1}. 
\begin{proof}
Due to the page limitation, we only provide the sketch of the proof in the main body. For more details, please refers to \url{https://arxiv.org/abs/2004.09125}. To prove whether Algorithm \ref{algo:QPRSGD_AR} converges under a non-convex objective $F$, a common practice is to achieve the following formula: 
\begin{equation*}
    \lim_{N \rightarrow \infty}\frac{1}{NK} \sum_{n=1}^N \sum_{k=0}^{K-1} \mathbb{E} \|\nabla F(\Bar{\bm{x}}_{n, k})\|_2^2 = 0. 
\end{equation*}
As the above term is accumulated for all $N \in \{1, ..., N\}$, with an optimal solution $F_*$, we have 
\begin{equation*}
    F_* - F(\tilde{x}_1) \leq \mathbb{E} [F(\tilde{x}_{N+1}) - F(\tilde{x}_1)] = \sum_{n=1}^N \mathbb{E} [F(\tilde{x}_{n+1}) - F(\tilde{x}_n)]
\end{equation*}
Next, we should find the boundary for $\mathbb{E} [F(\tilde{x}_{n+1}) - F(\tilde{x}_n)]$, where it will appear the accumulated term $\mathbb{E} \|\nabla F(\Bar{\bm{x}}_{n, k})\|_2^2$ for all $k \in \{0, ..., K-1\}$. According to $L$-smooth definition in Assumption \ref{assumption:1}, we have: 
\begin{equation} \label{eq:m9}
\begin{split}
    & \quad \mathbb{E} \left[F(\tilde{x}_{n+1}) - F(\tilde{x}_n)\right] \\
    &= \mathbb{E} \left[F(\tilde{x}_{n+1}) - F(\bar{x}_{n;K})\right] + \sum_{k=0}^{K-1} \mathbb{E} [F(\bar{x}_{n;k+1}) - F(\bar{x}_{n;k})] \\
    & \leq \mathbb{E} \left \langle \nabla F(\bar{x}_{n;K}), \tilde{x}_{n+1} - \bar{x}_{n;K} \right \rangle 
        + \frac{L}{2} \mathbb{E} \| \tilde{x}_{n+1} - \bar{x}_{n;K} \|_2^2\\
        &\quad + \sum_{k=0}^{K-1}\left(
            \mathbb{E} \left \langle \nabla F(\bar{x}_{n;k}), \bar{x}_{n;k+1} - \bar{x}_{n;k} \right \rangle
            + \frac{L}{2} \mathbb{E} \| \bar{x}_{n;k+1} - \bar{x}_{n;k} \|_2^2
        \right)\\
    & \overset{(a)}{\leq} \frac{C_1 L \gamma_n^2}{M}  \sum_{m=1} ^ M  \underbrace{ \mathbb{E} \left\| \sum_{j=0}^{K-1} \nabla f_m\left(x_{n;j}^{(m)}, \xi_j^{(m)}\right) \right\|_2^2}_{T_1} \\
    & \quad - \gamma_n \sum_{k=0}^{K-1} \underbrace{\mathbb{E} \left \langle \nabla 
    F\left(\bar{x}_{n;k}\right), \frac{1}{M} \sum_{m=1}^{M} \nabla F_m\left(x_{n;k}^{(m)}\right)
    \right \rangle}_{T_2} \\
    &\quad + \frac{L\gamma_n^2}{2M^2} \sum_{k=0}^{K-1} \underbrace{ \mathbb{E} \left\| \sum_{m=1}^{M} \nabla f_m\left(x_{n;k}^{(m)}, \xi_k^{(m)}\right) \right\|_2^2}_{T_3}
\end{split}
\end{equation}
where $(a)$ is according to unbiased feature of QSGD \cite{qsgd} and Remark \ref{remark:1}. Then, by analyzing $T_1$, $T_2$ and $T_3$, we have: 
\begin{align*}
    T_1 &\leq K \sigma^2 + 2 K L^2 \sum_{j=0}^{K-1} \left\| x_{n;j}^{(m)} - \bar{x}_{n;j} \right\|_2^2 + 2 K \sum_{j=0}^{K-1} \left\| \nabla F_m\left(\bar{x}_{n;j}\right) \right\|_2^2, \\
    T_2 &\leq -\frac{1}{2} \left\| \nabla F(\bar{x}_{n;k})\right\|_2^2
    - \frac{1}{2} \left\| \frac{1}{M} \sum_{m=1}^M \nabla F_m(\bar{x}_{n;k})\right\|_2^2\\
    &\quad + \frac{L^2}{2M} \sum_{m=1}^M \left\| \bar{x}_{n;k} - x_{n;k}^{(m)}\right\|_2^2, \\
    T_3 &\leq M \sigma^2 + \left\| \sum_{m=1}^M \nabla F_m\left(x_{n;k}^{(m)}\right)\right\|_2^2. 
\end{align*}
Given that the stepsize is a constant value and $F_*$ is the optimizer answer that we expect to obtain, knowing that 
\begin{align*}
    &\quad\frac{1}{M} \sum_{m=1}^M \sum_{k=0}^{K-1} \mathbb{E} \left\| x_{n;k}^{(m)} - \Bar{x}_{n;k} \right\|_2^2\\ 
    &\leq \frac{(K-1)K \gamma_n^2 \sigma^2 (M+1)}{2M [1 - 2(K+1)(K-2)L^2\gamma_n^2]}\\
    & \quad + \frac{(K-1)K\gamma_n^2}{1 - 2(K+1)(K-2)L^2\gamma_n^2} \sum_{k=0}^{K-1} \left\|\frac{1}{M} \sum_{m=1}^M \nabla F_m \left(x_{n;k}^{(m)}\right)\right\|_2^2\\
    & \quad + \frac{2(K-1)K\gamma_n^2}{1 - 2(K+1)(K-2)L^2\gamma_n^2} \sum_{k=0}^{K-1} \left(\kappa^2 + \|\nabla F(\Bar{x}_{n;k})\|_2^2\right), 
\end{align*}
we plug $T_1$, $T_2$ and $T_3$ into Equation \ref{eq:m9} and obtain 
\begin{align*}
    & \quad F_* - F(\tilde{x}_1) \leq \mathbb{E} [F(\tilde{x}_{N+1}) - F(\tilde{x}_1)]\\
    & = \sum_{n=1}^N \mathbb{E} [F(\tilde{x}_{n+1}) - F(\tilde{x}_n)] \\
    & \leq -\frac{\Bar{\gamma}\varepsilon}{2} \sum_{n=1}^N \sum_{k=0}^{K-1} \|\nabla F(\Bar{x}_{n; k})\|_2^2 + C_1 L \Bar{\gamma}^2 K \sigma^2 N + 2LK^2\Bar{\gamma}^2C_1\kappa^2N\\ 
    & \quad + \frac{L\Bar{\gamma}^2K\sigma^2N}{2M} + \left(2KL^3\Bar{\gamma}^2C_1 + \frac{L^2\Bar{\gamma}}{2}\right) \left( \frac{(K-1)K\gamma_n^2\sigma^2(M+1)}{2M\delta} \right. \\
    & \left. \qquad + \frac{2(K-1)K^2\gamma_n^2\kappa^2}{\delta} \right) N.
\end{align*}

By polishing the formula above, we can obtain Equation \ref{eq:9}. Practically, the stepsize $\Bar{\gamma}$ is corresponding to $N$ and $K$. The following corollaries present the setting of the stepsize such that we can obtain different convergence rate. 
\end{proof}
\end{theorem}

To clearly show the convergence result derived in Theorem \ref{theorem:1}, we select an appropriate stepsize in the following corollary to achieves a sublinear convergence rate and a linear speedup property. The following corollaries treat all variables as hyper-parameters except $N$, $M$ and $K$ and assume that $N$ is sufficiently large. 

\begin{corollary}
Under Theorem \ref{theorem:1}, take the stepsiz $\Bar{\gamma} = \frac{1}{LK\sqrt{N}}$. Then, the output of Algorithm \ref{algo:QPRSGD_AR} achieves the ergodic convergence rate $O\left(1/\sqrt{N}\right)$.


\end{corollary}

\begin{corollary}
Under Theorem \ref{theorem:1}, given that the bound of local updates $K \leq \sigma^2/2\kappa^2$, we take the stepsize 
$\Bar{\gamma} = \frac{1}{L\sqrt{NK}}$.
\noindent Then, the output of Algorithm \ref{algo:QPRSGD_AR} achieves the ergodic convergence rate $O\left(1/\sqrt{NK}\right)$.


\end{corollary}

\begin{corollary} \label{corollary:1.3}
Under Theorem \ref{theorem:1}, given that the bound of the number of workers in Equation \ref{eq:worker_bound} and the bound of local updates $K \leq \frac{\sigma^2}{2\kappa^2} \left(1+2\exp\left(-\frac{d}{4s^2}\right)\right)$, we take the stepsize $\Bar{\gamma} = \frac{\sqrt{M}}{L\sqrt{NK}}$. Then, the output of Algorithm \ref{algo:QPRSGD_AR} achieves the ergodic convergence rate as $O\left(1/\sqrt{NKM}\right)$.


\end{corollary}

The last corollary suggests that our algorithm achieves linear speedup with respect to the number of workers and the number of local updates in the best case. Compared to the best case in \cite{FedPAQ}, which achieves a convergence rate of $O\left(1/\sqrt{NK}\right)$ under PS, our result presents a distinct dominance. 

\noindent \textbf{Impact of the quantization level $s$} \quad Referred to the constraints of Corollary \ref{corollary:1.3}, larger $s$ contributes to a wider range for $M$ and $K$. However, it probably requires more bits in every transmission. Apparently, there exists a trade-off between communication costs and the convergence result. 

\subsection{Gossip Paradigm}

\subsubsection{Model Design}

\begin{figure*}
    \centering
    \includegraphics[width=1.0\textwidth]{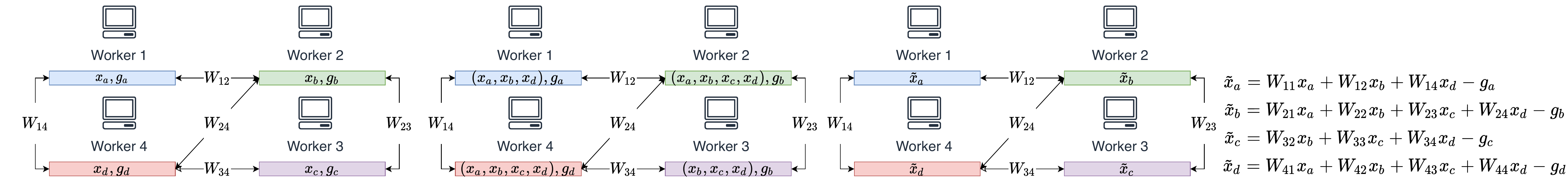}
    \caption{Gossip Paradigm with 4 workers. $W_{ij}$ indicates the weight between node $i$ and node $j$. }
    \label{fig:gossip}
\end{figure*}

RAR paradigm requires ring network topology or otherwise, it loses its dominance \cite{momentum_sgd}. Thus, we consider a more general case to implement QPRSGD in this section. Gossip paradigm satisfies all kinds of network topology, where it can be represented by an undirected graph with the value of: ($V$, $W$). $V \in \{1, ..., M\}$ denotes the set of $M$ workers. $W \in \mathbb{R}^{M \times M}$ is a doubly stochastic matrix, which satisfies (i) $W_{ij} \in [0, 1]$, (ii) $W = W^T$ and (iii) $\sum_j W_{ij} = 1$ for all $i$. Apparently, if there is no connection between node $i$ and node $j$, then both $W_{ij}$ and $W_{ji}$ are set to 0. Unlike the algorithms under PS or AR paradigm, workers transmit the parameters ahead of global update under gossip paradigm. Hence, as for Worker $m$ with $x_m$, the result integrates both the weighted average of the parameters with weights $\{W_{mj}, j \in \{1, ..., M\}\}$ and the stochastic gradient with respect to $x_m$. Figure \ref{fig:gossip} illustrates how gossip paradigm update the parameters with the network size of four workers. 


A recent research shows that the quantization method we use cannot directly compress the entire model parameters \cite{dcdpsgd}. Therefore, our proposed algorithm is to quantize and transmit the update of the model parameters, and each client tracks the update of its neighbours' parameters. As for Worker $m$, the steps between two successive synchronizations are presented as follows:  

\begin{itemize}
    \item \textbf{(Pull):} pull the parameter $\tilde{\bm{x}}^{(m)}$ from the last update as initial state $\bm{x}_{0}$; 
    \item \textbf{(Compute):} repeat the following steps for $K$ times: 
    \begin{itemize}
        \item Generate a realization of random variables $\xi_k, k=1, ..., K$, 
        \item Compute the gradient for the next iteration: $\bm{x}_{k+1} = \bm{x}_{k} - \gamma \nabla f(\bm{x}_k; \xi_k)$; 
    \end{itemize} 
    \item \textbf{(Store):} store the intermediate value by $\tilde{\bm{x}}_{0.5}^{(m)} = \sum_{j=1}^M W_{mj}\tilde{\bm{x}}^{(j)} - \bm{x}_{0} + \bm{x}_{K}$;
    \item \textbf{(Push):} push the quantized update $g^{(m)}$ by $g^{(m)} = \tilde{\bm{x}}^{(m)} - \tilde{\bm{x}}_{0.5}^{(m)}$;
    \item \textbf{(Update):} update the parameters $\tilde{\bm{x}}^{(j)}$ by $\tilde{\bm{x}}^{(j)} = \tilde{\bm{x}}^{(j)} - g^{(j)}$ for all connected neighbours (i.e. $W_{mj} \neq 0$). 
\end{itemize}
The full realization of G-QPRSGD is presented in Algorithm \ref{algo:G-QPRSGD}. 

\IncMargin{0em}
\begin{algorithm}[t!]
\SetKwData{Left}{left}\SetKwData{This}{this}\SetKwData{Up}{up}
\SetKwFunction{Union}{Union}\SetKwFunction{FindCompress}{FindCompress}
\SetKwInOut{Input}{Input}\SetKwInOut{Output}{Output}
\Input{Initial Point $\tilde{\bm{x}}_1^{(j)} = \bm{x}_1, \forall j \in \{1, ..., M\}$, stepsize series $\{\gamma_n\}$, weighted matrix $W$, the interval value $K$, and the number of total iterations $N$}
\BlankLine
\For{$  n\leftarrow 1$ \KwTo $N$}{
Initiate the first parameter of an epoch by $\bm{x}_{n;0}^{(m)} \leftarrow \tilde{\bm{x}}_n^{(m)}$\;
\For{$k\leftarrow 0$ \KwTo $K-1$}{
Randomly sample $\xi_k^{(m)}$ from local data $\mathcal{D}_m$\;
Compute local stochastic gradient $\nabla f_m(\bm{x}_{n;k}^{(m)}; \xi_{k}^{(m)})$ and update the local parameters via $\bm{x}_{n;k+1}^{(m)} \leftarrow \bm{x}_{n;k}^{(m)} - \gamma_n \nabla f_m(\bm{x}_{n;k}^{(m)}; \xi_{k}^{(m)})$\;
}

Temporarily update the local model with $\tilde{\bm{x}}_{n+\frac{1}{2}}^{(m)} \leftarrow \sum_{j=1}^M W_{mj}\tilde{\bm{x}}_n^{(j)} - \bm{x}_{n;0}^{(m)} + \bm{x}_{n; K}^{(m)}$ \;

    

Calculate and compress the update by $g_n^{(m)} \leftarrow  Q_s(\tilde{\bm{x}}_n^{(m)} - \tilde{\bm{x}}_{n+\frac{1}{2}}^{(m)})$\;

Send $g_n^{(m)}$ to its neighbours\;
Receive $g_n^{(j)}, \forall j \in \{1, 2, ..., M\} \wedge W_{mj} \neq 0$ and update the gradient through $\tilde{\bm{x}}_{n+1}^{(j)} \leftarrow \tilde{\bm{x}}_n^{(j)} - g_n^{(j)}$\;
}
\caption{G-QPRSGD (Worker $m$)}\label{algo:G-QPRSGD}
\end{algorithm} 

\subsubsection{Theoretical Analysis}

It is challenging to directly find the recursion formula of Algorithm \ref{algo:G-QPRSGD} for two successive local updates. We first treat multiple local updates as one single update and then use bounded noise to replace the quantization loss. Therefore, the recursive function can be represented as follows: 

\begin{align} \label{eq:9}
\begin{split}
    X_{n;t} &= X_{n;t-1} - \gamma_n G(X_{n; t-1}, \xi_{t-1})\\
    &= X_n - \gamma_n \sum_{j=0}^{t-1} G(X_{n; j}, \xi_j) \\
    X_{n+1} &= X_n W - \gamma_n \sum_{k=0}^{K-1} G(X_{n; k}, \xi_k) + C_n 
\end{split}
\end{align}
where 
\begin{itemize}
    \item[$\bullet$]  $G(X_{n; k}, \xi_k) = \left[\nabla f_1\left(\bm{x}_{n;k}^{(1)}, \xi_k^{(1)}\right) \quad ... \quad \nabla f_M\left(\bm{x}_{n;k}^{(M)}, \xi_k^{(M)}\right)\right]$ represents a stochastic matrix under a local update;
    \item[$\bullet$] $\Delta X_n = X_n W - \gamma_n \sum_{k=0}^{K-1} G(X_{n; k}, \xi_k) - X_n$ represents the gradient difference; 
    \item[$\bullet$] $C_n = Q_s(\Delta X_n) - \Delta X_n$ represents the noise of compression with quantization. 
\end{itemize}


It is universally acknowledged that Assumption \ref{assumption:1} is not sufficient to analyze the convergence rate under gossip paradigm. As a result, we make an additional assumption, which is commonly adopted in \cite{dsgd}.

\begin{assumption}[Spectral gap] \label{assumption:2}
Given the symmetric doubly stochastic matrix $W \in \mathbb{R}^{M \times M}$, we assume $\rho := \max \{|\lambda_2(W)|, |\lambda_M(W)|\} < 1$. 
\end{assumption}

\allowdisplaybreaks

By the given recursive function, the convergence rate of Algorithm \ref{algo:G-QPRSGD} depends on the stochastic matrix and the noise of compression. Before analyzing the convergence property, we find the relationship between the accumulated noise and the stochastic matrix. 

\begin{lemma} \label{lemma:2}
Suppose $1 - \left(d \mu^2\right)/\left(s^2 \left(1-\rho\right)^2\right) > 0$. With the fixed stepsize, under Assumption \ref{assumption:1} and Assumption \ref{assumption:2}, the second moment of the noise of compression is: 
\begin{align*}
     \mathbb{E} \| C_n \|_F^2 \leq \Bar{\gamma}^2 K \left[\frac{d\left(1-\rho\right)^2 + 2d\mu^2}{2\left(s^2\left(1-\rho\right)^2 - d\mu^2\right)}\right] \sum_{k=0}^{K-1} \left\| G\left(X_{n;k}, \xi_k\right) \right\|_F^2
\end{align*}
\end{lemma}

The following lemma introduces the boundary for the stochastic matrix $G\left(X_{n;k}, \xi_k\right)$. 

\begin{lemma} \label{lemma:3}
Suppose $1 - \left(d \mu^2\right)/\left(s^2 \left(1-\rho\right)^2\right) > 0$. Given the fixed stepsize $\Bar{\gamma}$ satisfying: 
\begin{align*}
    1 - 8L^2\Bar{\gamma}^2KD’_1D’_3 > 0 \quad \text{and} \quad 1-12\Bar{\gamma}^2L^2\left(K+1\right)\left(K-2\right) > 0. 
\end{align*}
Under Assumption \ref{assumption:1} and Assumption \ref{assumption:2}, we have 
\begin{align*}
    & \quad \sum^N_{n=1}\sum^K_{k=0}\left\|G\left(X_{n;k}, \xi_k\right)\right\|^2_F \\
    & \leq \frac{MNK\left(\sigma^2+4\kappa^2\right) + 4NL^2\sigma^2\Bar{\gamma}^2K\left(K-1\right)}{1-8L^2\Bar{\gamma}^2KD’_1D’_3}\\
    &\quad +\frac{8MNK\left(K-1\right)L^2\Bar{\gamma}^2[\sigma^2+\kappa^2\left(2K-1\right)]}{\left(1-8L^2\Bar{\gamma}^2KD’_1D’_3\right)[1-12\Bar{\gamma}^2L^2\left(K+1\right)\left(K-2\right)]} \\
    &\quad+\frac{4MK\left(K-1\right)L^2\Bar{\gamma}^2}{1-8L^2\Bar{\gamma}^2KD’_1D’_3}\sum^N_{n=1}\sum^{K-1}_{k=0}\left\|\frac{\partial F\left(X_{n;k}\right) \cdot 1_M}{M}\right\|^2_2\\
    &\quad +\frac{4K}{1-8L^2\Bar{\gamma}^2KD’_1D’_3}\sum^N_{n=1}\left\|\nabla F\left(\frac{X_n \cdot 1_M}{M}\right)\right\|^2_2
\end{align*}
where 
\begin{align*}
    D_1^\prime &:= \frac{2K + 24\Bar{\gamma}^2 L^2 - 1}{1 - 12\Bar{\gamma}^2 L^2 (K+1) (K-2)}, \\
    D_2^\prime &:= \frac{d(1-\rho)^2 + 2d\mu^2}{2[s^2(1-\rho)^2 - d\mu^2]},\  D_3^\prime := \frac{D_2^\prime}{1-\rho^2} + \frac{1}{(1-\rho)^2}, \\
\end{align*}
\end{lemma}

Under Assumption \ref{assumption:1} and Assumption \ref{assumption:2}, Theorem \ref{theorem:2} introduces the convergence result for Algorithm \ref{algo:G-QPRSGD} for non-convex objectives with a constant stepsize. 

\begin{theorem} \label{theorem:2}
The weighted matrix $W$ in Algorithm \ref{algo:G-QPRSGD} is a symmetric double stochastic matrix satisfying $d\mu^2 - s^2(1-\rho)^2 < 0$, and initial point is $\bm{x}_1 = \bm{0}$. Under Assumption \ref{assumption:1} and Assumption \ref{assumption:2}, by choosing the fixed stepsize $\Bar{\gamma}$ with which the following inequalities hold:  
\begin{gather*}
    1 - 12 \Bar{\gamma}^2 L^2 (K+1) (K-2) > 0, \quad  \\
    1 - 8 L^2 \Bar{\gamma}^2 K D_1^\prime D_3^\prime > 0, \quad D_6^\prime > 0. 
\end{gather*}
Then for all $N \geq 1$, we have 
\begin{equation}
\begin{split}
    & \quad \frac{1}{N}\left(\left(1 - \frac{4D_5^\prime}{M}\right) \sum_{n=1}^N \mathbb{E} \left\|\nabla F \left(\frac{X_n\cdot\bm{1}_M}{M}\right)\right\|_2^2 \right. \\
    & \quad \quad \quad \left. + \frac{D_6^\prime}{K} \sum_{n=1}^N \mathbb{E} \left\|\frac{\partial F(X_n) \cdot \bm{1}_M}{M}\right\|_2^2\right) \\
    & \leq \frac{8(K-1)L^2\Bar{\gamma}^2D_5^\prime}{1 - 12 \Bar{\gamma}^2 L^2 (K+1)(K-2)} \left(\sigma^2 + 2(2K-1)\kappa^2 \right) \\
    & \quad + D_5^\prime \left(\sigma^2 + 4\kappa^2 + \frac{4L^2\Bar{\gamma}^2(K-1)}{M}\right) + \frac{2(F(\bm{x}_1) - F_*)}{\Bar{\gamma}KN} \\
\end{split}
\end{equation}
where $D_1^\prime$, $D_2^\prime$ and $D_3^\prime$ are the same as the definition in Lemma \ref{lemma:3}, and
\begin{align*}
    D_4^\prime &:= 2 \Bar{\gamma}LKD_3^\prime + \frac{\Bar{\gamma}L(K-1)}{2} + \frac{D_2^\prime}{2M}, \\
    D_5^\prime &:= \frac{2KL\Bar{\gamma}D_4^\prime}{1 - 8 L^2 \Bar{\gamma}^2 K D_1^\prime D_3^\prime}, \\
    D_6^\prime &:= 1 - 4K(K-1)L^2\Bar{\gamma}^2D_5^\prime - 2\Bar{\gamma}KL. 
\end{align*}
\begin{proof}
Constrained by the space, we only give the outline of the proof. Similar to Theorem \ref{theorem:1}, we should find a convergence rate, which is negatively corresponding to $N$, such that it tends to 0 as the total number of synchronizations increase. Following the proof sketch in Theorem \ref{theorem:1}, we utilize L-smooth feature in Assumption \ref{assumption:1} and obtain: 
\allowdisplaybreaks
\begin{align} \label{eq:m12}
    &\nonumber \quad \mathbb{E} \left[F\left(\frac{X_{n+1} \cdot 1_M}{M}\right)\right] - \mathbb{E} \left[F\left(\frac{X_{n} \cdot 1_M}{M}\right)\right] \\
    &\nonumber \leq \mathbb{E} \left\langle \nabla F\left(\frac{X_{n} \cdot 1_M}{M}\right), \frac{X_{n+1} \cdot 1_M}{M} - \frac{X_{n} \cdot 1_M}{M} \right\rangle \\
    &\nonumber \quad + \frac{L}{2} \mathbb{E} \left\|\frac{X_{n+1} \cdot 1_M}{M} - \frac{X_{n} \cdot 1_M}{M}\right\|_2^2 \\
    &\nonumber = \mathbb{E} \left\langle \nabla F\left(\frac{X_{n} \cdot 1_M}{M}\right), \frac{C_{n} \cdot 1_M}{M} \right\rangle \\
    &\nonumber \quad - \mathbb{E} \left\langle \nabla F\left(\frac{X_{n} \cdot 1_M}{M}\right), \frac{\Bar{\gamma}}{M} \sum_{k=0}^{K-1} \sum_{m=1}^M \nabla f_m\left(x_{n;k}^{(m)}, \xi_k^{(m)}\right) \right\rangle \\
    &\nonumber \quad + \frac{L}{2} \mathbb{E} \left\| -\frac{\Bar{\gamma}}{M} \sum_{k=0}^{K-1} \sum_{m=1}^M \nabla f_m\left(x_{n;k}^{(m)}, \xi_k^{(m)}\right) + \frac{C_{n} \cdot 1_M}{M} \right\|_2^2 \\
    &\nonumber \overset{(a)}{=} - \Bar{\gamma} \sum_{k=0}^{K-1}  \mathbb{E} \left\langle \nabla F\left(\frac{X_n \cdot 1_M}{M}\right), \frac{1}{M} \sum_{m=1}^M \nabla F_m\left(x_{n;k}^{(m)}\right) \right\rangle \\
    & \nonumber\quad + \frac{L}{2} \mathbb{E} \left\|\frac{1}{M} \sum_{m=1}^M C_{n}^{(m)} \right\|_2^2 \\
    &\nonumber \quad + \frac{L \Bar{\gamma}^2}{2} \mathbb{E} \left\| \frac{1}{M}  \sum_{m=1}^M \sum_{k=0}^{K-1} \nabla f_m\left(x_{n;k}^{(m)}, \xi_k^{(m)}\right)\right\|_2^2\\
    &\nonumber \overset{(b)}{\leq} - \frac{\Bar{\gamma} K}{2} \mathbb{E} \left\| \nabla F\left(\frac{X_n \cdot 1_M}{M}\right) \right\|_2^2 - \frac{\Bar{\gamma}}{2} \sum_{k=0}^{K-1} \mathbb{E} \left\|\frac{\partial F(X_{n;k}) \cdot 1_M}{M} \right\|_2^2 \\
    &\nonumber \quad + \frac{L^2 \Bar{\gamma}}{2M} \sum_{k=0}^{K-1} \sum_{m=1}^{M} \mathbb{E} \left\| \frac{X_n \cdot 1_M}{M} - x_{n;k}^{(m)}\right\|_2^2 + \frac{L}{2M^2} \sum_{m=1}^M \mathbb{E} \left\|C_{n}^{(m)} \right\|_2^2  \\
    &\quad + \frac{L \Bar{\gamma}^2}{2} \mathbb{E} \left\| \frac{1}{M}  \sum_{m=1}^M \sum_{k=0}^{K-1} \nabla f_m\left(x_{n;k}^{(m)}, \xi_k^{(m)}\right)\right\|_2^2
\end{align}
where (a) is based on the expected value of compression noise is 0, i.e., 
\begin{equation*}
    \mathbb{E}_{C_n} \left(\frac{C_n \cdot 1_M}{M}\right) = 0, 
\end{equation*}
(b) is according to L-smooth assumption, Cauchy–Schwarz inequality, $\langle \mathbf{a}, \mathbf{b} \rangle = \left(\|\mathbf{a}\|_2^2 + \|\mathbf{b}\|_2^2 - \|\mathbf{a}-\mathbf{b}\|_2^2\right)/2$, 
and 
\begin{align*}
&\quad \left\|\frac{1}{M}\sum^M_{m=1}C^{\left(m\right)}_n\right\|^2_2 = \left\|\frac{1}{M}\sum^M_{m=1}\left(Q_s\left(\Delta X^{\left(m\right)}_n\right) - \Delta^{\left(m\right)}_n\right)\right\|^2_2\\
&= \frac{1}{M^2}\sum^M_{m=1}\left\|Q_s\left(\Delta X^{\left(m\right)}_n\right) - \Delta X^{\left(m\right)}_n\right\|^2_2 = \frac{1}{M^2}\sum^M_{m=1}\left\|C^{\left(m\right)}_n\right\|^2_2. 
\end{align*}
Next, knowing that 
\begin{align*}
    & \quad \sum^{K-1}_{k=0}\sum^M_{m=1}\mathbb{E}\left\|\frac{X_n \cdot 1_M}{M} - x^{\left(m\right)}_{n;k}\right\|^2_2 \\
    & \leq \frac{4K}{1-\rho^2}\mathbb{E}\|C_n\|^2_F\\ 
    &\quad + \left[\frac{4K^2\Bar{\gamma}^2}{\left(1-\rho\right)^2} + \Bar{\gamma}^2K\left(K-1\right)\right]\sum^{K-1}_{k=0}\mathbb{E}\left\|G\left(X_{n;k}, \xi_k\right)\right\|^2_F, 
\end{align*}
and 
\begin{align*}
    & \quad \mathbb{E} \left\| \frac{1}{M} \sum_{m=1}^{M} \sum_{k=0}^{K-1} \nabla f_m\left(x_{n;k}^{\left(m\right)}, \xi_k^{\left(m\right)}\right)\right\|_2^2 \\
    &\leq \frac{2K \sigma^2}{M} + 2K \sum_{k=0}^{K-1} \mathbb{E}\left\| \frac{\partial F\left(X_{n;k}\right) \cdot 1_M}{M} \right\|_2^2, 
\end{align*}
we sum up the inequalities \ref{eq:m12} for all $n \in \{1, ..., N\}$ and have 
\begin{align*}
    & \quad F_* - F(X_1) \leq \mathbb{E} \left[F\left(\frac{X_{n+1} \cdot 1_M}{M}\right)\right] - F(X_1) \\
    & \leq -\frac{\Bar{\gamma}K}{2} \sum_{n=1}^N \mathbb{E} \left\| \nabla F\left(\frac{X_{n} \cdot 1_M}{M}\right) \right\|_2^2\\
    & \quad - \left(\frac{\Bar{\gamma}}{2} - \Bar{\gamma}^2 LK\right) \sum_{n=1}^N \sum_{k=0}^{K-1} \mathbb{E} \left\| \frac{\partial F(X_{n;k}) \cdot 1_M}{M} \right\|_2^2 \\
    & \quad + \left(\frac{2\Bar{\gamma}L^2K}{M(1-P^2)} + \frac{L}{2M^2}\right) \sum_{n=1}^N \mathbb{E} \left\| C_n \right\|_F^2 \\
    & \quad + \frac{\Bar{\gamma}L^2}{2M} \left[\frac{4K^2\Bar{\gamma}^2}{(1-\rho)^2} + \Bar{\gamma}^2 K(K-1)\right] \sum_{n=1}^N \sum_{k=0}^{K-1} \mathbb{E} \left\| G\left(X_{n;k}, \xi_k\right) \right\|_F^2\\
    & \overset{(a)}{\leq} -\frac{\Bar{\gamma}K}{2}\left(1-D’_5\right)\sum^N_{n=1}\mathbb{E}\left\|\nabla F\left(\frac{X_n \cdot 1_M}{M}\right)\right\|^2_2 \\
    & \quad - \frac{\Bar{\gamma}}{2}D_6^\prime\sum^N_{n=1}E\left\|\frac{\partial F\left(X_n \cdot 1_M\right)}{M}\right\|^2_2 \\
	& \quad + \frac{\Bar{\gamma}^2KLD’_4N}{1-8L^2\Bar{\gamma}^2KD’_1D’_3}\left[1+\frac{4L^2\Bar{\gamma}^2\left(K-1\right)}{M} \right. \\
	& \qquad \left. + \frac{8\left(K-1\right)L^2\Bar{\gamma}^2}{1-12\Bar{\gamma}^2L^2\left(K+1\right)\left(K-2\right)}\right]\sigma^2 \\
	& \quad + \frac{4\Bar{\gamma}^2K^2LD’_4N}{1-8L^2\Bar{\gamma}^2KD’_1D’_3}\left[1+\frac{4\left(K-1\right)\left(2K-1\right)L^2\Bar{\gamma}^2}{1-12\Bar{\gamma}^2L^3\left(K+1\right)\left(K-2\right)}\right]\kappa^2
\end{align*}
where (a) follows Lemma \ref{lemma:2} and Lemma \ref{lemma:3}. 
\end{proof}

\end{theorem}

To have an intuitive insight of Theorem \ref{theorem:2}, we choose a constant stepsize and obtain the following corollary: 

\begin{corollary} \label{corollary:2}
Under Theorem \ref{theorem:2}, choose the stepsize 

\begin{equation}
    \Bar{\gamma} := \left(\sigma \sqrt{N/M} + 3KL\sqrt[\leftroot{-2}\uproot{2} 3]{D_2^\prime} + 16KLD_3^\prime + 6KL \right)^{-1} 
\end{equation}

\noindent Then we have the following convergence rate: 


\begin{equation*}
\begin{split}
    & \quad \frac{1}{N} \sum_{n=1}^N \mathbb{E} \left\|\nabla F\left(\frac{X_n \cdot \bm{1}_M}{M}\right)\right\|_2^2 
     \leq \frac{16(F(\bm{x}_1) - F_*)}{K\sqrt{NM}} \\
    & \quad + \frac{4(F(\bm{x}_1) - F_*)(3L\sqrt[\leftroot{-2}\uproot{2} 3]{D_2^\prime} + 16LD_3^\prime + 6L)}{N}
\end{split}
\end{equation*}

\noindent when $N$ and $K$ satisfy  

\begin{gather*}
    N \geq \frac{L^2M}{\sigma^2} \max \left(\frac{K^2M^2(4D_3^\prime+1)^2}{(D_2^\prime)^2}, 6(K-1)(2K-1)\right) \\[1ex]
    K \leq \sqrt{\frac{[F(\bm{x}_1) - F_*]\sigma^2}{L D_2^\prime (\sigma^2 + 4 \kappa^2)}}
\end{gather*}

\end{corollary}

\begin{figure*}[htbp]
\centering
\subfigure[]{
\begin{minipage}[t]{0.24\textwidth}
\centering
\includegraphics[width=1.1\textwidth]{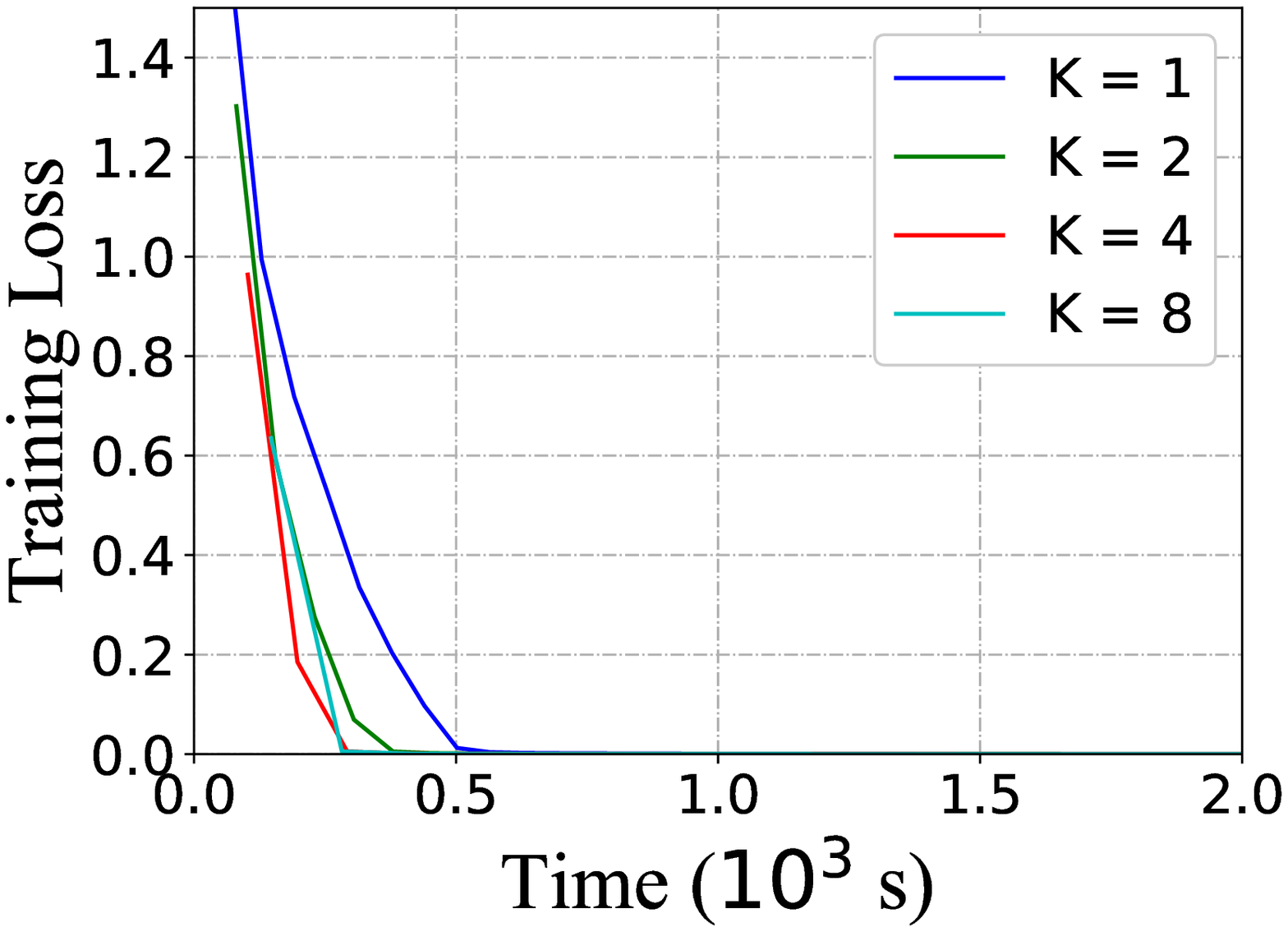}
\end{minipage}
}%
\subfigure[]{
\begin{minipage}[t]{0.24\textwidth}
\centering
\includegraphics[width=1.1\textwidth]{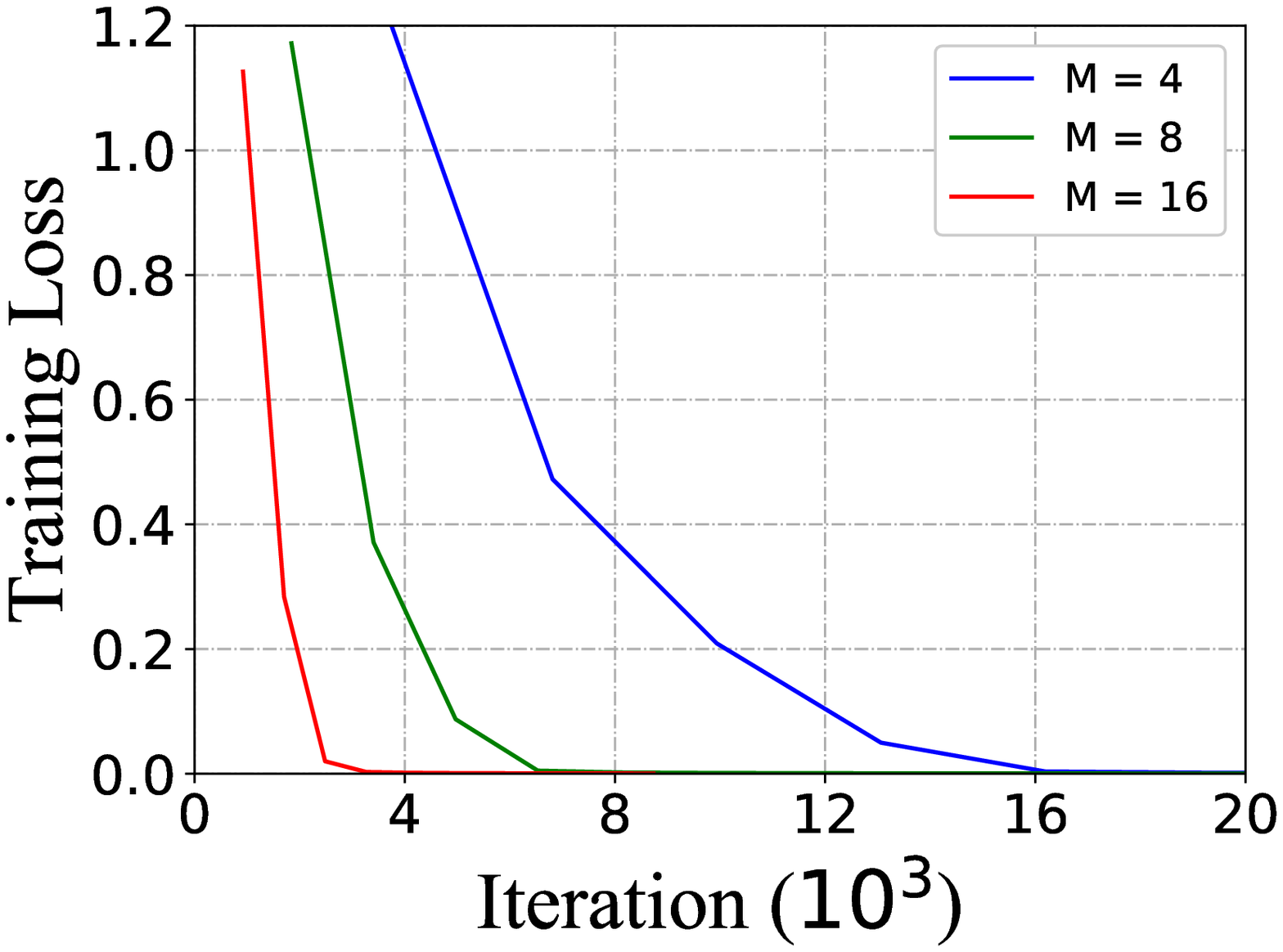}
\end{minipage}
}%
\subfigure[]{
\begin{minipage}[t]{0.24\textwidth}
\centering
\includegraphics[width=1.1\textwidth]{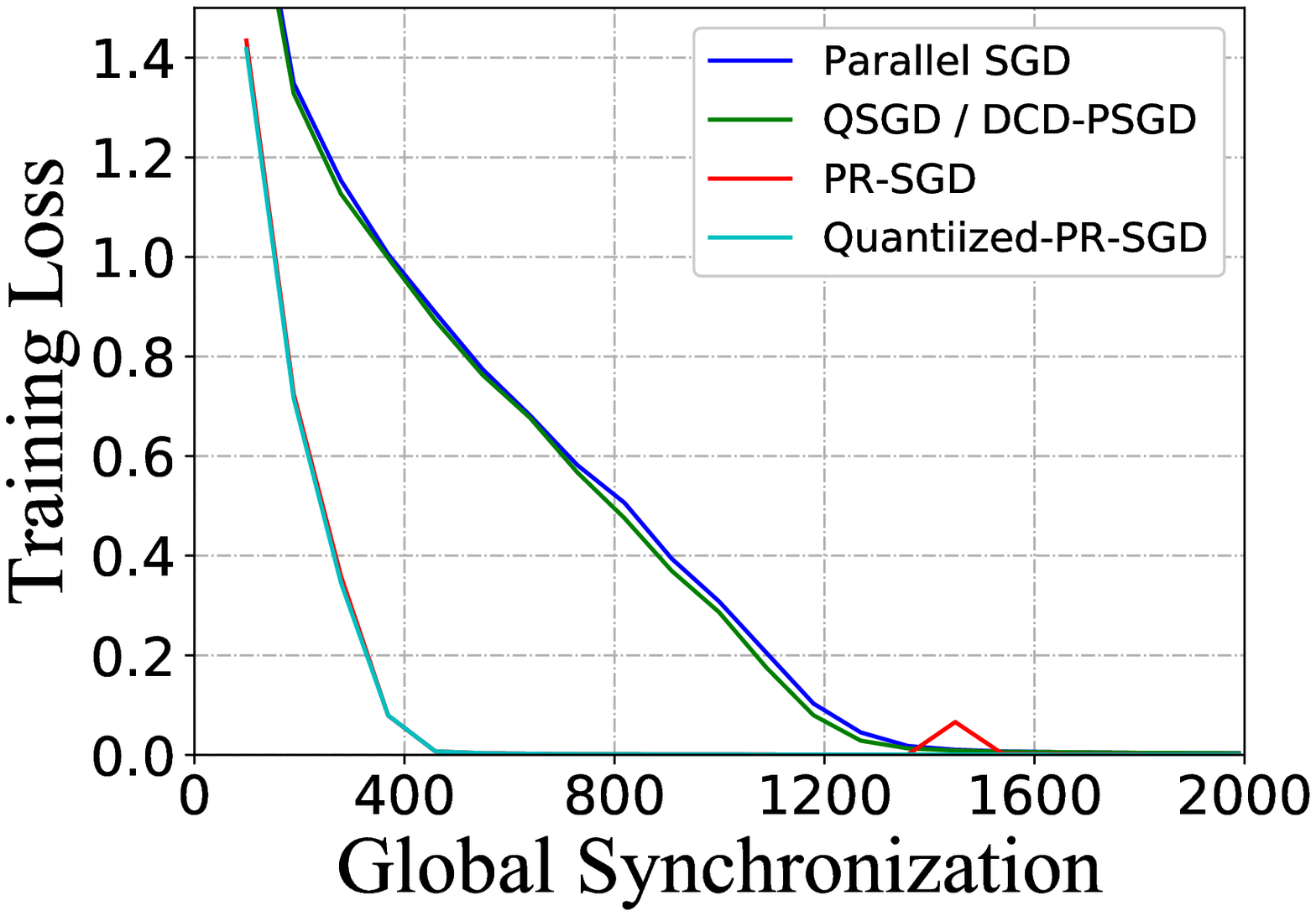}
\end{minipage}
}%
\subfigure[]{
\begin{minipage}[t]{0.24\textwidth}
\centering
\includegraphics[width=1.1\textwidth]{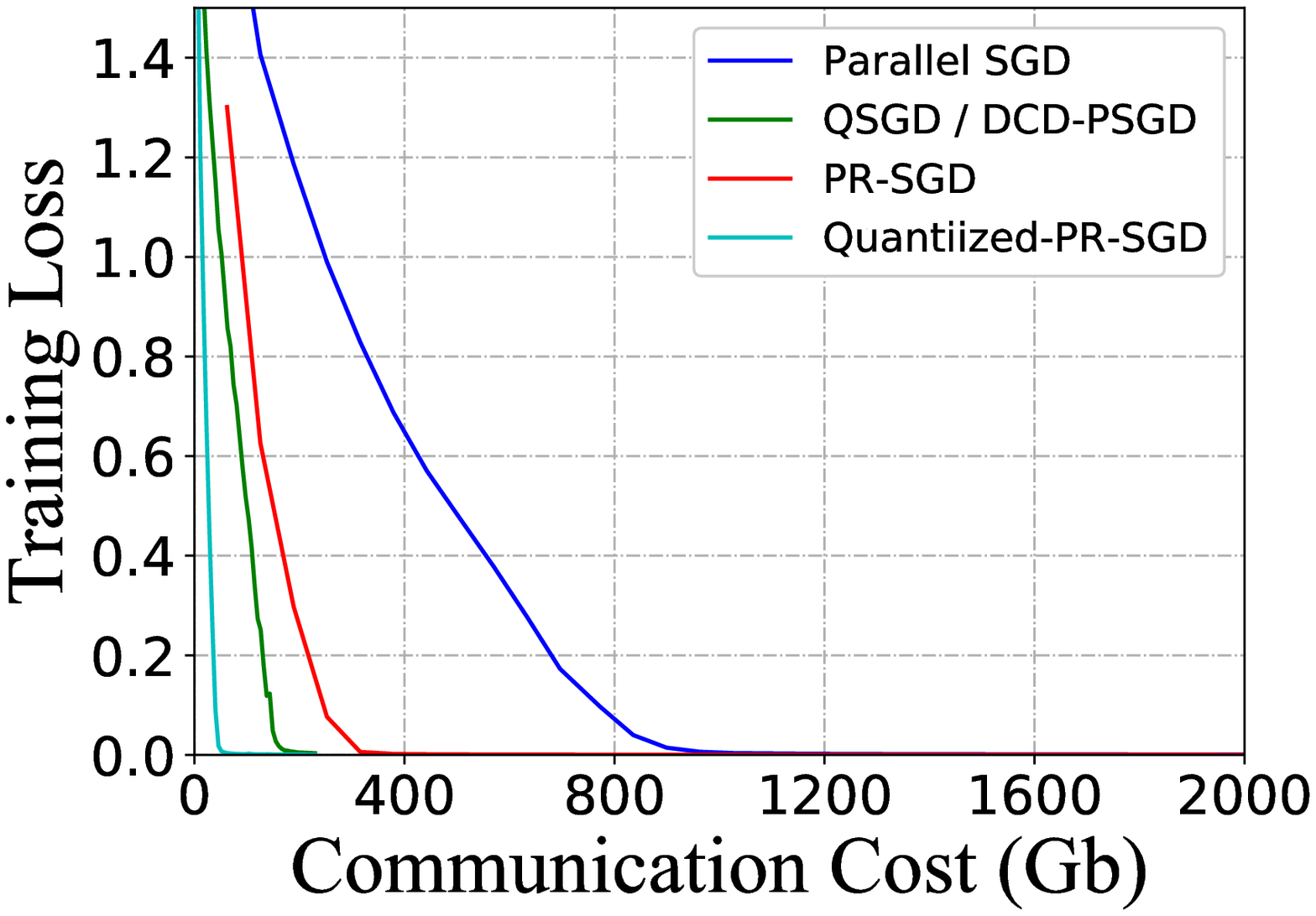}
\end{minipage}
}%
\quad
\subfigure[]{
\begin{minipage}[t]{0.24\textwidth}
\centering
\includegraphics[width=1.1\textwidth]{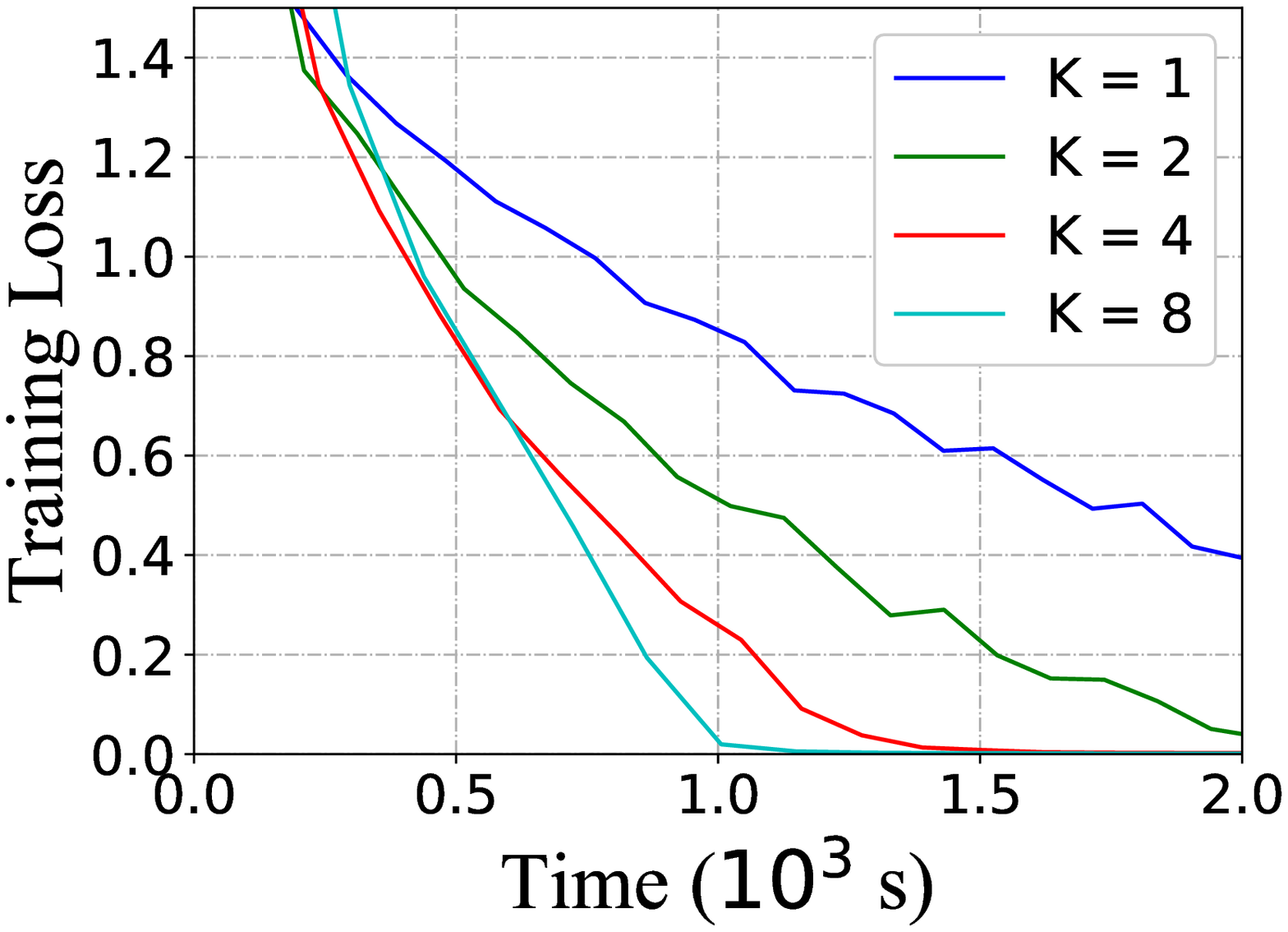}
\end{minipage}
}%
\subfigure[]{
\begin{minipage}[t]{0.24\textwidth}
\centering
\includegraphics[width=1.1\textwidth]{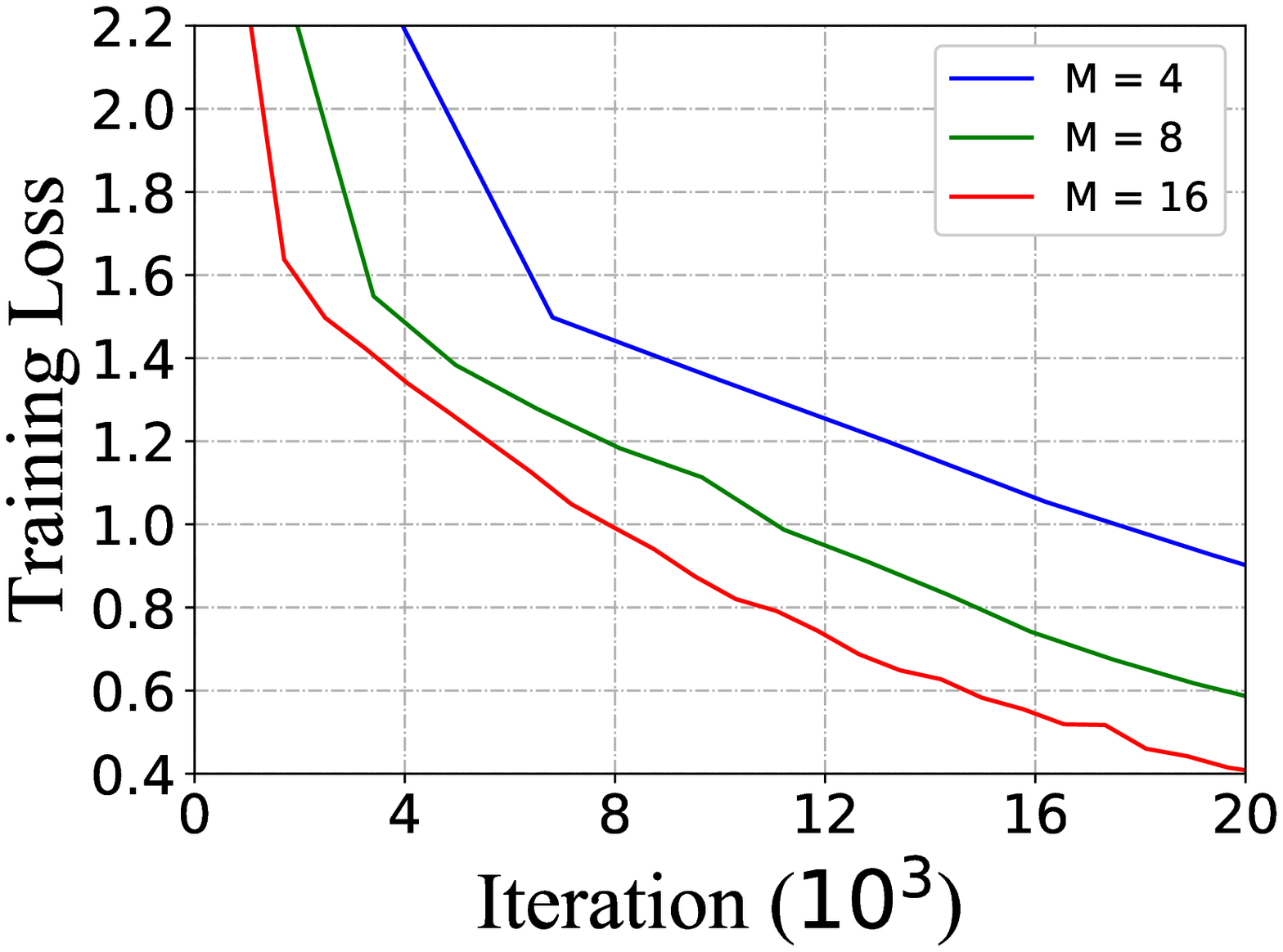}
\end{minipage}
}%
\subfigure[]{
\begin{minipage}[t]{0.24\textwidth}
\centering
\includegraphics[width=1.1\textwidth]{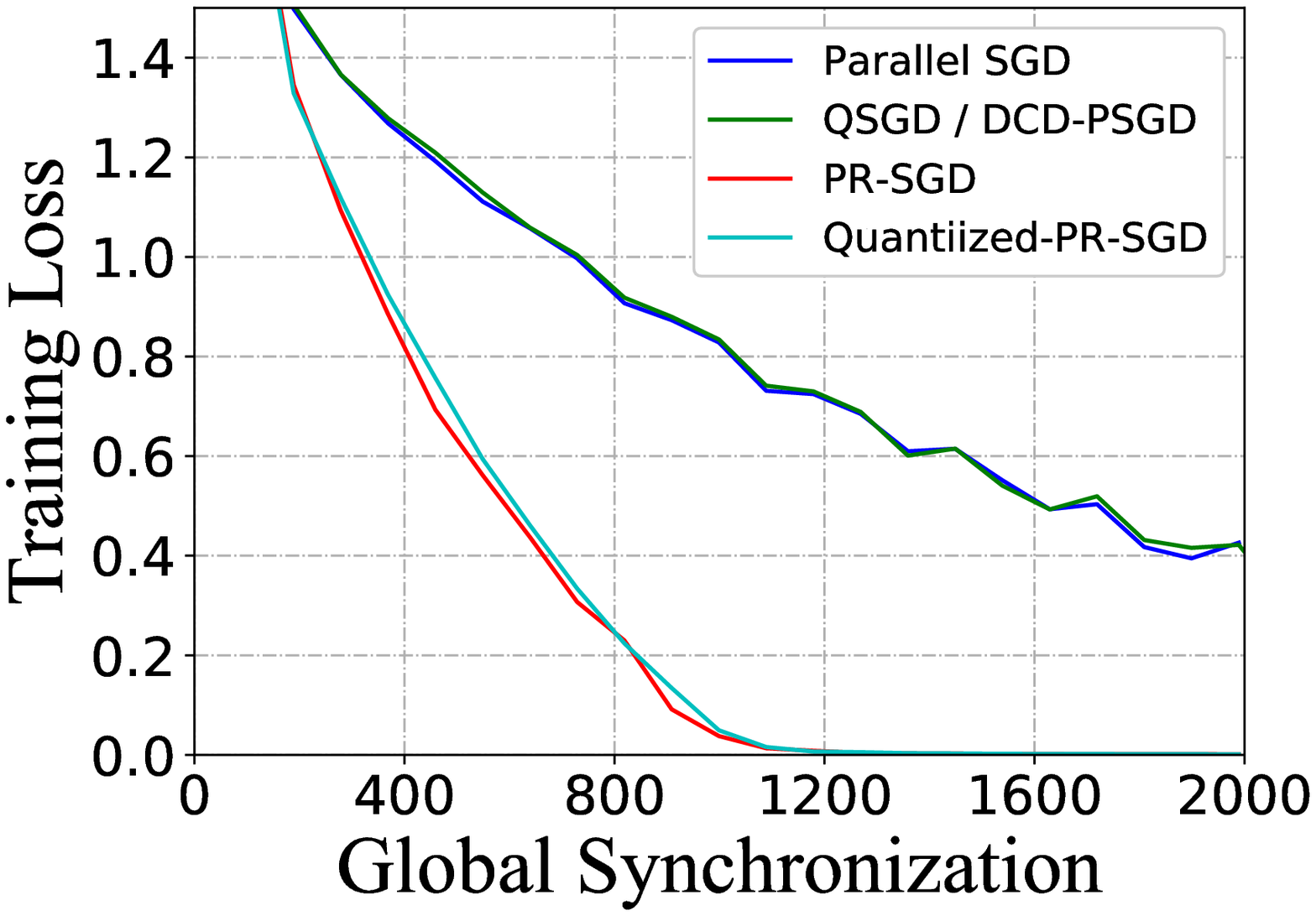}
\end{minipage}
}%
\subfigure[]{
\begin{minipage}[t]{0.24\textwidth}
\centering
\includegraphics[width=1.1\textwidth]{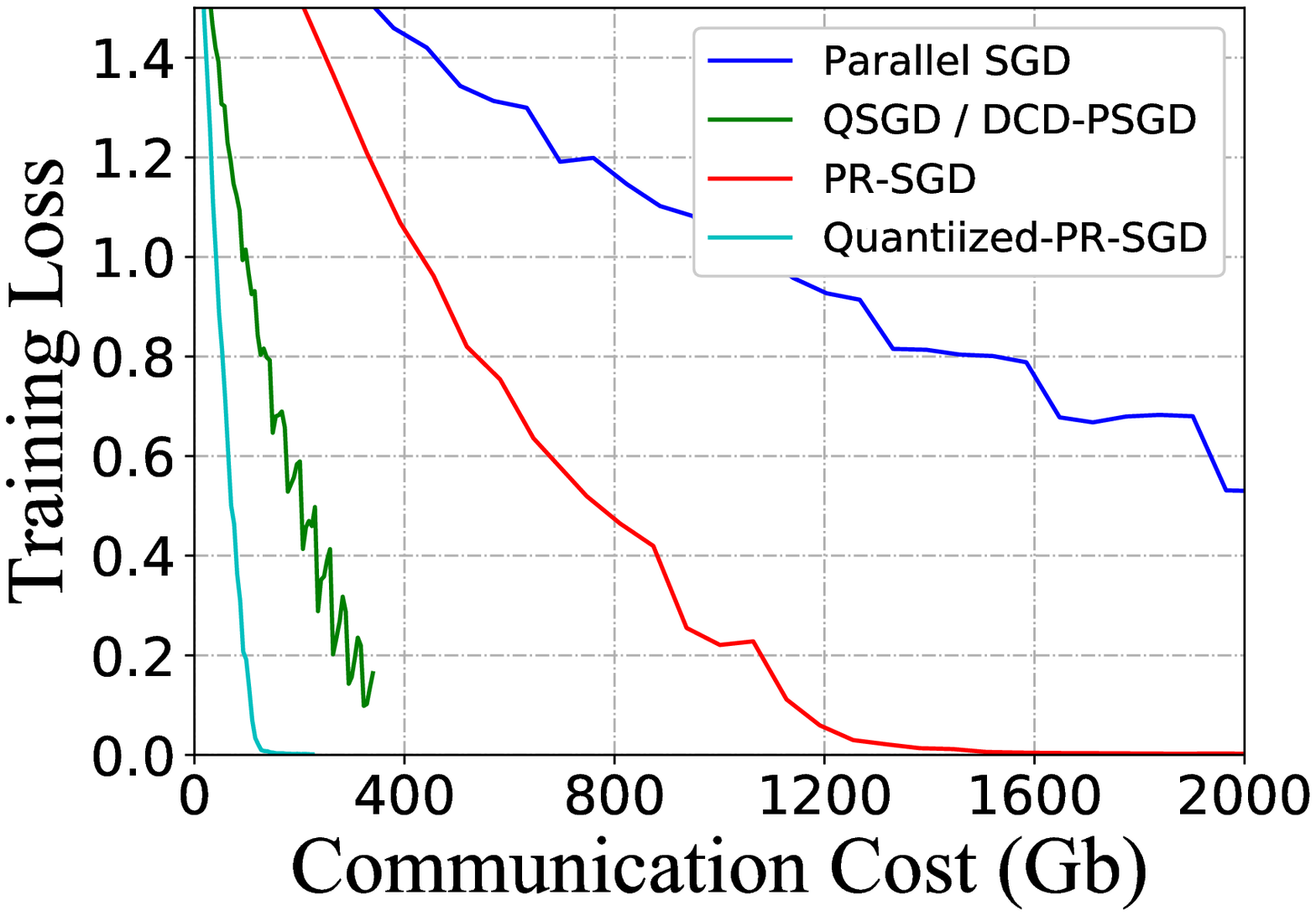}
\end{minipage}
}%
\centering
\caption{Experiment results from ResNet-34 on CIFAR-10 with RAR paradigm (top) and gossip paradigm (bottom).}
\label{fig:exp_general_case}
\end{figure*}

The result suggests that the convergent rate for Algorithm \ref{algo:G-QPRSGD} is $O(\frac{1}{K\sqrt{NM}} + \frac{1}{N})$. The following discussions interpret the tightness of our result. 

\noindent \textbf{Linear Speedup} \quad Apparently, the term $\frac{1}{K\sqrt{NM}}$ dominates the term $\frac{1}{N}$ when $N$ is sufficiently large, leading to a convergent rate of $O(1/K\sqrt{NM})$. Similar to PR-SGD \cite{local_sgd,prsgd}, the algorithm achieves linear speedup with respect to the number of workers. So far, the best case of PR-SGD is $O(1/\sqrt{NKM})$ \cite{prsgd}, indicating that the existence of gradient quantization leads to a better convergence rate. 


\noindent \textbf{Effect of K} \quad Considering that $N \times K$ is a constant, we notice that the optimal $K$ is not always 1. In fact, when $K$ is within a specified boundary, the larger value of $K$ is, the better convergence effect Algorithm \ref{algo:G-QPRSGD} has. 

\noindent \textbf{Consistence with DCD-PSGD} \quad Setting $K = 1$ to match the same scenario of DCD-PSGD. In this case, D-Quantized-PR-SGD admits the convergence rate of $O(\frac{1}{\sqrt{NM}} + \frac{1}{N})$. Apparently, the $N^{-2/3}$ term in DCD-PSGD \cite{dcdpsgd} is released, indicating that D-Quantized-PR-SGD is sightly better. 

\section{Experiments} \label{section:experiment}


To validate our analysis in the last section, we conduct a series of empirical studies for QPRSGD under RAR and gossip paradigm. We evaluate its performance for various values of $K$ and different scales of network. Besides, we compare it with the state-of-the-art algorithms (i.e. PR-SGD, QSGD and DCD-PSGD) and show the dominance of our proposed algorithm in terms of training time and communication cost. Classical PSGD is the baseline of all experiments. 

\subsection{Experimental Setup} \label{subsec:setup}

\noindent \textbf{Model and Dataset}
We train both ResNet-34 \cite{resnet} model on CIFAR-10 dataset \cite{cifar10}. ResNet-34 is composed of $16$ residual blocks with $3 \times 3$ filters and a final FC layer with a $10$-way softmax. The model size of ResNet-34 is approximately $85$ MB. CIFAR-10 is an image classification dataset constituted by $50000$ $32\times32$ colored training images and $10000$ testing images. 

\noindent \textbf{Implementation}
We implement QPRSGD under RAR and gossip paradigm on a cluster running with the operating system of GCC Linux Red Hat 4.8.5-16 and supporting OpenMPI:
\begin{itemize}
    \item \textbf{Ring network topology:} We develop a one-way ring for RAR paradigm with OpenMPI, where a node merely receives something from the last worker and sends to the next one. As for gossip paradigm, each node solely exchanges gradients with its two certain neighbors. 
\end{itemize}

We build up the environment on a job queue of the cluster, which is equipped with 32 CPU cores and 2 Tesla P100 Nvidia GPU cards for acceleration. In terms of communication compression, each element of a quantized gradient requires at most 12 bits. 

\noindent \textbf{Hyper-parameters}
Without special notation, local gradients update for 4 times (i.e. $K = 4$) before each global synchronization in both PR-SGD and Quantized-PR-SGD. Furthermore, the batch size and the learning rate are set to be $16$ and $0.1$, respectively.


\subsection{Results on General Cases} \label{subsec:epoch}

Without the constraints of network condition, Fig.~\ref{fig:exp_general_case} evaluates how the number of local updates and the number of local workers affect our convergence results under RAR paradigm and gossip paradigm. Besides, we compare QPRSGD with other well-performed algorithms by measuring the training time and the communication overhead. 


\subsubsection{Effect of $K$}

The results of RAR and gossip are shown in Fig. \ref{fig:exp_general_case} (a) and (e), respectively. As can be expected, larger $K$ probably leads to a faster and better convergence rate under both aggregation paradigms. With a larger $K$, QPRSGD takes less time to converge. For example, QPRSGD with $K=8$ takes nearly $30\%$ less time than with $K=4$ for the gossip paradigm. Besides, compared to RAR paradigm, gossip paradigm is more sensitive to the value of $K$. 

\subsubsection{Effect of $M$}

Fig. \ref{fig:exp_general_case} (b) and (f) presents the results of QRSGD with different number of nodes under RAR and gossip paradigm, respectively. When the number of worker $M$ increases, QPRSGD takes less iterations to converge. Besides, the results conform the conclusion that our algorithm preserves the linear speedup property. For example, the number of iterations reduces from $16\times10^3$ to $8\times10^3$ as the number of nodes increase from $4$ to $8$ for RAR paradigm. The similar conclusion also holds for the gossip paradigm.

\subsubsection{Training quality} 

As is shown Fig. \ref{fig:exp_general_case}(c) and (g), we evaluate the performance of Quantized-PR-SGD for training machine learning model ResNet and VGG with three different paradigms. Quantized-PR-SGD and PR-SGD convergence at a similar rate in terms of epoch, demonstrated as the overlapping lines (yellow \& red) in the chart. Similarly, QSGD and Parallel SGD (green \& blue) has a similar convergence rate in all experiments. 

The result indicates that the noise generated by quantization does not dramatically influence the convergence rate and that our approach remains the same convergence performance as PR-SGD and PSGD. It also demonstrates the epoch efficiency improvement of Quantized-PR-SGD comparing to Parallel SGD. For example, in Fig.~\ref{fig:exp_general_case}(c) the convergence speed of Quantized-PR-SGD is over 3$\times$ faster than PSGD.

\subsubsection{Communication Efficiency}
Our experiments indicate that QPRSGD has a significant reduction in communication cost comparing to other three algorithms. From Fig.~\ref{fig:exp_general_case}(d) and (h), our algorithm could achieve a communication cost reduction rate by over 90\% comparing to PR-SGD, and reduces communication cost by 50-60\% comparing to QSGD under both RAR and Gossip paradigms. In the mean time, though the communication cost is greatly reduced, our algorithm still preserves the same convergence rate as other algorithms, which is shown in Fig.~\ref{fig:exp_general_case}(c) and (g).

Under the same number of epoch, the total communication cost in PSGD is equal to the one in PR-SGD. However, a similar phenomenon does not appear by comparing QSGD and QPRSGD. We find that the total communication overhead of Quantized-PR-SGD is less than the one of QSGD. The result shows that our compression scheme requires less bits while approaching to the convergence.

\renewcommand\arraystretch{1.2}
\begin{table}[!b]
\centering
\begin{threeparttable}[t]
\begin{tabular}{|c|c|c|c|c|c|}
\hline
\multirow{3}{*}{\begin{tabular}[c]{@{}c@{}} Band.\tnote{1}\\ (Mbps)\end{tabular}} & \multirow{3}{*}{Method} & \multicolumn{2}{c|}{RAR}  & \multicolumn{2}{c|}{Gossip}   \\ \cline{3-6} 
 &    & \begin{tabular}[c]{@{}c@{}}Acc. \tnote{2}\\ (\%)\end{tabular} & \begin{tabular}[c]{@{}c@{}}Time \tnote{3}\\ (hour)\end{tabular} & \begin{tabular}[c]{@{}c@{}}Acc.\\ (\%)\end{tabular} & \begin{tabular}[c]{@{}c@{}}Time \\ (hour)\end{tabular} \\ \hline
\multirow{4}{*}{100}   

& PSGD  &  79.52   &  0.33  &  75.6   &  0.48\\ \cline{2-6} 
& QSGD  &  \textbf{79.55} &  0.36  &    74.52  &  0.51   \\ \cline{2-6} 
& PR-SGD &  79.06   &  \textbf{0.17} &  \textbf{79.15}  & \textbf{0.32} \\ \cline{2-6} 
 & QPRSGD &   79.1 & 0.21  &  75.01 &   0.42 \\ \hline
\multirow{4}{*}{20} 

& PSGD  & 70.6  &   0.63  &  56.57  & 2.82   \\ \cline{2-6} 
& QSGD & 69.27  &  0.42  & 60.52   & 2.03  \\ \cline{2-6}
& PR-SGD  & 78.84  & 0.3 &  65.02  & 0.91  \\ \cline{2-6}
& QPRSGD & \textbf{79.31} &  \textbf{0.23} &  \textbf{72.06}   &  \textbf{0.56} \\ \hline
\multirow{4}{*}{5} 

& PSGD  &  42.03 &   1.04 &  25.73 & 4.79\\ \cline{2-6} 
& QSGD   & 68.1  &  0.46  & 34.08  & 2.45 \\ \cline{2-6}
& PR-SGD  &  68.82 & 0.45 & 56.08 & 1.9   \\ \cline{2-6} 
& QPRSGD  &  \textbf{78.61}  & \textbf{0.25} & \textbf{69.05}  & \textbf{0.74} \\ \hline
\end{tabular}
\begin{tablenotes}
\item[1] Bandwidth
\item[2] Accuracy after training for 1200 seconds 
\item[3] Convergence time when the training loss $<$ 0.1
\end{tablenotes}
\end{threeparttable}
\caption{Experiment results with respect to various network bandwidths. }
\label{table:2}
\end{table}

\subsection{Results on Network-Intensive Cases} \label{subsec:convergence_time}

To better evaluate the convergence time of QPRSGD, we simulate network conditions under various bandwidths and compare their accuracy at a specific timestamp and their total time costs for convergence. In Table \ref{table:2}, we test various algorithms on ResNet model under bandwidth ranging from 5Mbps to 100Mbps. 

In ResNet model, QPRSGD shows its robustness as bandwidth decreases since it has less communication cost. Apart from the case with full bandwidth, it always has the fastest convergence speed. As for Ring AR Paradigm, the total communication cost of ResNet is relatively low considering its fast convergence rate. As a result, the PR-SGD could still converge in 60 minutes when bandwidth is limited to 5Mbps. For Gossip Paradigm, the convergence efficiency is lower, which requires more epoch to converge. As a result, the corresponding convergence time for Parallel SGD increases from 1 hour to 4 hours when bandwidth is limited to 5Mbps. However, QPRSGD remains almost the same convergence time when bandwidth changes. In general, QPRSGD runs up to 5.3x faster than PSGD under the low bandwidth environment.

\section{Conclusion} \label{section:conclusion}


This paper investigates quantized parallel restarted SGD for RAR and Gossip paradigm, which is the seamless combination of two famous techniques -- QSGD and PR-SGD. This novel SGD algorithm is analyzed from theoretical and empirical aspects. We find that the algorithm can achieve linear speedup with respect to the number of workers and the number of updates. Furthermore, it significantly saves the total communication overhead and preserves the convergence rate comparing to its original prototypes.


\appendix

\subsection{Some Lemmas for Quantization Method \cite{qsgd}}

\begin{lemma}[Unbiasedness] \label{lemma:a2}
For any vector $\bm{v} \in \mathbb{R}^d$, we have 

\begin{align*}
    \mathbb{E}[Q_s\left(\bm{v}\right)] = \bm{v}
\end{align*}

\begin{proof}
For each index $i$, 

\begin{align*}
    \nonumber\mathbb{E}[\zeta\left(v_i, s\right)/s] &= \frac{\ell}{s} \left(1 - P\right) + \frac{\ell+1}{s} P = \frac{\ell}{s} + \frac{1}{s} P \\
    &= \frac{\ell}{s} + \frac{1}{s} \left(\frac{|v_i|}{\|\bm{v}\|_2}s - \ell\right) = \frac{|v_i|}{\|\bm{v}\|_2}
\end{align*}
Thus, 
\begin{align*}
    \nonumber\mathbb{E}\left(v_i^\prime\right) &= \|\bm{v}\|_2 \cdot \text{sgn}\left(v_i\right) \cdot \mathbb{E}[\zeta\left(v_i, s\right)/s] \\
    &= \|\bm{v}\|_2 \cdot \text{sgn}\left(v_i\right) \cdot \frac{|v_i|}{\|\bm{v}\|_2} = v_i
\end{align*}
Obviously, the expected value of a quantized gradient is its original gradient. 

\end{proof}
\end{lemma}

\begin{lemma}[Second moment bound] \label{lemma:a3}
For any vector $\bm{v} \in \mathbb{R}^d$, we have 

\begin{align*}
    \mathbb{E}[\|Q_s\left(\bm{v}\right) - \bm{v}\|_2^2] \leq \frac{d}{4s^2} \|\bm{v}\|_2^2
\end{align*}

\begin{proof} 
In order to show the statement holds, we first find the expectation of the square of $\zeta\left(v_i, s\right)/s$, i.e. 

\begin{align*}
    \nonumber\mathbb{E}[\left(\zeta\left(v_i, s\right)/s\right)^2] &= \frac{\ell^2}{s^2} \left(1 - P\right) + \frac{\left(\ell + 1\right)^2}{s^2} P \\
    &= -\frac{1}{s^2} P^2 + \frac{1}{s^2} P + \frac{|v_i|^2}{\|\bm{v}\|_2^2} \\ 
    &= - \frac{1}{s^2} \left(P - \frac{1}{2}\right)^2 + \frac{1}{4s^2} + \frac{|v_i|^2}{\|\bm{v}\|_2^2} \\
    &\overset{(a)}{\leq} \frac{1}{4s^2} + \frac{|v_i|^2}{\|\bm{v}\|_2^2} 
\end{align*}
where $(a)$ holds because the range of $P$ is $[0, 1]$. Then, 
\begin{align*}
    \nonumber\mathbb{E}[\|Q_s\left(\bm{v}\right)\|_2^2] &= \sum_{i=1}^{d} \mathbb{E}[\|\bm{v}\|_2^2 \cdot \left(\zeta\left(v_i, s\right)/s\right)^2] \\
    &= \|\bm{v}\|_2^2 \cdot \sum_{i=1}^{d} \mathbb{E}[\left(\zeta\left(v_i, s\right)/s\right)^2] \\
    &\leq \|\bm{v}\|_2^2 \cdot \sum_{i=1}^{d} \left(\frac{1}{4s^2} + \frac{|v_i|^2}{\|\bm{v}\|_2^2}\right)\\
    &= \left(\frac{d}{4s^2} + 1\right)\|\bm{v}\|_2^2
\end{align*}
Therefore, with Lemma \ref{lemma:1}, we get
\begin{align*}
     \nonumber\mathbb{E}[\|Q_s\left(\bm{v}\right) - \bm{v}\|_2^2] &= \mathbb{E}[\|Q_s\left(\bm{v}\right) - \mathbb{E}[Q_s\left(\bm{v}\right)]\|_2^2]\\
     &= \mathbb{E}[\|Q_s\left(\bm{v}\right)\|_2^2] - \|\mathbb{E}[Q_s\left(\bm{v}\right)]\|_2^2 \\
     &\leq \left(\frac{d}{4s^2} + 1\right)\|\bm{v}\|_2^2 - \|\bm{v}\|_2^2 = \frac{d}{4s^2} \|\bm{v}\|_2^2
\end{align*}
\end{proof}
\end{lemma}

\subsection{Proof of Lemma \ref{lemma:1}} \label{appendix:B}

The communication framework is built with RAR paradigm. Thus, a vector $\bm{w}$ should be equally split into $M$ segments. In each communication, specific segment is compressed. For each segment $\bm{w}_i$, we have: 
\begin{equation*}
    \mathbb{E}\left[Q_s(\bm{w}_i)\right] = \bm{w}_i; \quad \mathbb{E}\|Q_s(\bm{w}_i)\|_2^2 \leq \left(\frac{d}{4s^2M}+1\right) \|\bm{w}_i\|_2^2
\end{equation*}
Denote $X$ is the compressed part and $x_i$ is the next segment. Therefore, by variance equation, we have: 
\begin{align*}
    & \quad \mathbb{E} \left\|Q_s\left(Q_s(X)+x_i\right) - \left(\mathbb{E}(X)+x_i\right) \right\|_2^2\\
    &= \mathbb{E} \left\|Q_s\left(Q_s(X)+x_i\right)\right\|_2^2 - \left\|\mathbb{E}(X)+x_i\right\|_2^2\\
    & \leq \left(\frac{d}{4s^2M}+1\right) \mathbb{E} \left\|Q_s(X)+x_i\right\|_2^2 - \left\|\mathbb{E}(X)+x_i\right\|_2^2 \\
    & = \left(\frac{d}{4s^2M}+1\right) \mathbb{E} \left\|Q_s(X)\right\|_2^2 + \frac{d}{2s^2M}\left\langle \mathbb{E}(X), x_i \right\rangle \\
    &\quad + \frac{d}{4s^2M} \|x_i\|_2^2 - \|\mathbb{E}(X)\|_2^2 \\
    & \overset{(a)}{\leq} \left(\frac{d}{4s^2M}+1\right) \mathbb{E} \left\|Q_s(X)\right\|_2^2 + \left(\frac{d}{4s^2M}-1\right)\|\mathbb{E}(X)\|_2^2 \\
    &\quad + \frac{d}{2s^2M}\|x_i\|_2^2
\end{align*}
where (a) is based on $\langle a, b \rangle \leq \left(\|a\|^2 + \|b\|^2\right)/2$. Thus, 
\begin{align*}
    &\quad\mathbb{E} \left\|Q_s\left(Q_s(X)+x_i\right)\right\|_2^2 \\
    &\leq \left(\frac{d}{4s^2M}+1\right) \mathbb{E} \left\|Q_s(X)\right\|_2^2 + \frac{d}{2s^2M}\|x_i\|_2^2 + \left\|\mathbb{E}(X)+x_i\right\|_2^2\\
    &\quad + \left(\frac{d}{4s^2M}-1\right)\|\mathbb{E}(X)\|_2^2
\end{align*}
For the sake of compression order having no influence on the final result, we have: 
\begin{align*}
    & \quad \mathbb{E} \left\| Q_s\left(...Q_s\left(Q_s\left(\bm{w}_i^{(1)}\right)+\bm{w}_i^{(2)}\right)+...+\bm{w}_i^{(M)}\right) \right\|_2^2 \\
    & \leq \left(\frac{d}{4s^2M}+1\right) \mathbb{E} \left\|Q_s\left(...Q_s\left(Q_s\left(\bm{w}_i^{(1)}\right)+\bm{w}_i^{(2)}\right)+...+\bm{w}_i^{(M-1)}\right)\right\|_2^2\\
    &\quad + \frac{d}{2s^2M}\left\|\bm{w}_i^{(M)}\right\|_2^2 + \mathbb{E}\left\|\bm{w}_i^{(1)} + ... + \bm{w}_i^{(M)}\right\|_2^2 \\
    & \quad + \left(\frac{d}{4s^2M}-1\right)\left\|\bm{w}_i^{(1)} + ... + \bm{w}_i^{(M-1)}\right\|_2^2 \\
    & \leq \left(\frac{d}{4s^2M}+1\right)^{M-1} \mathbb{E} \left\|Q_s\left(\bm{w}_i^{(1)}\right)\right\|_2^2\\
    & \quad + \frac{d}{2s^2M} \sum_{m=2}^M \left(\frac{d}{4s^2M} + 1\right)^{M-m} \left\|\bm{w}_i^{(m)}\right\|_2^2 + \left\|\sum_{m=1}^M \bm{w}_i^{(m)}\right\|_2^2 \\
    & \quad + \frac{d}{2s^2M} \sum_{m=2}^{M-1} \left(\frac{d}{4s^2M} + 1\right)^{M-m-1} \left\|\sum_{j=1}^m \bm{w}_i^{(j)}\right\|_2^2 \\
    & \leq \left[\left(\frac{d}{4s^2M}+1\right)^{M} + \left(\frac{d}{4s^2M}+1\right)^{M-2}\left(\frac{d}{4s^2M}-1\right)\right] \left\|\bm{w}_i^{(1)}\right\|_2^2 \\
    &\quad+ \frac{d}{2s^2M} \sum_{m=2}^M \left(\frac{d}{4s^2M} + 1\right)^{M-m} \left\|\bm{w}_i^{(m)}\right\|_2^2 + \left\|\sum_{m=1}^M \bm{w}_i^{(m)}\right\|_2^2 \\ 
    & \quad + \frac{d}{2s^2M} \sum_{m=2}^{M-1} \left(\frac{d}{4s^2M} + 1\right)^{M-m-1} \cdot m \sum_{j=1}^m \left\|\bm{w}_i^{(j)}\right\|_2^2 \\
    & \overset{(a)}{\leq} \left\|\sum_{m=1}^M \bm{w}_i^{(m)}\right\|_2^2 + \left[2\exp\left(\frac{d}{4s^2}\right) \right.\\
    & \left. \qquad + \frac{8s^2M}{d} \left[\left(\frac{d}{4s^2M}+1\right)^{M-1} - 1\right] - 2M\right] \sum_{m=1}^{M-1} \left\|\bm{w}_i^{(m)}\right\|_2^2 
\end{align*}
where (a) is because 
\begin{equation*}
    \begin{split}
        &\quad \frac{d}{2s^2M} \sum_{m=2}^{M-1} \left(\frac{d}{4s^2M} + 1\right)^{M-m-1} \cdot m \\
        &= \frac{8s^2M}{d} \left[\left(\frac{d}{4s^2M}+1\right)^{M-1} - 1\right] + 2 \left(\frac{d}{4s^2M}+1\right)^{M-2} - 2M
    \end{split}
\end{equation*}
and
\begin{equation*}
    \forall m \in \{2, ..., M\}, \left(\frac{d}{4s^2M}+1\right)^{M-m} \leq \left(\frac{d}{4s^2M}+1\right)^{M-2}
\end{equation*}
and
\begin{equation*}
    2\left(\frac{d}{4s^2M} + 1\right)^{M-1} \leq \left(\frac{d}{4s^2M}+1\right)^{M} + \left(\frac{d}{4s^2M}+1\right)^{M-1}; 
\end{equation*}
and function $f(M) = \left(\frac{d}{4s^2M} + 1\right)^M$ strictly decreases and thus, 
\begin{equation*}
  \left(\frac{d}{4s^2M}+1\right)^{M} \leq \lim_{M \rightarrow \infty} \left(\frac{d}{4s^2M}+1\right)^{M} = \exp \left(\frac{d}{4s^2}\right). 
\end{equation*}

\subsection{Proof of Theorem \ref{theorem:1}} \label{appendix:C}



Prior to validating Theorem \ref{theorem:1}, we introduce an important lemma below. 

\begin{lemma} \label{lemma:4}
Under Assumption \ref{assumption:1}, given that $1-2(K+1)(K-2)L^2\gamma_n^2 > 0$, we have: 

\begin{align*}
    &\quad\frac{1}{M} \sum_{m=1}^M \sum_{k=0}^{K-1} \mathbb{E} \left\| x_{n;k}^{(m)} - \Bar{x}_{n;k} \right\|_2^2\\ 
    &\leq \frac{(K-1)K \gamma_n^2 \sigma^2 (M+1)}{2M [1 - 2(K+1)(K-2)L^2\gamma_n^2]}\\
    & \quad + \frac{(K-1)K\gamma_n^2}{1 - 2(K+1)(K-2)L^2\gamma_n^2} \sum_{k=0}^{K-1} \left\|\frac{1}{M} \sum_{m=1}^M \nabla F_m \left(x_{n;k}^{(m)}\right)\right\|_2^2\\
    & \quad + \frac{2(K-1)K\gamma_n^2}{1 - 2(K+1)(K-2)L^2\gamma_n^2} \sum_{k=0}^{K-1} \left(\kappa^2 + \|\nabla F(\Bar{x}_{n;k})\|_2^2\right). 
\end{align*}

\end{lemma}

The lemma above accumulates the variance between global average parameters and local parameters throughout the local update period. According to $L$-smooth definition in Assumption \ref{assumption:1}, we have: 

\begin{align} \label{eq:11}
  &\nonumber \quad \mathbb{E} \left[F(\tilde{x}_{n+1}) - F(\tilde{x}_n)\right] \\
    &\nonumber= \mathbb{E} \left[F(\tilde{x}_{n+1}) - F(\bar{x}_{n;k})\right] + \sum_{k=0}^{K-1} \mathbb{E} [F(\bar{x}_{n;k+1}) - F(\bar{x}_{n;k})] \\
    &\nonumber \leq \mathbb{E} \left \langle \nabla F(\bar{x}_{n;k}), \tilde{x}_{n+1} - \bar{x}_{n;k} \right \rangle 
        + \frac{L}{2} \mathbb{E} \| \tilde{x}_{n+1} - \bar{x}_{n;k} \|_2^2\\
        &\nonumber\quad + \sum_{k=0}^{K-1}\left(
            \mathbb{E} \left \langle \nabla F(\bar{x}_{n;k}), \bar{x}_{n;k+1} - \bar{x}_{n;k} \right \rangle
            + \frac{L}{2} \mathbb{E} \| \bar{x}_{n;k+1} - \bar{x}_{n;k} \|_2^2
        \right) \\
    & \nonumber\overset{(a)}{=} \frac{L}{2} \mathbb{E} \left\| 
        -\frac{1}{M} Q_s\left(...Q_s\left(Q_s\left(g_n^{(1)}\right)+g_n^{(2)}\right)+...+g_n^{(M)}\right)
        + \frac{1}{M} \sum_{m=1}^{M} g_n^{(m)}
    \right\|_2^2 \\
    &\nonumber \quad + \sum_{k=0}^{K-1} 
        \mathbb{E} \left\langle \nabla F(\bar{x}_{n;k}), - \frac{\gamma_n}{M} \sum_{m=1}^{M} \nabla f_m\left(x_{n;k}^{(m)}, \xi_k^{(m)}\right) \right \rangle \\
        &\nonumber \quad + \frac{L}{2} \sum_{k=0}^{K-1} \mathbb{E} \left\| - \frac{\gamma_n}{M} \sum_{m=1}^{M} \nabla f_m \left(x_{n;k}^{(m)}, \xi_k^{(m)}\right) \right\|_2^2 \\
    &\nonumber \overset{(b)}{\leq} \frac{C_1 L \gamma_n^2}{M}  \sum_{m=1} ^ M  \underbrace{ \mathbb{E} \left\| \sum_{j=0}^{K-1} \nabla f_m(x_{n;j}^{(m)}, \xi_j^{(m)}) \right\|_2^2}_{T_1} \\
    &\nonumber \quad - \gamma_n \sum_{k=0}^{K-1} \underbrace{\mathbb{E} \left \langle \nabla 
    F\left(\bar{x}_{n;k}, \frac{1}{M} \sum_{m=1}^{M} \nabla F_m\left(x_{n;k}^{(m)}\right)\right)
    \right \rangle}_{T_2} \\
    &\quad + \frac{L\gamma_n^2}{2M^2} \sum_{k=0}^{K-1} \underbrace{ \mathbb{E} \left\| \sum_{m=1}^{M} \nabla f_m\left(x_{n;k}^{(m)}, \xi_k^{(m)}\right) \right\|_2^2}_{T_3}
\end{align}
where $(a)$ is according to unbiased feature of QSGD \cite{qsgd} and $(b)$ is based on the variance equation $\mathbb{E}[(x-\mathbb{E}(x))^2] = \mathbb{E}(x^2) - [\mathbb{E}(x)]^2$ and Lemma \ref{lemma:1}. Then, we analyze $T_1$, $T_2$ and $T_3$ one by one. 
\begin{align*}
    T_1 & = \| \sum_{j=0}^{K-1} \nabla f_m(x_{n;j}^{(m)}, \xi_j^{(m)}) \|_2^2 \\
    &= \| \sum_{j=0}^{K-1} [
            \nabla f_m(x_{n;j}^{(m)}, \xi_j^{(m)})
            - \nabla F_m(x_{n;j}^{(m)}) 
        ]
        + \sum_{j=0}^{K-1} \nabla F_m(x_{n;j}^{(m)})
    \|_2^2 \\
    & = \| \sum_{j=0}^{K-1} [
            \nabla f_m(x_{n;j}^{(m)}, \xi_j^{(m)})
            - \nabla F_m(x_{n;j}^{(m)}) 
        ] \|_2^2 \\
    &\quad + \| \sum_{j=0}^{K-1} \nabla F_m(x_{n;j}^{(m)}) \|_2^2 \\
    & = \sum_{j=0}^{K-1} \|
            \nabla f_m(x_{n;j}^{(m)}, \xi_j^{(m)})
            - \nabla F_m(x_{n;j}^{(m)}) 
        \|_2^2 \\
    &\quad + \| \sum_{j=0}^{K-1} \nabla F_m(x_{n;j}^{(m)}) \|_2^2 \\
    & \leq K \sigma^2 + K \sum_{j=0}^{K-1} \| \nabla F_m(x_{n;j}^{(m)}) \|_2^2 \\
    & \overset{(a)}{\leq} K \sigma^2 + 2 K L^2 \sum_{j=0}^{K-1} \| x_{n;j}^{(m)} - \bar{x}_{n;j} \|_2^2
        + 2 K \sum_{j=0}^{K-1} \| \nabla F_m(\bar{x}_{n;j}) \|_2^2
\end{align*}
where $(a)$ is based on Assumption \ref{assumption:1} and 
\begin{equation}
    \left\| \nabla F_m\left(x_{n;j}^{(m)}\right) \right\|_2^2 \leq 2 \left\| \nabla F_m\left(x_{n;j}^{(m)}\right) - \nabla F_m(\bar{x}_{n;j}) \right\|_2^2 + 2 \left\| \nabla F_m(\bar{x}_{n;j}) \right\|_2^2
\end{equation}
Similarly, 
\begin{align*}
    T_2 &= -\mathbb{E} \left \langle \nabla 
    F(\bar{x}_{n;k}, \frac{1}{M} \sum_{m=1}^{M} \nabla F_m(x_{n;k}^{(m)}))
    \right \rangle\\
    &\leq -\frac{1}{2} \left\| \nabla F(\bar{x}_{n;k})\right\|_2^2
    - \frac{1}{2} \left\| \frac{1}{M} \sum_{m=1}^M \nabla F_m(\bar{x}_{n;k})\right\|_2^2\\
    &\quad + \frac{L^2}{2M} \sum_{m=1}^M \left\| \bar{x}_{n;k} - x_{n;k}^{(m)}\right\|_2^2
\end{align*}
and 
\begin{align*}
    T_3 = \left\| \sum_{m=1}^{M} \nabla f_m\left(x_{n;k}^{(m)}, \xi_k^{(m)}\right) \right\|_2^2 \leq M \sigma^2 + \left\| \sum_{m=1}^M \nabla F_m\left(x_{n;k}^{(m)}\right)\right\|_2^2
\end{align*}
Next, applying $T_1$, $T_2$ and $T_3$ to Equation \ref{eq:11}, we obtain the following formula with the constraints mentioned in Theorem \ref{theorem:1}: 
\begin{align*}
    &\quad \mathbb{E} [F(\tilde{x}_{n+1}) - F(\tilde{x}_n)]\\
    & \leq -\frac{\gamma_n\varepsilon}{2}  \sum_{k=0}^{K-1} \|\nabla F(\Bar{x}_{n; k})\|_2^2 + C_1 L \gamma_n^2 K \sigma^2  + 2LK^2\gamma_n^2C_1\kappa^2\\
    &\quad + \frac{L\gamma_n^2K\sigma^2}{2M} + \left(2KL^3\gamma_n^2C_1 + \frac{L^2\gamma_n}{2}\right) \left( \frac{(K-1)K\gamma_n^2\sigma^2(M+1)}{2M\delta} \right. \\
    & \left. \qquad + \frac{2(K-1)K^2\gamma_n^2\kappa^2}{\delta} \right)
\end{align*}
Given that the stepsize is a constant value and $F_*$ is the optimizer answer that we expect to obtain, we then accumulate the formula above and get: 
\begin{align*}
    & \quad F_* - F(\tilde{x}_1) \leq \mathbb{E} [F(\tilde{x}_{N+1}) - F(\tilde{x}_1)] = \sum_{n=1}^N \mathbb{E} [F(\tilde{x}_{n+1}) - F(\tilde{x}_n)] \\
    & \leq -\frac{\Bar{\gamma}\varepsilon}{2} \sum_{n=1}^N \sum_{k=0}^{K-1} \|\nabla F(\Bar{x}_{n; k})\|_2^2 + C_1 L \Bar{\gamma}^2 K \sigma^2 N + 2LK^2\Bar{\gamma}^2C_1\kappa^2N\\ 
    & \quad + \frac{L\Bar{\gamma}^2K\sigma^2N}{2M} + \left(2KL^3\Bar{\gamma}^2C_1 + \frac{L^2\Bar{\gamma}}{2}\right) \left( \frac{(K-1)K\gamma_n^2\sigma^2(M+1)}{2M\delta} \right. \\
    & \left. \qquad + \frac{2(K-1)K^2\gamma_n^2\kappa^2}{\delta} \right) N
\end{align*}

\subsubsection{Proof of Lemma \ref{lemma:4}}

We firstly calculate for the bound of a single term $\mathbb{E} \| \Bar{x}_{n;t} - x_{n;t}^{(m)} \|_2^2$: 
\begin{align*}
    & \quad \mathbb{E} \| \Bar{x}_{n;t} - x_{n;t}^{(m)} \|_2^2 \\
    & = \mathbb{E} \| \frac{\gamma_n}{M} \sum_{m=1}^M \sum_{j=0}^{t-1} \nabla f_m(x_{n;j}^{(m)}, \xi_j^{(m)}) - \gamma_n \sum_{j=0}^{t-1} \nabla f_m(x_{n;j}^{(m)}, \xi_j^{(m)})\|_2^2 \\
    & = \mathbb{E} \| \frac{\gamma_n}{M} \sum_{m=1}^M \sum_{j=0}^{t-1} (\nabla f_m(x_{n;j}^{(m)}, \xi_j^{(m)}) - \nabla F_m(x_{n;j}^{(m)}))\\
    &\qquad + \frac{\gamma_n}{M} \sum_{m=1}^M \sum_{j=0}^{t-1} \nabla F_m(x_{n;j}^{(m)}) - \gamma_n \sum_{j=0}^{t-1} \nabla f_m(x_{n;j}^{(m)}, \xi_j^{(m)})\|_2^2  \\
    & = \mathbb{E} \| \frac{\gamma_n}{M} \sum_{m=1}^M \sum_{j=0}^{t-1} (\nabla f_m(x_{n;j}^{(m)}, \xi_j^{(m)}) - \nabla F_m(x_{n;j}^{(m)})) \|_2^2 \\
    &\quad + \mathbb{E}\|\frac{\gamma_n}{M} \sum_{m=1}^M \sum_{j=0}^{t-1} \nabla F_m(x_{n;j}^{(m)}) - \gamma_n \sum_{j=0}^{t-1} \nabla f_m(x_{n;j}^{(m)}, \xi_j^{(m)})\|_2^2 \\
    & =  \frac{\gamma_n^2}{M^2} \sum_{m=1}^M \sum_{j=0}^{t-1} \mathbb{E} \|\nabla f_m(x_{n;j}^{(m)}, \xi_j^{(m)}) - \nabla F_m(x_{n;j}^{(m)}) \|_2^2 \\
    & \quad + \mathbb{E}\|
        \frac{\gamma_n}{M} \sum_{m=1}^M \sum_{j=0}^{t-1} \nabla F_m(x_{n;j}^{(m)}) - \gamma_n \sum_{j=0}^{t-1} \nabla F_m(x_{n;j}^{(m)}) \\
        & \qquad - \gamma_n \sum_{j=0}^{t-1} [\nabla f_m(x_{n;j}^{(m)}, \xi_j^{(m)}) - \nabla F_m(x_{n;j}^{(m)})]
        \|_2^2 \\
    & =  \frac{\gamma_n^2}{M^2} \sum_{m=1}^M \sum_{j=0}^{t-1} \mathbb{E} \|\nabla f_m(x_{n;j}^{(m)}, \xi_j^{(m)}) - \nabla F_m(x_{n;j}^{(m)}) \|_2^2 \\
    & \quad + \gamma_n^2 \sum_{j=0}^{t-1} \mathbb{E} \|\nabla f_m(x_{n;j}^{(m)}, \xi_j^{(m)}) - \nabla F_m(x_{n;j}^{(m)})\|_2^2 \\
    & \quad + \mathbb{E}\|
        \frac{\gamma_n}{M} \sum_{m=1}^M \sum_{j=0}^{t-1} \nabla F_m(x_{n;j}^{(m)}) 
        - \gamma_n \sum_{j=0}^{t-1} \nabla F_m(x_{n;j}^{(m)})
        \|_2^2 \\
    & \leq \frac{\gamma_n^2}{M^2} \cdot M \cdot t \cdot \sigma^2 + \gamma_n^2 \cdot t \cdot \sigma^2 \\
    & \quad + \gamma_n^2 \mathbb{E}\|
        \frac{1}{M} \sum_{m=1}^M \sum_{j=0}^{t-1} \nabla F_m(x_{n;j}^{(m)}) 
        - \sum_{j=0}^{t-1} \nabla F_m(x_{n;j}^{(m)})
        \|_2^2 \\
    & = \gamma_n^2 \cdot t \cdot \sigma^2 (1 + \frac{1}{M}) \\ 
    &+ \gamma_n^2 \mathbb{E}\| \sum_{j=0}^{t-1} [
        \frac{1}{M} \sum_{m=1}^M  \nabla F_m(x_{n;j}^{(m)}) 
        - \nabla F_m(x_{n;j}^{(m)}) ]
        \|_2^2 \\
    & \leq (1 + \frac{1}{M}) \gamma_n^2 t \sigma^2\\
    & \quad + 
        \gamma_n^2 t \sum_{j=0}^{t-1} \mathbb{E}\|
        \frac{1}{M} \sum_{m=1}^M  \nabla F_m(x_{n;j}^{(m)}) 
        - \nabla F_m(x_{n;j}^{(m)})
        \|_2^2 \\
    & \leq (1 + \frac{1}{M}) \gamma_n^2 t \sigma^2 \\
    &\quad + \gamma_n^2 t \sum_{j=0}^{t-1} (
            2 \mathbb{E}\| \frac{1}{M} \sum_{m=1}^M  \nabla F_m(x_{n;j}^{(m)}) \|_2^2
            + 2 \mathbb{E}\| \nabla F_m(x_{n;j}^{(m)}) \|_2^2
        ) \\
    & \leq (1 + \frac{1}{M}) \gamma_n^2 t \sigma^2
        + 2 \gamma_n^2 t \sum_{j=0}^{t-1} \mathbb{E} \| \frac{1}{M} \sum_{m=1}^M \nabla F_m(x_{n;j}^{(m)}) \|_2^2 \\
    & \quad
        + 4 \gamma_n^2 t L^2 \sum_{j=0}^{t-1} \mathbb{E} \| \Bar{x}_{n;j} - x_{n;j}^{(m)} \|_2^2
        + 4 \gamma_n^2 t \sum_{j=0}^{t-1} \mathbb{E} \| \nabla F_m(\Bar{x}_{n;j}) \|_2^2
\end{align*}
Then, we sum up when t is from $0$ to $K-1$,
\begin{align*}
    & \quad \sum_{t=0}^{K-1} \mathbb{E} \| \Bar{x}_{n;t} - x_{n;t}^{(m)} \|_2^2 \\
    & \leq (1 + \frac{1}{M}) \gamma_n^2 \sigma^2 \cdot \sum_{t=0}^{K-1} t\\
    &\quad + 2 \gamma_n^2 \sum_{t=0}^{K-1} \sum_{j=0}^{t-1} t \mathbb{E} \| \frac{1}{M} \sum_{m=1}^M \nabla F_m(x_{n;j}^{(m)}) \|_2^2 \\
    & \quad
        + 4 \gamma_n^2 L^2 \sum_{t=0}^{K-1} \sum_{j=0}^{t-1} t \mathbb{E} \| \Bar{x}_{n;j} - x_{n;j}^{(m)} \|_2^2 \\
    & \quad + 4 \gamma_n^2 \sum_{t=0}^{K-1} \sum_{j=0}^{t-1} t \mathbb{E} \| \nabla F_m(\Bar{x}_{n;j}) \|_2^2 \\
    & \leq (1 + \frac{1}{M}) \gamma_n^2 \sigma^2 \frac{(K-1)K}{2} \\
    &  \quad + 2 \gamma_n^2 \frac{(K-1)K}{2} \sum_{t=0}^{K-1} \mathbb{E} \| \frac{1}{M} \sum_{m=1}^M \nabla F_m(x_{n;t}^{(m)}) \|_2^2 \\
    & \quad
        + 4 \gamma_n^2 L^2 \frac{(K+1)(K-2)}{2} \sum_{t=0}^{K-1} \mathbb{E} \| \Bar{x}_{n;j} - x_{n;j}^{(m)} \|_2^2 \\
    &\quad + 4 \gamma_n^2 \frac{(K-1)K}{2} \sum_{t=0}^{K-1} \mathbb{E} \| \nabla F_m(\Bar{x}_{n;t}) \|_2^2 \\
    & = \frac{(K-1)K \gamma_n^2 \sigma^2 (M+1)}{2M}\\
    & \quad + \gamma_n^2 (K-1)K \sum_{t=0}^{K-1} \mathbb{E} \| \frac{1}{M} \sum_{m=1}^M \nabla F_m(x_{n;t}^{(m)}) \|_2^2 \\
    & \quad
        + 2(K+1)(K-2) \gamma_n^2 L^2 \sum_{t=0}^{K-1} \mathbb{E} \| \Bar{x}_{n;j} - x_{n;j}^{(m)} \|_2^2 \\
    & \quad + 2(K-1)K \gamma_n^2 \sum_{t=0}^{K-1} \mathbb{E} \| \nabla F_m(\Bar{x}_{n;t}) \|_2^2
\end{align*}
Given that $1 - 2(K+1)(K-2)L^2 \gamma_n^2 > 0$, we have,
\begin{align*}
    & \quad \sum_{t=0}^{K-1} \mathbb{E} \| \Bar{x}_{n;t} - x_{n;t}^{(m)} \|_2^2 \\
    & \leq \frac{(K-1)K \gamma_n^2 \sigma^2 (M+1)}{2M[1 - 2(K+1)(K-2)L^2 \gamma_n^2 L^2]}\\
    & \quad + \frac{\gamma_n^2 (K-1)K}{1 - 2(K+1)(K-2)L^2 \gamma_n^2 L^2} \sum_{t=0}^{K-1} \mathbb{E} \| \frac{1}{M} \sum_{m=1}^M \nabla F_m(x_{n;t}^{(m)}) \|_2^2 \\
    & \quad
        + \frac{2(K-1)K \gamma_n^2}{1 - 2(K+1)(K-2)L^2 \gamma_n^2 L^2} \sum_{t=0}^{K-1} \mathbb{E} \| \nabla F_m(\Bar{x}_{n;t}) \|_2^2
\end{align*}

\subsection{Proof of Lemma \ref{lemma:2}}

With Lemma \ref{lemma:a2}, we attempt to bound the noise directly: 
\begin{align*}
    & \quad \sum_{n=1}^{N} \sum_{m=1}^M \mathbb{E} \| C_n^{\left(m\right)} \|_2^2 \leq \frac{d}{4s^2} \sum_{n=1}^N \sum_{m=1}^M \mathbb{E} \| \Delta X_n^{\left(m\right)} \|_2^2 \\
    & = \frac{d}{4s^2} \sum_{n=1}^N \sum_{m=1}^M \mathbb{E} \| X_n \left(W-I\right) e^{\left(m\right)} - \gamma_n \sum_{k=0}^{K-1} \nabla f_m\left(X_{n;k}^{\left(m\right)}, \xi_k^{\left(m\right)}\right)\|_2^2 \\
    & \overset{(a)}{\leq} \frac{d}{2s^2} \sum_{n=1}^N \sum_{m=1}^M \mathbb{E} \left\| X_n \left(W - I\right) e^{\left(m\right)} \right\|_2^2\\
    &\quad + \frac{d}{2s^2} \sum_{n=1}^N \sum_{m=1}^M \mathbb{E} \| \Bar{\gamma} \sum_{k=0}^{K-1} \nabla f_m\left(X_{n;k}^{\left(m\right)}, \xi_k^{\left(m\right)}\right) \|_2^2
\end{align*}
where (a) is from $\|a+b\|^2_2 \leq 2\|a\|^2_2$. Considering one of the element of the first term $\left\| X_n \left(W - I\right) e^{\left(m\right)} \right\|_2^2$, 
\begin{align*}
    & \quad \sum_{m=1}^M \mathbb{E} \| X_n\left(W - I\right) e^{\left(m\right)}\|_2^2 = \mathbb{E} \| X_n \left(W - I\right) \|_F^2 \\
    & = \mathbb{E} \| X_n P \left(\Lambda - I\right) P^T\|_F^2 \\
    & = \mathbb{E} \left\| X_n P 
    \begin{bmatrix}
    0 & & & & \\
    & \lambda_2 - 1 & & & \\
    & & \lambda_3 - 1 & & \\
    & & & \ddots & \\
    & & & & \lambda_n - 1
    \end{bmatrix}
    \right\|_F^2 \\
    & = \sum_{m=1}^M \left(\lambda_m - 1\right)^2 \mathbb{E} \| \tilde{x}_n^{\left(m\right)} \|_2^2 \\
    &\leq \mu^2 \sum_{m=2}^M \mathbb{E} \| \tilde{x}_n^{\left(m\right)} \|_2^2
    \quad \left(\mu := \max_{i \in [2, 3, \dots, M]} |\lambda_i - 1|\right)
\end{align*}
Then, the accumulative noise can be written as: 
\begin{align*}
    &\quad \sum_{n=1}^{N} \sum_{m=1}^M \mathbb{E} \| C_n^{\left(m\right)} \|_2^2 \\
    & \leq \frac{d \mu^2}{2s^2} \sum_{n=1}^N \sum_{m=2}^M \mathbb{E} \| \tilde{x}_n^{\left(m\right)}\|_2^2 \\
    & \quad + \frac{d}{2s^2} \sum_{n=1}^N \sum_{m=1}^M \mathbb{E} \| \Bar{\gamma} \sum_{k=0}^{K-1} \nabla f_m\left(x_{n;k}^{\left(m\right)}, \xi_k^{\left(m\right)}\right) \|_2^2
\end{align*}
With the recursive formula, we have the following recurrence one: 
\begin{align*}
    X_n & = X_{n-1} W - \gamma_n \sum_{k=0}^{K-1} G\left(X_{n-1;k}, \xi_k\right) + C_n\\
    &= -\sum_{\varepsilon=1}^{n-1} \Bar{\gamma} \sum_{k=0}^{K-1} G\left(X_{\varepsilon; k, \xi_k}\right) W^{n-\varepsilon-1} + \sum_{\varepsilon=1}^{n-1} C_{\varepsilon}W^{n-\varepsilon-1}
\end{align*}
Therefore, 
\begin{align*}
    &\quad\sum_{m=2}^M \mathbb{E} \| \tilde{x}_n^{\left(m\right)} \|_2^2 \\
    &= \sum_{m=2}^M \mathbb{E} \| -\sum_{\varepsilon=1}^{n-1} \Bar{\gamma} \sum_{k=0}^{K-1} G\left(X_{\varepsilon; k, \xi_k}\right) W^{n-\varepsilon-1}e^{\left(m\right)} \\
    & \qquad + \sum_{\varepsilon=1}^{n-1} C_{\varepsilon}W^{n-\varepsilon-1}e^{\left(m\right)} \|_2^2 \\
    & \leq 2 \sum_{m=2}^M \mathbb{E} \| \sum_{\varepsilon = 1}^{n-1} C_{\varepsilon}W^{n-\varepsilon-1}e^{\left(m\right)}\|_2^2\\ 
    & \quad + 2\sum_{m=2}^M \mathbb{E} \| \sum_{\varepsilon=1}^{n-1} \Bar{\gamma} \sum_{k=0}^{K-1} G\left(X_{\varepsilon; k, \xi_k}\right) W^{n-\varepsilon-1}e^{\left(m\right)} \|_2^2
\end{align*}
Next, we should consider the bound of each terms in the formula above: 
\begin{align*}
    & \quad \sum_{n-1}^N \sum_{m=2}^M \mathbb{E} \left\| \sum_{\varepsilon=1}^{n-1} C_{\varepsilon}W^{n-\varepsilon-1}e^{\left(m\right)}\right\|_2^2 \\
    & = \sum_{n=1}^N \mathbb{E} \left\| \left(\sum_{\varepsilon=1}^{n-1} C_{\varepsilon}\right) P
    \begin{bmatrix}
    0 & & & & \\
    & \lambda_2 - \varepsilon - 1 & & & \\
    & & & \ddots & \\
    & & & & \lambda_n - \varepsilon - 1
    \end{bmatrix}
    \right\|_F^2 \\
    & \leq \sum_{n=1}^N \mathbb{E} \left\| \sum_{\varepsilon=1}^{n-1} \rho^{n-\varepsilon-1} C_{\varepsilon}\right\|_F^2\\ 
    &= \frac{1}{\left(1-\rho\right)^2} \sum_{n=1}^N \| C_n \|_F^2; 
\end{align*}
and
\begin{align*}
    & \quad \sum_{n=1}^N \sum_{m=2}^M \mathbb{E} \| \sum_{\varepsilon=1}^{n-1} \Bar{\gamma} \sum_{k=0}^{K-1} G\left(X_{\varepsilon; k, \xi_k}\right) W^{n-\varepsilon-1}e^{\left(m\right)} \|_2^2\\
    & \leq \frac{1}{\left(1-\rho\right)^2} \sum_{n=1}^N \sum_{m=1}^M \left\| \Bar{\gamma} \sum_{k=0}^{K-1} \nabla f_m\left(x_{n;k}^{\left(m\right)}, \xi_k^{\left(m\right)}\right)\right\|_2^2
\end{align*}
Therefore, the accumulative noise can be represented as: 
\begin{align*}
    &\quad \sum_{n=1}^N \sum_{m=2}^M \mathbb{E} \left\| C_n^{\left(m\right)} \right\|_2^2\\
    &\leq \frac{d \mu^2}{s^2 \left(1-\rho\right)^2} \sum_{n=1}^N \sum_{m=1}^M \mathbb{E} \left\|C_n^{\left(m\right)}\right\|_2^2 \\
    & \quad + \left(\frac{d}{2s^2} + \frac{d\mu^2}{s^2\left(1-\rho\right)^2}\right)\sum_{n=1}^N \sum_{m=1}^M \left\| \Bar{\gamma} \sum_{k=0}^{K-1} \nabla f\left(x_{n;k}^{\left(m\right)}, \xi_k^{\left(m\right)}\right)\right\|_2^2 
\end{align*}
Notice that $\sum_{n=1}^N \sum_{m=2}^M \mathbb{E} \left\| C_n^{\left(m\right)} \right\|_2^2$ is in the both sides with different coefficients, we have: 
\begin{align*}
    & \quad \left[1 - \frac{d \mu^2}{s^2 \left(1-\rho\right)^2}\right] \sum_{n=1}^N \sum_{m=1}^M \mathbb{E} \left\| C_n^{\left(m\right)}\right\|_2^2\\
    & \leq \left(\frac{d}{2s^2} + \frac{d \mu^2}{s^2 \left(1-\rho\right)^2}\right) \sum_{n=1}^N \sum_{m=1}^M \left\| \Bar{\gamma} \sum_{k=0}^{K-1} \nabla f_m\left(x_{n;k}^{\left(m\right)}, \xi_k^{\left(m\right)}\right) \right\|_2^2 \\
    & \leq \Bar{\gamma}^2 K \left(\frac{d}{2s^2} + \frac{d\mu^2}{s^2\left(1-\rho\right)^2}\right)\sum_{n=1}^N \sum_{m=1}^M \sum_{k=0}^{K-1} \left\| \nabla f_m\left(x_{n;k}^{\left(m\right)}, \xi_k^{\left(m\right)}\right) \right\|_2^2 \\
    & = \Bar{\gamma}^2 K \left(\frac{d}{2s^2} + \frac{d\mu^2}{s^2\left(1-\rho\right)^2}\right)\sum_{n=1}^N \sum_{k=0}^{K-1} \| G\left(X_{n;k}, \xi_k\right) \|_F^2
\end{align*}

\noindent Thus, with $1 - \frac{d \mu^2}{s^2 \left(1-\rho\right)^2} > 0$, the theorem is proved.

\subsection{Proof of Lemma \ref{lemma:3}}

Prior to proving Lemma \ref{lemma:3}, we introduce a lemma first. 

\begin{lemma} \label{lemma:10}
Given the fixed size $\Bar{\gamma}$, under Assumption \ref{assumption:1} and Assumption \ref{assumption:2}, we have

\begin{align*}
    &\quad \sum^N_{n=1}\sum^{K-1}_{k=0} \sum^M_{m=1}E\left\|\frac{X_n\cdot 1_M}{M}-\tilde{x}^{\left(m\right)}_n\right\|^2_2\\
    &\leq  \frac{2K}{1-\rho^2}\sum^N_{n=1}E\left\|C_n\right\|^2_F + \frac{2 K^2 \Bar{\gamma}^2}{\left(1-\rho\right)^2} \sum_{n=1}^N \sum_{k=0}^{K-1} \left\| G\left(X_{n;k}, \xi_k\right) \right\|_F^2
\end{align*}
\end{lemma}

We initially establish a formula for $\|G\left(x_{n;k}, \xi_k\right)\|^2_F$: 
\begin{align*}
    &\quad \|G\left(x_{n;k}, \xi_k\right)\|^2_F = \sum^{M}_{m=1}\|\nabla f_m\left(x^{\left(m\right)}_{n;k}, \xi^{(m)}_k\right)\|^2_2 \\
    & = \sum^{M}_{m=1}\|\nabla f_m\left(x^{(m)}_{n;k}, \xi^{\left(m\right)}_k\right) - \nabla F_m\left(X^{\left(m\right)}_{n;k}\right)\|^2_2\\
    & \quad + \sum^{M}_{m=1}\|\nabla F_m\left(x^{(m)}_{n;k}\right)\|^2_2 \\
    & = M \sigma^2 + \sum^{M}_{m=1} \|\nabla F_m\left(x^{\left(m\right)}_{n;k}\right)\|^2_2 \\
    & \overset{(a)}{\leq} M \sigma^2 + 4L^2\sum^{M}_{m=1}\|x^{\left(m\right)}_{n;k} - \frac{X_{n;k} \cdot 1_M}{M}\|^2_2 + 4\kappa^2 M  \\
    & \quad + 4ML^2\|\frac{X_{n;k} \cdot 1_M}{M} - \frac{X_{n} \cdot 1_M}{M}\|^2_2 + 4\|\nabla F\left(\frac{X_{n} \cdot 1_M}{M}\right)\|^2_2 \\
    & = M \sigma^2 + 4L^2\sum^{M}_{m=1}\|x^{\left(m\right)}_{n;k} - \frac{X_{n;k} \cdot 1_M}{M}\|^2_2 + 4\kappa^2 M \\
    & \quad + 4ML^2 \Bar{\gamma}^2 \|\frac{1}{M}\sum^M_{m=1}\sum^{k-1}_{t=0}\nabla  f_m\left(X^{\left(m\right)}_{n;t}, \xi^{\left(m\right)}_t\right)\|^2_2 \\
    & \quad + 4\|\nabla F\left(\frac{X_{n} \cdot 1_M}{M}\right)\|^2_2 \\
    & \leq M \sigma^2 + 4L^2\sum^{M}_{m=1}\|x^{\left(m\right)}_{n;k} - \frac{X_{n;k} \cdot 1_M}{M}\|^2_2 + 4\kappa^2 M \\
    & \quad + 4ML^2 \Bar{\gamma}^2 \left(\frac{2k\sigma^2}{M} + 2k\sum^{k-1}_{t=0}\|\frac{\partial F\left(X_{n;t}\right)\cdot 1_M}{M}\|^2_2\right) \\
    &\quad + 4\|\nabla F\left(\frac{X_{n} \cdot 1_M}{M}\right)\|^2_2 \\
    & = M \sigma^2 + 4L^2\sum^{M}_{m=1}\|x^{\left(m\right)}_{n;k} - \frac{X_{n;k} \cdot 1_M}{M}\|^2_2 + 4\kappa^2 M + 8kL^2\sigma^2\Bar{\gamma}^2 \\
    & \quad + 8MkL^2\Bar{\gamma}\sum^{k-1}_{t=0}\|\frac{\partial F\left(X_{n;t}\right)\cdot 1_M}{M}\|^2_2 + 4\|\nabla F\left(\frac{X_{n} \cdot 1_M}{M}\right)\|^2_2
\end{align*}
where (a) follows Assumption \ref{assumption:1} and 
\begin{align*}
    &\quad \left\|\nabla F_m\left(x^{\left(m\right)}_{n;k}\right)\right\|^2_2\\
    &= \left\| \nabla F_m\left(x^{\left(m\right)}_{n;k}\right) - \nabla F_m\left(\frac{X_{n;k} \cdot 1_M}{M}\right) \right. \\
    &\left. \qquad+ \nabla F_m\left(\frac{X_{n;k} \cdot 1_M}{M}\right) - \nabla F\left(\frac{X_{n;k} \cdot 1_M}{M}\right) \right. \\
    &\left. \qquad + \nabla F\left(\frac{X_{n;k} \cdot 1_M}{M}\right) - \nabla F\left(\frac{X_n \cdot 1_M}{M}\right) + \nabla F\left(\frac{X_n \cdot 1_M}{M}\right)\right\|^2_2 \\
    &\leq 4 \left\|\nabla F_m\left(x^{\left(m\right)}_{n;k}\right) - \nabla F_m\left(\frac{X_{n;k} \cdot 1_M}{M}\right)\right\|_2^2\\
    & \quad + 4 \left\| \nabla F_m\left(\frac{X_{n;k} \cdot 1_M}{M}\right) - \nabla F\left(\frac{X_{n;k} \cdot 1_M}{M}\right) \right\|_2^2\\
    & \quad + 4 \left\| \nabla F\left(\frac{X_{n;k} \cdot 1_M}{M}\right) - \nabla F\left(\frac{X_n \cdot 1_M}{M}\right)\right\|_2^2\\
    &\quad + 4 \left\|\nabla F\left(\frac{X_n \cdot 1_M}{M}\right)\right\|_2^2. 
\end{align*}
Therefore, 
\begin{equation}
\label{equation:l13}
\begin{split}
    & \quad \sum^N_{n=1}\sum^{K-1}_{k=0}\|G\left(X_{n;k}, \xi_k\right)\|^2_F \\
    & \leq NKM\sigma^2 + 4L^2\underbrace{\sum^N_{n=1}\sum^{K-1}_{k=0}\sum^M_{m=1}\|x^{\left(m\right)}_{n;k} - \frac{X_{n;k} \cdot 1_M}{M}\|^2_2}_{T_4} + 4\kappa^2MNK \\
    & \quad + 8NL^2\sigma^2\sum^{K-1}_{k=0}k\Bar{\gamma}^2 + 8ML^2\Bar{\gamma}^2\sum^N_{n=1}\sum^{K-1}_{k=0}\sum^{k-1}_{t=0}k\|\frac{\partial F\left(X_{n;t}\right)\cdot 1_M}{M}\|^2_2 \\
    & \quad + 4K\sum^N_{n=1}\|\nabla F\left(\frac{X_{n} \cdot 1_M}{M}\right)\|^2_2
\end{split}
\end{equation}
We then find the bound for $T_4$ step by step. 
\begin{align*}
&\quad \sum^M_{m=1}\|x^{\left(m\right)}_{n;k} - \frac{X_{n;k} \cdot 1_M}{M}\|^2_2 \\
& = \sum^M_{m=1}\|\tilde{x}^{\left(m\right)}_n - \Bar{\gamma}\sum^{k-1}_{j=0}\nabla f_m\left(x^{\left(m\right)}_{n;j}, \xi^{\left(m\right)}_j\right) - \frac{X_{n} \cdot 1_M}{M}\\
&\qquad + \frac{\Bar{\gamma}}{M}\sum^{k=1}_{j=0}\sum^M_{i=1}\nabla f_i\left(x^{\left(i\right)}_{n;j}, \xi^{\left(i\right)}_j\right)\|^2_2 \\
& \overset{(a)}{\leq} 2\sum^M_{m=1}\|\tilde{x}^{\left(m\right)}_n - \frac{X_{n} \cdot 1_M}{M}\|^2_2\\
&\quad + 2\sum^M_{m=1}\|\Bar{\gamma}\sum^{k-1}_{j=0}\left(\nabla f_m\left(x^{\left(m\right)}_{n;j}, \xi^{\left(m\right)}_j\right) - \frac{1}{M}\sum^M_{i=1}\nabla f_i\left(x^{\left(i\right)}_{n;j}, \xi^{\left(i\right)}_j\right)\right)\|^2_2 \\
& \overset{(b)}{\leq} 2\sum^M_{m=1}\|\tilde{x}^{\left(m\right)}_n - \frac{X_{n} \cdot 1_M}{M}\|^2_2 + 4k\sigma^2\Bar{\gamma}^2M\\
& \quad + 4\Bar{\gamma}^2\left(6kL^2\sum^M_{m=1}\sum^{k-1}_{j=0}\|\frac{X_{n} \cdot 1_M}{M} - x^{\left(m\right)}_{n;j}\|^2_2 + 3k^2\kappa^2M\right)
\end{align*}
where (a) follows $\|a+b\|^2_2 \leq 2\|a\|^2_2 + 2\|b\|^2_2$, (b) follows the variance bound of Assumption 1. Let us consider two cases -- $k = 0$ and $k \geq 1$. When $k = 0$:
\begin{align*}
    \sum_{m=1}^M \| \tilde{x}_n^{\left(m\right)} - \frac{X_n \cdot 1_M}{M} \|_2^2 \leq 2 \sum_{m=1}^M \| \tilde{x}_n^{\left(m\right)} - \frac{X_n \cdot 1_M}{M}\|_2^2
\end{align*} When $k \geq 1$:
\begin{align*}
    & \quad \sum_{k=1}^{K-1} \sum_{m=1}^M \| x_{n;k}^{\left(m\right)} - \frac{X_{n;k} \cdot 1_M}{M} \|_2^2 \\
    & \leq 2 \sum_{k=1}^{K-1} \sum_{m=1}^M \| \tilde{x}_n^{\left(m\right)} - \frac{X_n \cdot 1_M}{M}\|_2^2 + 4 \sigma^2 \Bar{\gamma}^2 M \sum_{k=1}^{K-1}k \\
    & \quad + 24 \Bar{\gamma}^2 L^2 \sum_{k=1}^{K-1} k \sum_{m=1}^M \sum_{j=0}^{k-1} \| \frac{X_{n;j} \cdot 1_M}{M} - x_{n;j}^{\left(m\right)}\|_2^2 \\
    &\quad + 12 \Bar{\gamma}^2 \kappa^2 M \sum_{k=0}^{K-1} k^2 \\
    & \leq 2 \left(K-1\right) \sum_{m=1}^M \| \tilde{x}_n^{\left(m\right)} - \frac{X_n \cdot 1_M}{M}\|_2^2 + 2 \sigma^2 \Bar{\gamma}^2 M K\left(K-1\right)  \\
    & \quad + 24 \Bar{\gamma}^2 L^2 \frac{K\left(K-1\right)}{2} \sum_{m=1}^M \| \tilde{x}_{n}^{\left(m\right)} - \frac{X_{n} \cdot 1_M}{M} \|_2^2 \\
    & \quad + 12 \Bar{\gamma}^2 L^2 \left(K+1\right)\left(K-2\right) \sum_{k=1}^{K-1} \sum_{m=1}^M \| x_{n;k}^{\left(m\right)} - \frac{X_{n;k} \cdot 1_M}{M} \|_2^2\\
    &\quad + 2\Bar{\gamma}^2 \kappa^2 M\left(K-1\right)K\left(2K-1\right) \\
    & = \left(2K-2 + 12K\left(K-1\right)\Bar{\gamma}^2L^2\right)\sum_{m=1}^M \| \tilde{x}_n^{\left(m\right)} - \frac{X_n \cdot 1_M}{M}\|_2^2\\
    &\quad+ 2 \sigma^2 \Bar{\gamma}^2 MK\left(K-1\right) + 2\Bar{\gamma}^2 \kappa^2 M\left(K-1\right)K\left(2K-1\right)\\
    &\quad + 12\Bar{\gamma}^2L^2\left(K+1\right)\left(K-2\right) \sum_{k=1}^{K-1} \sum_{m=1}^M \| x_{n;k}^{\left(m\right)} - \frac{X_{n;k} \cdot 1_M}{M}\|_2^2
\end{align*}
Then, we have:
\begin{align*}
    & \quad [1-12\Bar{\gamma}^2L^2\left(K+1\right)\left(K-2\right)] \sum_{k=1}^{k-1} \sum_{m=1}^M \| X_{n;k}^{\left(m\right)} - \frac{x_{n;k} \cdot 1_M}{M} \|_2^2 \\
    & \leq \left(2K - 2 +12K\left(K-1\right) \Bar{\gamma}^2 L^2\right) \sum_{m=1}^M \| \tilde{x}_n^{\left(m\right)} - \frac{X_n \cdot 1_M}{M} \|_2^2 \\
    & \quad + 2 \sigma^2 \Bar{\gamma}^2 MK\left(K-1\right) + 2\Bar{\gamma}^2\kappa^2 M\left(K-1\right)K\left(2K-1\right)
\end{align*}
Hence, suppose that $1 - 12 \Bar{\gamma}^2 L^2 \left(K+1\right)\left(k-2\right) > 0$, we have: 
\begin{align*}
    & \quad \sum_{k=1}^{K-1} \sum_{m=1}^M \| x_{n;k}^{\left(m\right)} - \frac{X_{n;k} \cdot 1_M}{M}\|_2^2 \\
    & \leq \frac{2K-2+12K\left(K-1\right)\Bar{\gamma}^2L^2}{1-12\Bar{\gamma}^2L^2\left(K+1\right)\left(K-2\right)} \sum_{m=1}^M \| \tilde{x}_n^{\left(m\right)} - \frac{X_n \cdot 1_M}{M}\|_2^2\\
    &\quad + \frac{2 \sigma^2 \Bar{\gamma}^2 MK\left(K-1\right) + 2\Bar{\gamma}^2 \kappa^2 M\left(K-1\right)K\left(2K-1\right)}{1-12 \Bar{\gamma}^2 L^2 \left(K+1\right)\left(K-2\right)} \\
    & \leq \frac{2K-1+24\Bar{\gamma}^2L^2}{1-12\Bar{\gamma}^2L^2\left(K+1\right)\left(K-2\right)} \sum_{n=1}^N \sum_{m=1}^M \| \tilde{x}_n^{\left(m\right)} - \frac{X_n \cdot 1_M}{M}\|_2^2\\
    & \quad + \frac{2 \sigma^2 \Bar{\gamma}^2 MK\left(K-1\right) + 2\Bar{\gamma}^2 \kappa^2 M\left(K-1\right)K\left(2K-1\right)}{1-12 \Bar{\gamma}^2 L^2 \left(K+1\right)\left(K-2\right)}
\end{align*}
Next, we accumulate the inequality above for a total round of $N$ and obtain following inequality: 
\begin{align*}
    T_4  &= \sum_{n=1}^N \sum_{k=0}^{K-1} \sum_{m=1}^M \| x_{n;k}^{\left(m\right)} - \frac{X_{n;k} \cdot 1_M}{M}\|_2^2 \\
    & \leq \frac{2K-1+24\Bar{\gamma}^2L^2}{1-12\Bar{\gamma}^2L^2\left(K+1\right)\left(K-2\right)} \sum_{n=1}^N \sum_{m=1}^M \| \tilde{x}_n^{\left(m\right)} - \frac{X_n \cdot 1_M}{M}\|_2^2\\
    & \quad + \frac{2 \sigma^2 \Bar{\gamma}^2 MK\left(K-1\right) + 2\Bar{\gamma}^2 \kappa^2 M\left(K-1\right)K\left(2K-1\right)}{1-12 \Bar{\gamma}^2 L^2 \left(K+1\right)\left(K-2\right)} N \\
    & \overset{(a)}{\leq} \left[\frac{2D’_1D’_2}{1-\rho^2} + \frac{2D_1}{\left(1-\rho\right)^2}\right] \Bar{\gamma}^2 K \sum_{n=1}^N \sum_{k=0}^{K-1} \| G\left(X_{n;k}, \xi_k\right) \|_F^2\\
    &\quad + \frac{2 \Bar{\gamma}^2 MK\left(K-1\right)[\sigma^2 + \kappa^2\left(2K-1\right)]N}{1-12 \Bar{\gamma}^2 L^2 \left(K+1\right)\left(K-2\right)}
\end{align*}
where (a) assembles both Lemma \ref{lemma:10} and Lemma \ref{lemma:2}, and 
\begin{align*}
    D’_1 & := \frac{\left(2K+24\Bar{\gamma}^2L^2-1\right)}{1-12\Bar{\gamma}^2L^2\left(K+1\right)\left(K-2\right)};\\
    D'_2 & := \frac{d\left(1-\rho\right)^2 + 2d \mu^2}{2\left(s^2\left(1-\rho\right)^2 - d\mu^2\right)}. 
\end{align*}
Then, plugging $T_4$ back to Equation \ref{equation:l13}, we have: 
\begin{align*}
    & \quad \sum^N_{n=1}\sum^{K-1}_{k=0}\|G\left(X_{n;k}, \xi_k\right)\|^2_F \\
    & \leq NKM\sigma^2 \\ 
    &+ 4L^2\left(\left[\frac{2D_1^\prime D_2^\prime}{1-\rho^2} + \frac{2D_1^\prime}{\left(1-\rho\right)^2}\right] \Bar{\gamma}^2K\sum^N_{n=1}\sum^K_{k=0}\|G\left(X_{n;k}, \xi_k\right)\|^2_F \right. \\
    &\left. \qquad + \frac{2\Bar{\gamma}^2MK\left(K-1\right)[\sigma^2+\kappa^2\left(2K-1\right)]N}{1-12\Bar{\gamma}^2L^2\left(K+1\right)\left(K-2\right)}\right) + 4\kappa^2MNK\\
    & \quad + 4NL^2\sigma^2\Bar{\gamma}^2K\left(K-1\right) + 4K\sum^N_{n=1} \|\nabla F\left(\frac{X_n \cdot 1_M}{M}\right)\|^2_2 \\
    & \quad + 4ML^2\Bar{\gamma}^2(K-1)K\sum^N_{n=1}\sum^{K-1}_{k=0} \|\frac{\partial F\left(X_{n;t}\right) \cdot 1_M}{M}\|^2_2 
\end{align*}
Let $D’_3 = \frac{D’_2}{1-\rho^2} + \frac{1}{\left(1-\rho\right)^2}$,
\begin{align*}
    &\quad\left(1 - 8L^2\Bar{\gamma}^2KD’_1D’_3\right) \sum^N_{n=1}\sum^K_{k=0}\|G\left(X_{n;k}, \xi_k\right)\|^2_F \\
    &\leq NKM\sigma^2 + 4\kappa^2MNK + 4NL^2\sigma^2\Bar{\gamma}^2K\left(K-1\right) \\
    &\quad+ \frac{8L^2\Bar{\gamma}^2MK\left(K-1\right)[\sigma^2+\kappa^2\left(2K-1\right)]N}
    {1-12\Bar{\gamma}^2L^2\left(K+1\right)\left(K-2\right)}\\
    &\quad + 4K\sum^N_{n=1} \|\nabla F\left(\frac{X_n \cdot 1_M}{M}\right)\|^2_2
\end{align*}
\noindent Note that $1 - 8L^2\Bar{\gamma}^2KD’_1D’_3 > 0$, the lemma is proved. 

\subsubsection{Proof of Lemma \ref{lemma:10}}

We firstly find the recursive function for $\frac{X_n \cdot 1_M}{M}$ and $\tilde{x}_n^{\left(m\right)}$: 
\begin{gather*}
    \frac{X_n \cdot 1_M}{M} = - \sum_{\varepsilon=1}^{n-1}[\Bar{\gamma} \sum_{k=0}^{K-1} G\left(X_{\varepsilon; k}, \xi_k\right)] \cdot \frac{1_M}{M} + \sum_{\varepsilon=1}^{n-1} C_\varepsilon \cdot \frac{1_M}{M} \\
    \tilde{x}_n^{\left(m\right)} = -\sum_{\varepsilon=1}^{n-1} \Bar{\gamma} \sum_{k=0}^{K-1} G\left(X_{\varepsilon; k}, \xi_k\right) W^{n-\varepsilon-1} e^{\left(m\right)} + \sum_{\varepsilon=1}^{n-1} C_\varepsilon W^{n - \varepsilon - 1} e^{\left(m\right)}
\end{gather*}
Then, we should find the bound for the term $\left\| \frac{X_n \cdot 1_M}{M} - \tilde{x}_n^{\left(m\right)} \right\|_2^2$ by accumulating throughout the workers: 
\begin{align*}
    & \quad \sum_{m=1}^M \mathbb{E} \left\| \frac{X_n \cdot 1_M}{M} - \tilde{x}_n^{\left(m\right)} \right\|_2^2 \\
    & = \sum_{m=1}^M \mathbb{E} \left\| \left(\sum_{\varepsilon=1}^{n-1} C_{\varepsilon} \cdot \frac{1_M}{M} - \sum_{\varepsilon=1}^{n-1}C_{\varepsilon} W^{n-\varepsilon-1} e^{\left(m\right)}\right) \right. \\  
    & \left. \qquad - \Bar{\gamma} \sum_{\varepsilon=1}^{n-1}  \sum_{t=0}^{K-1} \left( G\left(X_{\varepsilon;t}, \xi_t\right) \cdot \frac{1_M}{M} - G\left(X_{\varepsilon;t}, \xi_t\right)  W^{n-\varepsilon-1} e^{\left(m\right)}\right)\right\|_2^2 \\
    & \overset{(a)}{\leq} 2 \sum_{m=1}^M \left( \mathbb{E} \| \sum_{\varepsilon=1}^{n-1} C_{\varepsilon} \cdot \left(\frac{1_M}{M} - W^{n-\varepsilon-1} e^{\left(m\right)}\right) \|_2^2 \right.\\
    & \qquad \left. + \mathbb{E} \| \Bar{\gamma} \sum_{\varepsilon=1}^{n-1}  \sum_{t=0}^{K-1}  G\left(X_{\varepsilon;t}, \xi_t\right) \cdot \left(\frac{1_M}{M} - W^{n-\varepsilon-1} e^{\left(m\right)}\right)\|_2^2\right) \\
    & \overset{(b)}{\leq} 2 \sum_{m=1}^M  \sum_{\varepsilon=1}^{n-1} \mathbb{E} \| C_{\varepsilon} \cdot \left(\frac{1_M}{M} - W^{n-\varepsilon-1} e^{\left(m\right)}\right)\|_2^2 \\
    & \quad + 2 \sum_{m=1}^M \mathbb{E} \| \Bar{\gamma} \sum_{\varepsilon=1}^{n-1}  \sum_{t=0}^{K-1} G\left(X_{\varepsilon;t}, \xi_t\right) \cdot \left(\frac{1_M}{M} - W^{n-\varepsilon-1} e^{\left(m\right)}\right)\|_2^2 \\
    & \overset{(c)}{\leq} 2 \sum_{\varepsilon=1}^{n-1} \mathbb{E} \| \rho^{n - \varepsilon - 1} C_\varepsilon \|_F^2\\
    &\quad + 2 \mathbb{E} \left(\Bar{\gamma} \sum_{\varepsilon=1}^{n-1} \sum_{t=0}^{K-1} \rho^{n-\varepsilon-1} \| G\left(X_{\varepsilon; t}, \xi_t\right)\|_F\right)^2
\end{align*}
where (a) follows $\|a-b\|^2_2 \leq 2\|a\|^2_2 + 2\|b\|^2_2$, (b) is according to the fact that the expected compression loss is 0 and all nodes work independently, (c) is based on Lemma \ref{lemma:9}. With Lemma \ref{lemma:10}, we have: 
\begin{equation*}
    \sum_{n=1}^N \sum_{\varepsilon=1}^{n-1} \mathbb{E} \| \rho^{n-\varepsilon-1} C_\varepsilon\|_F^2 \leq \frac{1}{1-\rho^2} \sum_{n=1}^N \mathbb{E} \| C_n \|_F^2
\end{equation*}
and
\begin{align*}
    &\quad \sum_{n=1}^M \mathbb{E}\left(\Bar{\gamma} \sum_{\varepsilon=1}^{n-1} \sum_{t=0}^{K-1} \rho^{n-\varepsilon-1} \| G\left(X_{\varepsilon;t}, \xi_t\right)\|_F\right)^2\\
    & \leq \frac{1}{\left(1-\rho\right)^2} \sum_{n=1}^N \left(\Bar{\gamma} \sum_{t=0}^{K-1} \| G\left(X_{n;t}, \xi_t\right)\|_F \right)^2\\
    & \leq \frac{\Bar{\gamma}^2K}{\left(1-\rho\right)^2} \sum_{n=1}^N \sum_{t=0}^{K-1} \| G\left(X_{n;t}, \xi_t\right)\|_F^2
\end{align*}
Therefore, we have 
\begin{align*}
     & \quad \sum_{n=1}^N \sum_{k=0}^{K-1} \sum_{m=1}^M \mathbb{E} \| \frac{X_n \cdot 1_M}{M} - \tilde{X}_n^{\left(m\right)} \|_2^2 \\
     &\leq 2 \sum_{n=1}^N \sum_{k=0}^{K-1} \sum_{\varepsilon=1}^{n-1} \mathbb{E} \| \rho^{n - \varepsilon - 1} C_\varepsilon \|_F^2 \\
     &\quad + 2 \sum_{n=1}^N \sum_{k=0}^{K-1} \mathbb{E} \left(\Bar{\gamma} \sum_{\varepsilon=1}^{n-1} \sum_{t=0}^{K-1} \rho^{n-\varepsilon-1} \| G\left(X_{\varepsilon; t}, \xi_t\right) \|_F\right)^2 \\
    & \leq \frac{2K}{1 - \rho^2} \sum_{n=1}^N \mathbb{E} \| C_n \|_F^2 + \frac{2 K^2 \Bar{\gamma}^2}{\left(1-\rho\right)^2} \sum_{n=1}^N \sum_{k=0}^{K-1} \| G\left(X_{n;k}, \xi_k\right) \|_F^2
\end{align*}

\subsection{Proof of Theorem \ref{theorem:2}}


\begin{lemma}[Lemma 5 in \cite{dcdpsgd}] \label{lemma:8}
For any matrix $X_t \in \mathbb{R}^{N\times n}$, decompose the confusion matrix $W$ as $W = \sum_{i=1}^n \lambda_i \bm{v}^{\left(i\right)} \left(\bm{v}^{T}\right) = P \Lambda P^T$, where $P = \left(\bm{v}^{\left(1\right)}, \bm{v}^{\left(1\right)}, \dots, \bm{v}^{\left(n\right)}\right) \in \mathbb{R}^{N \times n}$, $\bm{v}^{\left(i\right)}$ is the normalized eigenvector of $\lambda_i$ and $\Lambda$ is a diagnal matrix with $\lambda_i$ be its i-th element. We have
\begin{align*}
    \sum_{m=1}^M \mathbb{E} \| X_n \cdot W^t \cdot e^{\left(m\right)} - X_n \cdot \frac{1_M}{M}\|_2^2 \leq \|  \rho^{2t} X_n \|_F^2 = \rho^{2t} \| X_n \|_F^2
\end{align*}
\end{lemma}

\begin{lemma} [Lemma 6 in \cite{dcdpsgd}] \label{lemma:9}
Given two non-negative sequences $\{a_t\}^{\infty}_{t=1}$ and $\{b_t\}^{\infty}_{t=1}$ that satisfying 
\begin{align*}
    a_{\varepsilon} = \sum_{\eta=1}^{\varepsilon} \rho^{\varepsilon - \eta} b_\eta
\end{align*}
\noindent with $\rho \in [0,1)$, we have
\begin{align*}
    \sum_{\varepsilon=1}^{k} a_{\varepsilon} \leq \frac{1}{1-\rho} \sum_{\varepsilon=1}^k b_{\varepsilon}; \quad 
    \sum_{\varepsilon=1}^{k} a_{\varepsilon}^2 \leq \frac{1}{\left(1-\rho\right)^2} \sum_{\varepsilon=1}^k b_{\varepsilon}^2
\end{align*}
\end{lemma}

\begin{lemma} \label{lemma:11}
Given the fixed stepsize $\Bar{\gamma}$, under Assumption \ref{assumption:1} and Assumption \ref{assumption:2}, we have:
\begin{align*}
    & \quad \sum^N_{n=1}\sum^{K-1}_{k=0}\sum^M_{m=1}E\|\frac{X_n \cdot 1_M}{M} - x^{\left(m\right)}_{n;k}\|^2_2 \\
    & \leq \frac{4K}{1-\rho^2}\sum^N_{n=1}E\|C_n\|^2_F\\ 
    &\quad + \left[\frac{4K^2\Bar{\gamma}^2}{\left(1-\rho\right)^2} + \Bar{\gamma}^2K\left(K-1\right)\right]\sum^N_{n=1}\sum^{K-1}_{k=0}\|G\left(X_{n;k}, \xi_k\right)\|^2_F
\end{align*}
\end{lemma}

\begin{lemma} \label{lemma:14}

Under Assumption \ref{assumption:1} and Assumption \ref{assumption:2}, given the constant stepsize $\Bar{\gamma}$, we have: 
\begin{align*}
    & \quad \mathbb{E} \| \frac{1}{M} \sum_{m=1}^{M} \sum_{k=0}^{K-1} \nabla f_m\left(x_{n;k}^{\left(m\right)}, \xi_k^{\left(m\right)}\right)\|_2^2 \\
    &\leq \frac{2K \sigma^2}{M} + 2K \sum_{k=0}^{K-1} \| \frac{\partial F\left(X_{n;k}\right) \cdot 1_M}{M} \|_2^2
\end{align*}
\end{lemma}

By L-smooth feature in Assumption \ref{assumption:1}, we have: 
\begin{align*}
    & \quad \mathbb{E} [F(\frac{X_{n+1} \cdot 1_M}{M})] - \mathbb{E} [F(\frac{X_{n} \cdot 1_M}{M})] \\
    & \leq \mathbb{E} \langle \nabla F(\frac{X_{n} \cdot 1_M}{M}, \frac{X_{n+1} \cdot 1_M}{M} - \frac{X_{n} \cdot 1_M}{M}) \rangle \\
    & \quad + \frac{L}{2} \|\frac{X_{n+1} \cdot 1_M}{M} - \frac{X_{n} \cdot 1_M}{M}\|_2^2 \\
    & = \mathbb{E} \langle \nabla F(\frac{X_{n} \cdot 1_M}{M}), -\frac{\Bar{\gamma} \sum_{k=0}^{K-1} \sum_{m=1}^M \nabla f_m(x_{n;k}^{(m)}, \xi_k^{(m)})}{M} + \frac{C_{n} \cdot 1_M}{M} \rangle \\
    & \quad + \frac{L}{2} \| -\frac{\Bar{\gamma} \sum_{k=0}^{K-1} \sum_{m=1}^M \nabla f_m(x_{n;k}^{(m)}, \xi_k^{(m)})}{M} + \frac{C_{n} \cdot 1_M}{M} \|_2^2 \\
    & \overset{(a)}{=} - \frac{\Bar{\gamma}}{M} \sum_{k=0}^{K-1} \sum_{m=1}^M \mathbb{E} \langle \nabla F(\frac{X_n \cdot 1_M}{M}), \nabla F_m(x_{n;k}^{(m)}) \rangle \\
    & \quad + \frac{L}{2} \| -\frac{\Bar{\gamma} \sum_{k=0}^{K-1} \sum_{m=1}^M \nabla f_m(x_{n;k}^{(m)}, \xi_k^{(m)})}{M} + \frac{C_{n} \cdot 1_M}{M} \|_2^2 \\
    & \leq - \frac{\Bar{\gamma} K}{2} \mathbb{E} \| \nabla F(\frac{X_n \cdot 1_M}{M}) \|_2^2 - \frac{\Bar{\gamma}}{2} \sum_{k=0}^{K-1} \mathbb{E} \| \frac{\partial F(X_{n;k} \cdot 1_M)}{M} \|_2^2 \\
    & \quad + \frac{\Bar{\gamma}}{2M} \sum_{k=0}^{K-1} \sum_{m=1}^{M} \mathbb{E} \| \nabla F_m(\frac{X_n \cdot 1_M}{M}) - \nabla F_m(x_{n;k}^{(m)})\|_2^2\\
    & \quad + \frac{L \Bar{\gamma}^2}{2} \| \frac{1}{M}  \sum_{m=1}^M \sum_{k=0}^{K-1} \nabla f_m(x_{n;k}^{(m)}, \xi_k^{(m)})\|_2^2 + \frac{L}{2}\| \frac{1}{M} \sum_{m=1}^M C_{n}^{(m)} \|_2^2 \\
    & \overset{(b)}{\leq} - \frac{\Bar{\gamma} K}{2} \mathbb{E} \| \nabla F(\frac{X_n \cdot 1_M}{M}) \|_2^2 - \frac{\Bar{\gamma}}{2} \sum_{k=0}^{K-1} \mathbb{E} \| \frac{\partial F(X_{n;k} \cdot 1_M)}{M} \|_2^2 \\
    & \quad + \frac{\Bar{\gamma} L^2}{2M} \sum_{k=0}^{K-1} \sum_{m=1}^{M} \mathbb{E} \| \frac{X_n \cdot 1_M}{M} - x_{n;k}^{(m)}\|_2^2\\
    &\quad + \frac{L \Bar{\gamma}^2}{2} \| \frac{1}{M}  \sum_{m=1}^M \sum_{k=0}^{K-1} \nabla f_m(X_{n;k}^{(m)}, \xi_k^{(m)})\|_2^2 + \frac{L}{2M^2} \sum_{m=1}^M\|   C_{n}^{(m)} \|_2^2 
\end{align*}
where (a) is based on the expected value of compression noise is 0, (b) follows that 
\begin{align*}
&\quad \|\frac{1}{M}\sum^M_{m=1}C^{\left(m\right)}_n\|^2_2 = \|\frac{1}{M}\sum^M_{m=1}\left(Q_s\left(\Delta X^{\left(m\right)}_n\right) - \Delta^{\left(m\right)}_n\right)\|^2_2\\
&= \frac{1}{M^2}\sum^M_{m=1}\|Q_s\left(\Delta X^{\left(m\right)}_n\right) - \Delta X^{\left(m\right)}_n\|^2_2 = \frac{1}{M^2}\sum^M_{m=1}\|C^{\left(m\right)}_n\|^2_2
\end{align*}
Next, we should find the accumulative bound: 
\begin{align*}
    & \quad F_* - F(X_1) \leq \mathbb{E} [F(\frac{X_{n+1} \cdot 1_M}{M})] - F(X_1) \\
    & \leq -\frac{\Bar{\gamma}K}{2} \sum_{n=1}^N \| \nabla F(\frac{X_{n} \cdot 1_M}{M}) \|_2^2 - \frac{\Bar{\gamma}}{2} \sum_{n=1}^N \sum_{k=0}^{K-1} \| \frac{\partial F(X_{n;k}) \cdot 1_M}{M} \|_2^2 \\
    & \quad + \frac{\Bar{\gamma} L^2}{2M} \sum_{n=1}^N  \sum_{k=0}^{K-1} \sum_{m=1}^{M} \mathbb{E} \| \frac{X_n \cdot 1_M}{M} - x_{n;k}^{(m)}\|_2^2\\
    &\quad + \frac{L \Bar{\gamma}^2}{2} \sum_{n=1}^N \| \frac{1}{M}  \sum_{m=1}^M \sum_{k=0}^{K-1} \nabla f_m(x_{n;k}^{(m)}, \xi_k^{(m)})\|_2^2 + \frac{L}{2M^2} \sum_{n=1}^N  \sum_{m=1}^M\|   C_{n}^{(m)} \|_2^2 \\
    & \overset{(a)}{\leq} -\frac{\Bar{\gamma}K}{2} \sum_{n=1}^N \| \nabla F(\frac{X_{n} \cdot 1_M}{M}) \|_2^2 - \frac{\Bar{\gamma}}{2} \sum_{n=1}^N \sum_{k=0}^{K-1} \| \frac{\partial F(X_{n;k}) \cdot 1_M}{M} \|_2^2 \\
    & \quad + \frac{2\Bar{\gamma}L^2K}{M(1-P^2)} \sum_{n=1}^N \mathbb{E} \|C_n\|_F^2 + \frac{L}{2M^2} \sum_{n=1}^N\|C_{n}^{(m)} \|_F^2\\ 
    & \quad + \frac{\Bar{\gamma}L^2}{2M} (\frac{4K^2 \Bar{\gamma}^2}{(1-\rho)^2} + \Bar{\gamma}^2 K(K-1)) \sum_{n=1}^N \sum_{k=0}^{K-1} \| G(X_{n;k}, \xi_k) \|_F^2 \\
    & \quad + \frac{\Bar{\gamma^2 L}}{2} (\frac{2K\sigma^2 N}{M} + 2K \sum_{n=1}^N \sum_{k=0}^{K-1} \| \frac{\partial F(X_{n;k}) \cdot 1_M}{M}\|_2^2) \\
    & = -\frac{\Bar{\gamma}K}{2} \sum_{n=1}^N \| \nabla F(\frac{X_{n} \cdot 1_M}{M}) \|_2^2\\
    & \quad - (\frac{\Bar{\gamma}}{2} - \Bar{\gamma}^2 LK) \sum_{n=1}^N \sum_{k=0}^{K-1} \| \frac{\partial F(X_{n;k}) \cdot 1_M}{M} \|_2^2 \\
    & \quad + (\frac{2\Bar{\gamma}L^2K}{M(1-P^2)} + \frac{L}{2M^2}) \sum_{n=1}^N \| C_n \|_F^2 \\
    & \quad + \frac{\Bar{\gamma}L^2}{2M} [\frac{4K^2\Bar{\gamma}^2}{(1-\rho)^2} + \Bar{\gamma}^2 K(K-1)] \sum_{n=1}^N \sum_{k=0}^{K-1} \| G(X_{n;k}, \xi_k) \|_F^2\\
    & \overset{(b)}{\leq} -\frac{\Bar{\gamma}K}{2} \sum_{n=1}^N \| \nabla F(\frac{X_{n} \cdot 1_M}{M}) \|_2^2\\
    & \quad - (\frac{\Bar{\gamma}}{2} - \Bar{\gamma}^2 LK) \sum_{n=1}^N \sum_{k=0}^{K-1} \| \frac{\partial F(X_{n;k}) \cdot 1_M}{M} \|_2^2 \\
    & \quad + \frac{\Bar{\gamma}L^2}{2M} [\frac{4K^2\Bar{\gamma}^2}{(1-\rho)^2} + \Bar{\gamma}^2 K(K-1)] \sum_{n=1}^N \sum_{k=0}^{K-1} \| G(X_{n;k}, \xi_k) \|_F^2 \\
    & \quad + (\frac{2\Bar{\gamma}L^2K}{M(1-\rho^2)} + \frac{L}{2M^2}) \Bar{\gamma}^2 K D_2^\prime \sum_{n=1}^N \sum_{k=0}^{K-1} \| G(X_{n;k}, \xi_k)\|_F^2 \\
    & \overset{(c)}{\leq} -\frac{\Bar{\gamma}K}{2} \sum_{n=1}^N \| \nabla F(\frac{X_{n} \cdot 1_M}{M}) \|_2^2\\ 
    &\quad - (\frac{\Bar{\gamma}}{2} - \Bar{\gamma}^2 LK) \sum_{n=1}^N \sum_{k=0}^{K-1} \| \frac{\partial F(X_{n;k}) \cdot 1_M}{M} \|_2^2 \\
    & \quad + \left(\frac{\Bar{\gamma}^3L^2K}{2M}[\frac{4K}{(1-\rho)^2} + K - 1] + \frac{\Bar{\gamma}^2 KL}{M}(\frac{2 \Bar{\gamma}LK}{1-\rho^2} + \frac{1}{2M})D_2 \right) \\
    & \qquad \cdot \left(\frac{MNK(\sigma^2 + 4\kappa^2) + 4NL^2\sigma^2 \Bar{\gamma}^2K(K-1)}{1-8L^2\Bar{\gamma}^2KD’_1D’_3}\right) \\
    & \quad + \frac{8MNK(K-1)L^2\Bar{\gamma}^2[\sigma^2 + \kappa^2(2K-1)]}{(1-8L^2\Bar{\gamma^2}KD’_1D’_3)[1-12\Bar{\gamma}^2L^2(K+1)(K-2)]}\\
    & \quad + \frac{4MK(K-1)L^2\Bar{\gamma}^2}{1-8L^2\Bar{\gamma^2}KD’_1D’_3} \sum_{n=1}^N \sum_{k=0}^{K-1} \| \frac{\partial F(X_{n;k}) \cdot 1_M}{M} \|_2^2 \\
    & \quad + \frac{4K}{1-8L^2\Bar{\gamma^2}KD’_1D’_3} \sum_{n=1}^N \| \nabla F(\frac{X_n \cdot 1_M}{M})\|_2^2) \\
    & = -\frac{\Bar{\gamma}K}{2}(1 - \frac{8KL\Bar{\gamma}}{(1-8L^2\Bar{\gamma^2}KD’_1D’_3)M} [\frac{\Bar{\gamma}L}{2}[\frac{4K}{(1-\rho)^2} + K - 1] \\
    & \quad \quad + (\frac{2\Bar{\gamma}LK}{1-\rho^2} + \frac{1}{2M}) D_2 ]) \sum_{n=1}^N \mathbb{E} \| \nabla F(\frac{X_n \cdot 1_M}{M})\|_2^2 \\
    & \quad -\frac{\Bar{\gamma}}{2}(1 - 2\Bar{\gamma}LK -  \frac{8K^2(K-1)L^3\Bar{\gamma}^3}{1-8L^2\Bar{\gamma^2}KD’_1D’_3} [\frac{\Bar{\gamma}L}{2}[\frac{4K}{(1-\rho)^2} + K - 1] \\
    & \quad \quad + (\frac{2\Bar{\gamma}LK}{1-\rho^2} + \frac{1}{2M}) D_2 ]) \sum_{n=1}^N \sum_{k=0}^{K-1} \mathbb{E} \| \frac{\partial F(X_{n;k}) \cdot 1_M}{M}\|_2^2 \\
    & \quad + \Bar{\gamma}^2KL (\frac{\Bar{\gamma}L}{2}[\frac{4K}{(1-\rho)^2} + K - 1] + (\frac{2\Bar{\gamma}LK}{1-\rho^2} + \frac{1}{2M} )D_2) \\
    & \quad \quad [(\sigma^2 + 4\kappa^2) + \frac{4L^2\sigma^2\Bar{\gamma}^2(K-1)}{M}\\
    &\qquad + \frac{8(K-1)L^2\Bar{\gamma}^2[\sigma^2 + \kappa^2(2K-1)]}{(1-12\Bar{\gamma}^2L^2(K+1)(K-2))}] \frac{NK}{1-8L^2\Bar{\gamma^2}KD’_1D’_3} 
\end{align*}
where (a) follows Lemma \ref{lemma:3}, (b) follows Lemma \ref{lemma:14} and (c) follows Lemma \ref{lemma:2}.
Suppose $D^\prime_6 > 0$, therefore,
\begin{align*}
    & \quad F_* - F\left(X_1\right) \\
    & \leq -\frac{\Bar{\gamma}K}{2}\left(1-D’_5\right)\sum^N_{n=1}E\|\nabla F\left(\frac{X_n \cdot 1_M}{M}\right)\|^2_2 \\
    & \quad - \frac{\Bar{\gamma}}{2}D_6^\prime\sum^N_{n=1}E\|\frac{\partial F\left(X_n \cdot 1_M\right)}{M}\|^2_2 \\
	& \quad + \frac{\Bar{\gamma}^2KLD’_4N}{1-8L^2\Bar{\gamma}^2KD’_1D’_3}[1+\frac{4L^2\Bar{\gamma}^2\left(K-1\right)}{M} \\
	& \qquad + \frac{8\left(K-1\right)L^2\Bar{\gamma}^2}{1-12\Bar{\gamma}^2L^2\left(K+1\right)\left(K-2\right)}]\sigma^2 \\
	& \quad + \frac{4\Bar{\gamma}^2K^2LD’_4N}{1-8L^2\Bar{\gamma}^2KD’_1D’_3}[1+\frac{4\left(K-1\right)\left(2K-1\right)L^2\Bar{\gamma}^2}{1-12\Bar{\gamma}^2L^3\left(K+1\right)\left(K-2\right)}]\kappa^2
\end{align*}

\subsubsection{Proof of Lemma \ref{lemma:11}}

We firstly find the bound by accumulating all workers from $\{1, ..., M\}$: 
\begin{align*}
    & \quad \sum_{m=1}^M \mathbb{E} \left\| \frac{X_n \cdot 1_M}{M} - x_{n;k}^{\left(m\right)} \right\|_2^2\\
    & = \sum_{m=1}^M \mathbb{E} \| \frac{X_n \cdot 1_M }{M} - \tilde{x}_n^{\left(m\right)} + \Bar{\gamma} \sum_{j=0}^{k-1} \nabla f_m\left(x_{n;j}^{\left(m\right)}, \xi_j^{\left(m\right)}\right) \|_2^2 \\
    & \overset{(a)}{\leq} 2 \sum_{m=1}^M \mathbb{E} \| \frac{X_n \cdot 1_M }{M} - \tilde{x}_n^{\left(m\right)} \|_2^2  + 2 \sum_{m=1}^M \mathbb{E} \| \Bar{\gamma} \sum_{j=0}^{k-1} \nabla f_m\left(x_{n;j}^{\left(m\right)}, \xi_j^{\left(m\right)}\right) \|_2^2 \\
    & \overset{(b)}{\leq} 2 \sum_{m=1}^M \mathbb{E} \| \frac{X_n \cdot 1_M }{M} - \tilde{x}_n^{\left(m\right)} \|_2^2 \\  
    &+ 2 \Bar{\gamma}^2 k \sum_{m=1}^M  \sum_{j=0}^{k-1} \mathbb{E} \| \nabla f_m\left(x_{n;j}^{\left(m\right)}, \xi_j^{\left(m\right)}\right) \|_2^2 \\
    & = 2 \sum_{m=1}^M \mathbb{E} \| \frac{X_n \cdot 1_M }{M} - \tilde{x}_n^{\left(m\right)} \|_2^2  + 2 \Bar{\gamma}^2 k \sum_{j=0}^{k-1} \mathbb{E} \| G\left(x_{n;j}^{\left(m\right)}, \xi_j^{\left(m\right)}\right) \|_F^2
\end{align*}
where (a) follows $\|a+b\|^2_2 \leq 2\|a\|^2_2 + 2\|b\|^2_2$ and (b) is based on Cauthy-Schwarz inequality. Next, we find the bound which sums up all iterations. Therefore, 
\begin{align*}
     & \quad \sum_{n=1}^N \sum_{k=0}^{K-1} \sum_{m=1}^M \mathbb{E} \| \frac{X_n \cdot 1_M}{M} - x_{n;k}^{\left(m\right)} \|_2^2\\
     & \leq 2 \sum_{n=1}^N \sum_{k=0}^{K-1} \sum_{m=1}^M \mathbb{E} \| \frac{X_n \cdot 1_M }{M} - \tilde{x}_n^{\left(m\right)} \|_2^2\\ 
     & \quad + 2 \Bar{\gamma}^2 \sum_{n=1}^N \sum_{k=0}^{K-1} k \sum_{j=0}^{k-1} \mathbb{E} \| G\left(x_{n;j}^{\left(m\right)}, \xi_j^{\left(m\right)}\right) \|_F^2 \\
     & = 2 \sum_{n=1}^N \sum_{k=0}^{K-1} \sum_{m=1}^M \mathbb{E} \| \frac{X_n \cdot 1_M }{M} - \tilde{x}_n^{\left(m\right)} \|_2^2\\ 
     &\quad + 2 \Bar{\gamma}^2 \sum_{n=1}^N \sum_{k=0}^{K-1} \frac{\left(K+k\right)\left(K-k-1\right)}{2} \mathbb{E} \| G\left(x_{n;j}^{\left(m\right)}, \xi_j^{\left(m\right)}\right) \|_F^2 \\
     & \overset{(a)}{\leq} 2 \left(\frac{2K}{1 - \rho^2} \sum_{n=1}^N \mathbb{E} \| C_n \|_F^2 + \frac{2 K^2 \Bar{\gamma}^2}{\left(1-\rho\right)^2} \sum_{n=1}^N \sum_{k=0}^{K-1} \| G\left(x_{n;k}, \xi_k\right) \|_F^2\right) \\
     & \quad + \frac{2\Bar{\gamma}^2 K \left(K-1\right)}{2}  \sum_{n=1}^N \sum_{k=0}^{K-1}  \mathbb{E} \| G\left(x_{n;j}^{\left(m\right)}, \xi_j^{\left(m\right)}\right) \|_F^2 
\end{align*}
\noindent
where (a) comes from Lemma \ref{lemma:10}.

\subsubsection{Proof of Lemma \ref{lemma:14}}

\begin{align*}
    & \quad \mathbb{E} \| \frac{1}{M} \sum_{m=1}^{M} \sum_{k=0}^{K-1} \nabla f_m\left(x_{n;k}^{\left(m\right)}, \xi_k^{\left(m\right)}\right)\|_2^2 \\
    & = \mathbb{E} \| \frac{1}{M} \sum_{m=1}^{M} \sum_{k=0}^{K-1} \left(\nabla f_m\left(x_{n;k}^{\left(m\right)}, \xi_k^{\left(m\right)}\right) - \nabla F_m\left(x_{n;k}^{\left(m\right)}\right)\right)\\
    &\qquad + \frac{1}{M} \sum_{m=1}^M \sum_{k=0}^{K-1} \nabla F_m\left(x_{n;k}^{\left(m\right)}\right) \|_2^2 \\
    & \leq 2 \mathbb{E} \| \frac{1}{M} \sum_{m=1}^{M} \sum_{k=0}^{K-1} \left(\nabla f_m\left(x_{n;k}^{\left(m\right)}, \xi_k^{\left(m\right)}\right) - \nabla F_m\left(x_{n;k}^{\left(m\right)}\right)\right) \|_2^2\\
    & \quad + 2 \mathbb{E} \| \frac{1}{M} \sum_{m=1}^M \sum_{k=0}^{K-1} \nabla F_m\left(x_{n;k}^{\left(m\right)}\right) \|_2^2 \\
    & \leq 2 \frac{1}{M^2} \sum_{m=1}^{M} \sum_{k=0}^{K-1} \mathbb{E} \| \nabla f_m\left(x_{n;k}^{\left(m\right)}, \xi_k^{\left(m\right)}\right) - \nabla F_m\left(x_{n;k}^{\left(m\right)}\right) \|_2^2\\
    & \quad + 2 K \sum_{k=0}^{K-1} \mathbb{E} \| \frac{\partial F\left(x_{n;k}\right) \cdot 1_M}{M} \|_2^2 \\
    & = \frac{2K \sigma^2}{M} + 2K \sum_{k=0}^{K-1} \| \frac{\partial F\left(x_{n;k}\right) \cdot 1_M}{M} \|_2^2
\end{align*}

\subsubsection{Proof of Corollary \ref{corollary:2}}

We selected the stepsize 
\begin{align*}
    \Bar{\gamma}=\left(\sigma\sqrt{N/M} + 3KL^3\sqrt{D_2} + 16KLD’_3 + 6KL\right)^{-1}
\end{align*}
Then, we have:
\begin{align*}
& 1 - 12\Bar{\gamma}^2L^2\left(K+1\right)\left(K-2\right) \geq \frac{2}{3}, \quad D’_1\Bar{\gamma}L \leq \frac{2}{3}, \quad D’_5 \leq \frac{1}{2}  \\
& D’_3KL\Bar{\gamma} \leq \frac{1}{16}, \quad 1 - 8L^2\Bar{\gamma}^2KD’_1D’_3 \geq \frac{2}{3}, \quad D’_6 > 0
\end{align*}
Therefore,
\begin{align*} 
& \quad \frac{1}{2N}\sum^N_{n=1} E\|\nabla F\left(\frac{X_n \cdot 1_M}{M}\right)\|^2_2 \\
& \leq 3\Bar{\gamma}KLD'_4 [1 + \frac{4L^2\Bar{\gamma}^2\left(K-1\right)}{M} + 12\left(K-1\right)L^2\Bar{\gamma}^2]\sigma^2\\
&\quad+ 12\Bar{\gamma}KLD'_4\kappa^2 \left(1 + 6\left(K-1\right)\left(2K-1\right)L^2\Bar{\gamma}^2\right) \\
& \quad + \frac{2\left(F\left(X_1\right)-F_*\right)\left(\sigma\sqrt{N/M} + 3KL^3\sqrt{D_2} + 16KLD'_3 + 6KL\right)}{KN} \\
& = 3\Bar{\gamma}KLD’_4 \left[1 + \frac{4L^2\Bar{\gamma}^2\left(K-1\right)}{M} + 12\left(K-1\right)L^2\Bar{\gamma}^2\right]\sigma^2\\
&\quad+ 4\left[1 + 6\left(K-1\right)\left(2K-1\right)L^2\Bar{\gamma}^2\right]\kappa^2 + \frac{2[F\left(X_1\right)-F_*]\sigma }{K\sqrt{NM}}\\
&\quad+ \frac{2\left(F\left(X_1\right)-F_*\right)\left(\sigma\sqrt{N/M} + 3L^3\sqrt{D_2} + 16LD’_3 + 6L\right)}{N} 
\end{align*}
With the given range of $N$, the following inequality holds: 
\begin{enumerate}
    \item 
    \begin{align*}
        3\Bar{\gamma}KLD’_4 &= 3\Bar{\gamma}KL\left(2\Bar{\gamma}KLD’_3 + \frac{\Bar{\gamma}L\left(K-1\right)}{2} + \frac{D’_2}{2M}\right)\\
        &\leq 6\Bar{\gamma}^2K^2L^2D’_3 + \frac{3\Bar{\gamma}^2L^2K^2}{2} + \frac{3\Bar{\gamma}KLD’_2}{2M} \\
        & = \frac{3K^2L^2M}{\sigma^2N}\left(2D’_3 + \frac{1}{2}\right) + \frac{3KLD’_2}{2\sigma\sqrt{NM}} \leq \frac{3KLD’_2}{\sigma\sqrt{NM}} 
    \end{align*}
    \item 
    \begin{align*}
        & \quad \frac{4L^2\Bar{\gamma}^2\left(K-1\right)}{M} + 12\left(K-1\right)L^2\Bar{\gamma}^2\\
        &= 4\left(K-1\right) L^2\Bar{\gamma}^2\left(\frac{1}{M} + 3\right) \frac{4\left(K-1\right)L^2}{\sigma^2N/M}\left(\frac{1}{M} + 3\right)\\
        & = \frac{4\left(K-1\right)L^2}{\sigma^2N}\left(1 + 3M\right) \leq 1
    \end{align*}
    \item
    \begin{align*}
        6\left(K-1\right)\left(2K-1\right)L^2\Bar{\gamma}^2 \leq \frac{6\left(K-1\right)\left(2K-1\right)L^2}{\sigma^2N/M} \leq 1 
    \end{align*}
\end{enumerate}
Therefore, the following inequality is derived: 
\begin{align*}
&\quad \frac{1}{2N}\sum^N_{n=1} E\|\nabla F\left(\frac{X_n \cdot 1_M}{M}\right)\|^2_2\\
& \leq \frac{6KLD’_2}{\sigma\sqrt{NM}}\left(\sigma^2 + 4\kappa^2\right) + \frac{2[F\left(X_1\right)-F_*]\sigma}{K\sqrt{NM}}\\
& \quad + \frac{2\left(F\left(X_1\right)-F_*\right)\left(3L^3\sqrt{D_2}+16LD’_3+6L\right)}{N}
\end{align*}

\subsection{Communication Cost of QPRSGD}

For the sake of reducing the number of bits in communication, we leverage Elias gamma coding \cite{eliascoding} to compress the vector. In this part, We first introduce several lemmas and a alternative compression scheme is subsequently proposed.  

\begin{lemma} \label{lemma:15}
For any vector $\bm{v} \in \mathbb{R}^d$, the expected number of non-zero values in a vector should be:
\begin{align*}
    E[\| Q_s\left(\bm{v}\right)\|_0] \leq \min \left(s^2 + s\sqrt{d}, d\right) 
\end{align*}

\begin{proof}
Let $u=\bm{v}/\|\bm{v}\|_2$. Let $I\left(u\right)$ be the set of index $i$ where $|u_i| < 1/s$. Since 
\begin{align*}
    \left(d-|I\left(u\right)|\right)/s^2 \leq \sum_{i \notin I\left(u\right)}u_i^2 \leq 1
\end{align*}
\noindent the inequality $d-|I\left(u\right)| \leq s^2$ holds. Then, with the definition \ref{equation:definition_ell}, the probability that $Q_s\left(v_i\right)$ is non-zero value is $|u_i|s$ for all $i \in I\left(u\right)$ and therefore, we have  
\begin{align*}
    E[\|Q_s\left(v\right)\|_0] &= d - |I\left(u\right)| + \sum_{i \in I\left(u\right)} s|u_i|\\ 
    &\leq s^2 + s\|u\|_1 \leq s^2 + s\sqrt{d}
\end{align*}
\noindent Besides, it is easy to notice that $\bm{v}$ is a $d$-dimension vector so that: 
\begin{align*}
    E[\|Q_s\left(v\right)\|_0] \leq d
\end{align*}
\end{proof}
\end{lemma}

\begin{lemma} \label{lemma:16}
Let $\bm{v} \in \mathbb{R}^d$ be a vector so that for all $i$, $v_i$ is a positive integer and moreover, $\|v\|_\rho^\rho \leq \rho$. Then, 
\begin{align*}
\sum_{i=1}^d |Elias\left(v_i\right)| \leq \frac{2d}{\rho}\log\frac{\rho}{d}+d
\end{align*}
\begin{proof}
\begin{align*}
    \sum_{i=1}^d|Elias\left(v_i\right)| & = \sum_{i=1}^d\left(2\log \left(v_i\right) + 1\right) = 2 \sum_{i=1}^d \log \left(v_i\right) + d \\ 
    & \overset{(a)}{\leq} \frac{2}{\rho} d \log \left(\frac{1}{d} \sum_{i=1}^d v_i^\rho\right) + d = \frac{2d}{\rho} \log \frac{\rho}{d} + d
\end{align*}
\noindent where (a) holds on account for Jensen's inequality. 
\end{proof}
\end{lemma}

\subsubsection{Compression Schemes} 

For any integer $k$, we use Elias gamma encoding \cite{eliascoding}, denoted as Elias$\left(k\right)$, to generate its code. The encoding process is simple: let bin$\left(k\right)$ be the binary representation of $k$ and $len$ be the length of bin$\left(k\right)$, the code Elias$\left(k\right)$ would simply be $len$ of zeros added before bin$\left(k\right)$. Therefore, the encoding length |Elias$\left(k\right)$| = O$\left(2log\left(k\right) - 1\right)$ = O$\left(log\left(k\right)\right)$. Such encoding scheme is used to encode positive integer whose upper-bound is unknown, since the actual length of its binary representation could be calculated by the number of 0s before the first 1 received.

For compressed vector $Q_s\left(\textbf{v}\right) = [v'_1, v'_2, \dots, v'_d]$, we have $v'_i = \|\textbf{v}\|_2 \cdot sgn\left(v_i\right) \cdot \zeta\left(v_i, s\right) / s$. We use the following process to implement the encoding: firstly, we put the 32-bit full precision of $\|\textbf{v}\|$ in the beginning of the transmission code. For $i = 1 \dots d$, we use 1 bit to represent $sgn\left(v_i\right)$ and O$\left(log\left(\zeta\left(v_i, s\right)\right)\right)$ bits for Elias$\left(\zeta\left(v_i, s\right)+1\right)$, which are concatenated after the end of previous code in order. The decoding scheme could be processed in the similar way: we first read the 32-bit precision of $\|\textbf{v}\|$, then keep reading $sgn\left(v_i\right)$ and Elias$\left(\zeta\left(v_i, s\right)+1\right)$ until the end of the coding. 

\setcounter{theorem}{2}
\setcounter{corollary}{2}
\begin{theorem}
For any vector $\bm{v} \in \mathbb{R}^d$, in compression scheme 2, the upper bound of the expected communication cost is 
\begin{align*}
    \mathbb{E}[|Code\left(Q_s\left(\bm{v}\right)\right)|] \leq F + 2d + d \log \frac{s^2 + 2s\sqrt{d} + 1 + d/4}{d}
\end{align*}
\begin{proof}
Let $\bm{y} = \left(y_1, y_2, ..., y_d\right)$. Then, we have:  
\begin{align*}
    \mathbb{E}[\|\bm{y}+1\|_2^2] \leq \mathbb{E}[\|\bm{y}\|_2^2] + 2\mathbb{E}[\|\bm{y}\|_1] + 1 \leq \frac{d}{4} + s^2 + 2s\sqrt{d} + 1
\end{align*}
\noindent Therefore, 
\begin{align*}
    \mathbb{E}[|Code\left(Q_s\left(\bm{v}\right)\right)|] &= F + \sum_{i=1}^d \left(1+|Elias\left(y_i+1\right)|\right) \\
    & \leq F + 2d + d \log \frac{s^2 + 2s\sqrt{d} + 1 + d/4}{d}
\end{align*}
\end{proof}
\end{theorem}


\bibliographystyle{IEEEtran}
\bibliography{manuscript}

\begin{IEEEbiography}[{\includegraphics[width=1in,height=1.25in,clip,keepaspectratio]{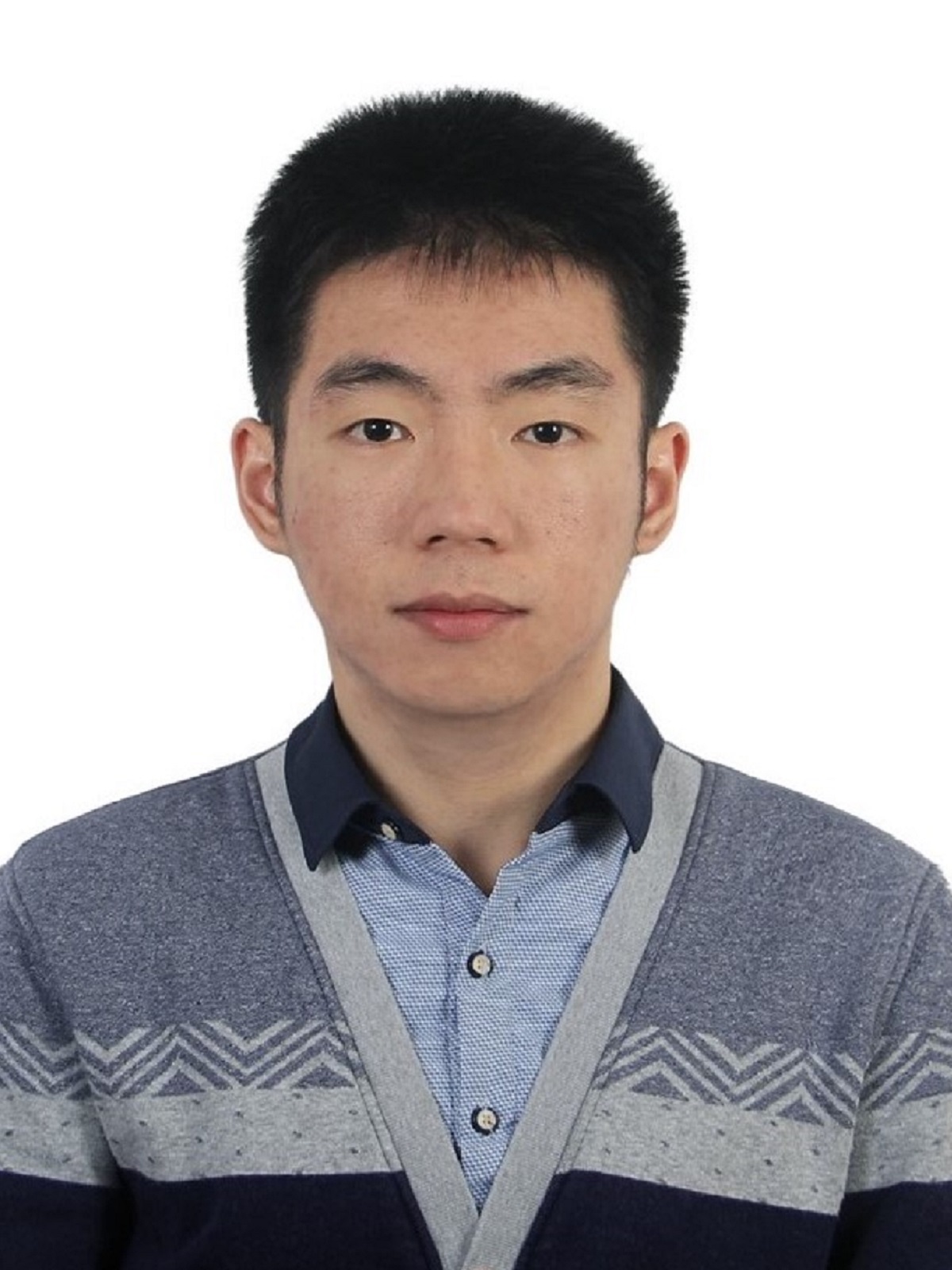}}]{Feijie Wu} is currently an M.Phil. student in the Department of Computing at The Hong Kong Polytechnic University. He obtained his B.Sc. degree in Computing from The Hong Kong Polytechnic University in 2020. During his undergraduate study, he also minored in Applied Mathematics and participated in a half-year exchange program in Technical Computer Science at the University of Twente, Netherlands. Besides, he has been working as a visiting research assistant at The University of British Columbia and The Chinese University of Hong Kong, Shenzhen. His recent research interests include federated learning, blockchain and game. He was a recipient of Best Student Paper from ACM BSCI'19.
\end{IEEEbiography}
\begin{IEEEbiography}[{\includegraphics[width=1in,height=1.25in,clip,keepaspectratio]{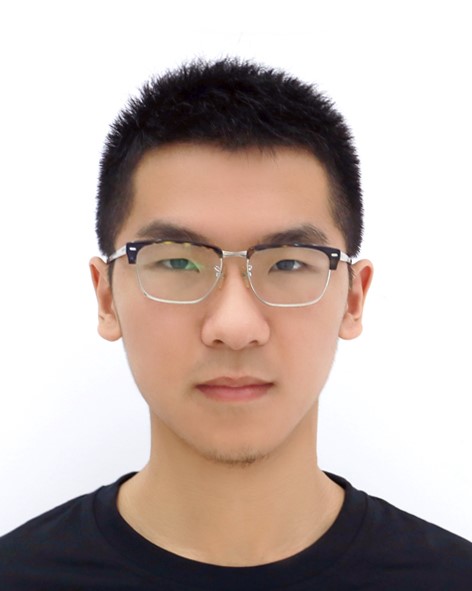}}]{Shiqi He} is currently a master student at the department of computer science, The University of British Columbia, under the supervision of Dr. Ivan Beschastnikh. He received the bachelor's degree from the department of Computing, The Hong Kong Polytechnic University in 2020. His research interests include distributed machine learning and machine learning security.
\end{IEEEbiography}
\begin{IEEEbiography}[{\includegraphics[width=1in,height=1.25in,clip,keepaspectratio]{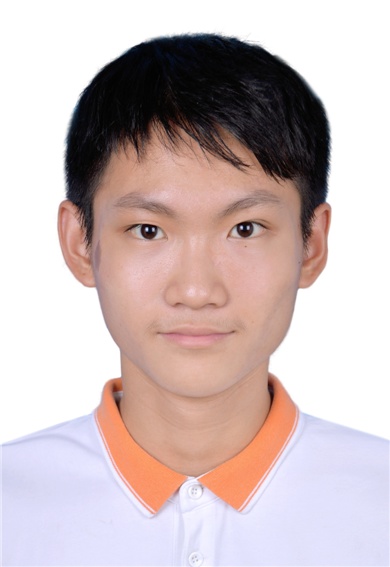}}]{Yutong Yang} is currently taking his Master Programme in National University of Singapore, and obtained Bachelor degree in Computer Science from the Hong Kong Polytechnic University. Besides, he serves as a student assistant at the Hong Kong Polytechnic University. His research interest mainly lies in distributed learning.
\end{IEEEbiography}
\begin{IEEEbiography}[{\includegraphics[width=1in,height=1.25in,clip,keepaspectratio]{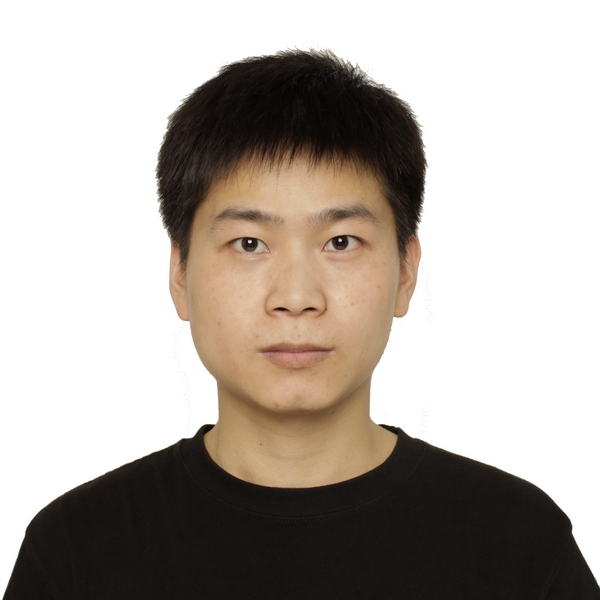}}]{Haozhao Wang} is currently a Ph.D. candidate in the School of Computer Science and Technology at Huazhong University of Science and Technology and a research assistant in the Department of Computing at The Hong Kong Polytechnic University. His research interests include Distributed Machine Learning and Federated Learning.
\end{IEEEbiography}
\begin{IEEEbiography}[{\includegraphics[width=1in,height=1.25in,clip,keepaspectratio]{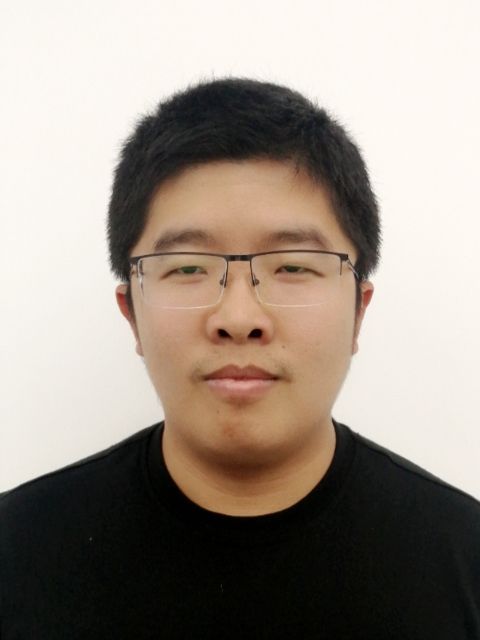}}]{Zhihao Qu} received his B.S. and Ph.D. degree in computer science from Nanjing University, Nanjing, China, in 2009, and 2018, respectively. He is currently an assistant researcher in the College of Computer and Information at Hohai University and in the Department of Computing at The Hong Kong Polytechnic University. His research interests are mainly in the areas of wireless networks, edge computing, and distributed machine learning.
\end{IEEEbiography}
\begin{IEEEbiography}[{\includegraphics[width=1in,height=1.25in,clip,keepaspectratio]{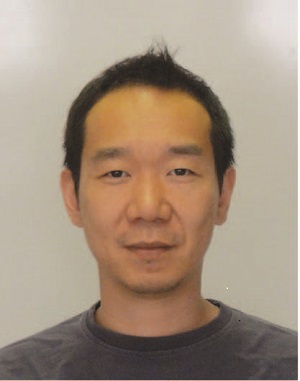}}]{Song Guo} (M’02-SM’11-F’19) received the Ph.D. degree in computer science from the University of Ottawa and was a professor with the University of Aizu. He is a full professor with the Department of Computing, The Hong Kong Polytechnic University. His research interests are mainly in the areas of big data, cloud computing, mobile computing, and distributed systems. He is the recipient of the 2019 IEEE TCBD Best Conference Paper Award, 2018 IEEE TCGCC Best Magazine Paper Award, 2019 \& 2017 IEEE Systems Journal Annual Best Paper Award, and other 6 Best Paper Awards from IEEE/ACM conferences. His work was also recognized by the 2016 Annual Best of Computing: Notable Books and Articles in Computing in ACM Computing Reviews. He is an IEEE Fellow (Computer Society) and the Editor-in-Chief of IEEE Open Journal of the Computer Society. He was a Distinguished Lecturer of IEEE Communications Society (ComSoc) and served in the IEEE ComSoc Board of Governors. He has been named on editorial board of a number of prestigious international journals like IEEE Transactions on Parallel and Distributed Systems, IEEE Transactions on Cloud Computing, IEEE Transactions on Emerging Topics in Computing, etc. He has also served as chairs of organizing and technical committees of many international conferences.
\end{IEEEbiography}
\begin{IEEEbiography}[{\includegraphics[width=1in,height=1.25in,clip,keepaspectratio]{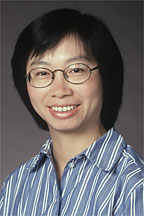}}]{Weihua Zhuang}
    (M'93--SM'01--F'08) has been with the Department of Electrical and Computer Engineering, University of Waterloo, Canada, since 1993, where she is a Professor and a Tier I Canada Research Chair in Wireless Communication Networks. Her current research focuses on resource allocation and QoS provisioning in wireless networks, and on smart grid. She is a co-recipient of several best paper awards from IEEE conferences. Dr. Zhuang was the Editor-in-Chief of IEEE Transactions on Vehicular Technology (2007-2013), and the Technical Program Chair/Co-Chair of the IEEE VTC Fall 2017/2016. She is a Fellow of the IEEE, a Fellow of the Canadian Academy of Engineering, a Fellow of the Engineering Institute of Canada, and an elected member in the Board of Governors and VP Publications of the IEEE Vehicular Technology Society.
    \end{IEEEbiography}    
\end{document}